\documentclass[ngerman,twoside,numbers=noenddot]{scrartcl}
\KOMAoptions{DIV=12, BCOR=0mm, paper=a4, fontsize=12pt, parskip=half, twoside=true, titlepage=true}
\usepackage[onehalfspacing]{setspace} 
\usepackage{authblk}
\usepackage[headsepline,automark,pagestyleset=KOMA-Script,markcase=ignoreuppercase]{scrlayer-scrpage}
\clearmainofpairofpagestyles
\ohead{\pagemark}
\ihead{\headmark}
\usepackage[left=3cm,right=3cm,top=1.5cm,bottom=3cm,head=35.9889pt,includehead]{geometry}
\usepackage[utf8]{inputenc} 
\usepackage[english]{babel} 
\usepackage[expansion=true, protrusion=true]{microtype}
\usepackage{amsmath}
\numberwithin{equation}{section}
\usepackage{calc}
\usepackage{mathtools}
\usepackage{amsfonts}
\usepackage{amssymb} 
\usepackage{amsthm}
\usepackage{upgreek}
\usepackage{bbold}
\usepackage{enumitem}
\usepackage[dvipsnames]{xcolor}
\usepackage{units}
\usepackage{graphicx} 
\usepackage[hypcap]{caption} 
\usepackage{booktabs} 
\usepackage{flafter} 
\usepackage[section]{placeins}
\usepackage{tikz-cd}
\usepackage{slashed}
\usepackage{hyperref} 
\hypersetup{colorlinks=true, breaklinks=true, citecolor=red, linkcolor=blue, menucolor=red, urlcolor=cyan, bookmarksopen=false, bookmarksopenlevel=0, plainpages=false,hypertexnames=false}
\newtheorem{theorem}{Theorem}[section]
\newtheorem{proposition}[theorem]{Proposition}
\newtheorem{lemma}[theorem]{Lemma}
\newtheorem{corollary}[theorem]{Corollary}

\newtheorem{remark}[theorem]{Remark}

\theoremstyle{definition}
\newtheorem{definition}[theorem]{Definition}
\newtheorem{example}[theorem]{Example}

\usepackage[T1]{fontenc}
\usepackage{lmodern}
\usepackage{pdfpages}

\usepackage[autostyle]{csquotes}
\usepackage[numbers]{natbib}

\newcommand{\bigslant}[2]{{\left.\raisebox{.2em}{$#1$}\middle/\raisebox{-.2em}{$#2$}\right.}}
\newcommand\restr[2]{{\left.#1\right|_{#2} }}

\newcommand{\bismut}{\nabla^\pm}
\newcommand{\bismutthird}{\nabla^{\pm 1/3}}

\newcommand{\divergence}{\mathrm{div}}
\newcommand{\grad}{\operatorname{grad}}
\newcommand{\extd}[0]{\mathrm{d}}
\newcommand{\id}[0]{\mathrm{id}}

\newcommand{\idon}[1]{\id_{#1}}
\newcommand{\riem}{\mathrm{Rm}}
\newcommand{\ric}{\mathrm{Rc}}

\newcommand{\rscal}{\mathrm{Sc}}
\newcommand{\vol}{\mathrm{vol}}
\newcommand{\genmet}{\mathcal{G}}
\newcommand{\genmetHS}{\mathcal{H}}

\newcommand{\genscal}{\mathcal{S}c}
\newcommand{\genfullric}{\overline{\mathcal{R}c}}
\newcommand{\genric}{\mathcal{R}c}
\newcommand{\genriem}{\mathcal{R}m}

\newcommand{\genext}{\mathcal{K}}
\newcommand{\genextEnd}{\mathcal{A}}

\newcommand{\gennormalext}{\mathcal{L}}
\newcommand{\gennormalshape}{\mathcal{B}}
\newcommand{\genmean}{\mathcal{T}}

\newcommand{\scalbrack}[1]{\left\langle #1 \right\rangle}
\newcommand{\scalprod}[2]{\scalbrack{#1, #2}}
\newcommand{\scalprodmap}{\scalprod{\cdot}{\cdot}}
\newcommand{\brackmap}{[\cdot,\cdot]}

\newcommand{\End}{\mathrm{End}}

\newcommand{\Sym}{\mathrm{Sym}}

\newcommand{\reals}{\mathbb{R}}
\newcommand{\absolute}[1]{\left\lvert #1 \right\rvert}

\newcommand{\tr}{\mathrm{tr}\:}
\newcommand{\trwith}[1]{\mathrm{tr}_{#1}\:}

\newcommand{\spanspace}{\mathrm{span}\:}

\newcommand{\sym}{\mathrm{sym}}
\newcommand{\antisym}{\mathrm{antisym}}
\newcommand{\sigmasym}{{\sigma\text{-}\mathrm{sym}}}
\newcommand{\sigmaantisym}{{\sigma\text{-}\mathrm{antisym}}}

\begin{document}
\pagenumbering{Roman}
\thispagestyle{empty}
\cleardoublepage
\setcounter{page}{1}
\thispagestyle{empty}

\cleardoublepage
\pagenumbering{arabic}

\title{Exterior Generalised Geometry}

\date{\today}

\author{Vicente Cort\'es and Oskar Schiller}

\maketitle

\tableofcontents
\section{Abstract}

It is the aim of this paper to transfer to generalised geometry tools employed in the study of semi-Riemannian immersions, specializing at times to semi-Riemannian hypersurfaces. 

Given an exact Courant algebroid $E \to M$ and an immersion $\iota\colon N \hookrightarrow M$, there is a well-known construction of an exact Courant algebroid $\iota^! E \to N$, the pullback of $E$. This paper explains the pullback of generalised metrics and divergence operators. Assuming $N$ is a hypersurface, it develops the notion of generalised exterior curvature, introducing the generalised second fundamental form and the generalised mean curvature. Generalised versions of the Gauß-Codazzi equations are obtained. As an application, the constraint equations for the initial value formulation of the generalised Einstein equations are established in the formalism of generalised geometry. Further applications include a generalised geometry version of the fundamental theorem for hypersurfaces and the result that generalised K\"ahler and hyper-K\"ahler structures restrict to submanifolds compatible with the generalised almost complex structure. In particular, we characterise exact semi-Riemannian Courant algebroids which are flat with respect to the canonical generalised connection. These play the role of the ambient space in the fundamental theorem mentioned above.

\textit{MSc classification}: 53D18 (Generalized geometry a la Hitchin); 35Q76 (Einstein's equations).\

\textit{Key words:} generalised metrics, submanifold theory, (generalised) Einstein equations, energy and momentum constraints in supergravity

\section{Introduction}

Submanifold geometry is foundational to modern semi-Riemannian geometry. It is the goal of this work to provide the corresponding language in generalised geometry, so that questions asked in the classical setting may be discussed in the new setting of generalised geometry. 

Potential applications in generalised geometry are numerous. Our main motivation is the formulation and study of the generalised Einstein equations as an initial value problem, which is the main theme of the advanced PhD project of the second author. As a preliminary step and first application of our theory, we determine here the constraint equations for the above initial value problem. More specifically, analogues of the classical Hamiltonian and momentum constraints in general relativity are  established in Corollaries  \ref{genenergyconstrcor} and \ref{genmomentumconstrcor}. By decomposing these equations in terms of classical data and assuming the divergence operator to be associated with a (dilaton) function we finally recover the formulas \cite[Proposition 5.11.]{CarlosPhDthesis}, which were obtained using the language of gerbes for a description of NS-NS supergravity. In addition the constraint equations for the Killing spinor equations in those supergravity theories were obtained in \cite{CarlosPhDthesis}. Working in the framework of geometric structures on Courant algebroids, we expect that our approach will provide additional insight in the structure of the moduli spaces of initial data. 

Restricting attention from generalised Ricci flat to flat generalised ambient metrics we arrive at a realisation problem analogous to the fundamental problem of hypersurface theory, which we solve in Theorem \ref{fth}. Other geometric structures and corresponding evolution equations on Courant algebroids can be studied with the extrinsic methods developed here. That is left for future investigation. 

Our theory is well adapted to the study of geometric structures subject to integrability conditions characterised by the existence of compatible generalised connections, such as the ones studied in \cite{cortesdavid}. As an example, we show in Corollary~\ref{GK_preserved:cor} that a generalised K\"ahler structure on an  exact Courant algebroid $E$ over $M$ induces a generalised K\"ahler structure over any semi-Riemannian submanifold $N$ for which the induced Courant algebroid $E_N$ is invariant under the generalised almost complex structure. A similar result holds for generalised hyper-K\"ahler structures.

We begin in section \ref{pullbackCAsection} by reviewing the pullback Courant algebroid, which was introduced in \cite[Lemma 3.7]{BursztynCavalcantiGualtieri2007} and \cite[Proposition 7.1.]{gualtieriPhDannals}. 

Given an exact Courant algebroid (CA) $E$ over a manifold $M$, and an immersion $\iota\colon N \hookrightarrow M$, the \emph{pullback CA} $\iota^! E$ is a canonical exact CA over $N$. Denoting by $K = \mathrm{Ann}(TN) \subset \iota^*(T^*M)$ the annihilator of $TN$ and by $K^\perp \subset \iota^*E$ the orthogonal bundle, the pullback CA is defined as the quotient vector bundle
\begin{equation*}
    \iota^! E \coloneqq \bigslant{K^\perp}{K}
\end{equation*}
equipped with a canonical anchor, bracket, and inner product. To obtain a canonical realisation of this quotient bundle, we introduce the notion of a \emph{Courant transversal bundle}, which is a subbundle $\mathcal{T} \subset \iota^* E$ such that the sequence
\begin{equation*}
    0 \longrightarrow K \longrightarrow \mathcal{T} \longrightarrow \nu \longrightarrow 0
\end{equation*}
is exact for some subbundle $\nu \subset \iota^*TM$ transversal to $TN$. It is immediate that the orthogonal bundle $E_N \coloneqq \mathcal{T}^\perp \subset \iota^*E$ with respect to the inner product is a concrete realisation of the quotient bundle $K^\perp / K$. Indeed, we find in Proposition \ref{inducedCAprop} that one can canonically equip $E_N$ with the structure of an exact CA, so that $\iota^! E \cong E_N$ on the level of exact CAs. 

If $E$ is equipped with a generalised metric $\genmet$, there are two reasons for working with the concrete realisation of the pullback bundle $\iota^!E$ as the orthogonal complement $E_N = \mathcal{T}^\perp$ of a Courant transversal bundle. One, a generalised metric provides (at least if the immersion is semi-Riemannian) a canonical choice of Courant transversal bundle: the Courant normal bundle $\mathcal{N}$. Two, with this canonical choice $E_N = \mathcal{N}^\perp$ is the orthogonal complement of $\mathcal{T}$ with respect to both the inner product and the generalised metric. As a consequence, definitions as well as calculations simplify. 

In section \ref{inheritedgeometrysection}, we consider the ambient CA $E$ to be equipped with a semi-Riemannian generalised metric $\genmet$, and consider the canonical realisation of the pullback CA as $E_N = \mathcal{N}^\perp$. We show that the restriction
\begin{equation*}
    \genmetHS \coloneqq \restr{\genmet}{E_N \times E_N}
\end{equation*}
defines a generalised metric $\genmetHS$ on the $E_N$. Furthermore, if $\divergence = \divergence^\genmet - \scalbrack{e, \cdot}$ is a divergence operator on $E$, then 
\begin{equation*}
    \divergence_N \coloneqq \divergence^\genmetHS - \scalbrack{e^\parallel, \cdot}
\end{equation*}
with $e^\parallel \coloneqq \pi^\parallel\iota^* e \in \Gamma(E_N)$ defines a canonical divergence operator on $E_N$. Here $\pi^\parallel : \iota^*E = E_N \oplus \mathcal N\to E_N$ denotes the canonical projection. Finally, if $D$ is a generalised connection $E$, then setting for arbitrary $u,v \in \Gamma_{\mathrm{loc}}(E_N)$ and extensions $\Tilde{u},\Tilde{v} \in\Gamma_{\mathrm{loc}}(E)$
\begin{equation}\label{inheritedconndefeq}
    D^N_u v \coloneqq \pi^\parallel(D_{\Tilde{u}}\Tilde{v})
\end{equation}
one obtains a generalised connection on $E_N$. The following Proposition summarizes results obtained in Lemma \ref{gausslemma1} and Proposition \ref{divdiscompatibilityprop}.
\begin{proposition}\label{genconninheritanceprop}
    If $D$ is torsion-free, then so is $D^N$. If $D$ is $\genmet$-metric, then $D^N$ is $\genmetHS$-metric. However, compatibility of $D$ with $\divergence$ does not generally imply compatibility of $D$ with $\divergence_N$.
\end{proposition}
Given a pair $(\genmet, \divergence)$ consisting of a generalised metric $\genmet$ and a divergence operator $\divergence$, consider a generalised LC connection $D \in \mathcal{D}^0(\genmet,\divergence)$. There are two ways in which the triple $(\genmet,\divergence ,D)$ induces canonically a generalised connection on $N$. Directly, $D$ induces a generalised connection $D^N$ by restriction and projection as defined in (\ref{inheritedconndefeq}). Indirectly, the pair $(\genmet, \divergence)$ induces the pair $(\genmetHS, \divergence_N)$ and then the canonical generalised LC connection $\Tilde{D}^N \in \mathcal{D}^0(\genmetHS, \divergence_N)$. Crucially, by Proposition \ref{genconninheritanceprop}, the two definitions are inequivalent; In general, $D^N \neq \Tilde{D}^N$. But which definition is to be preferred?

We find that many calculations simplify if one works with $D^N$, and hence mostly work with $D^N$. However, in some instances, $\Tilde{D}^N$ contains the relevant information. An example for this is the formulation of the generalised Einstein equations as an initial value problem. 

The generalised Einstein equations are equations precisely for a pair $(\genmet,\divergence)$ consisting of a generalised metric and a divergence operator. A generalised LC connection is only introduced as a tool, but it is not the fundamental object of consideration. Hence, considering a hypersurface $\Sigma \hookrightarrow M$ and the canonical realisation of the pullback CA $E_\Sigma$, one is interested in which generalised metric $\genmetHS$ and divergence operator $\divergence_\Sigma$ it inherits. As an auxiliary object, one may then consider the canonical generalised LC connection $\Tilde{D}^\Sigma$.

In such a scenario, where the generalised metric and the divergence operator are the objects of study, we find that it is useful to introduce the space $\mathcal{D}^0(\genmet, \divergence; N)$ of generalised LC connections that project to generalised LC connections with divergence $\divergence_N$ on the immersed manifold $N$. We prove in Lemma \ref{HScanconnlemma} that this set is non-empty.

In section \ref{genextgeomsection}, we introduce the \emph{generalised second fundamental form}, an object that captures exterior curvature in generalised geometry. Let $(E,\genmet)$ be a semi-Riemannian exact CA with base manifold $M$. Denote by $g$ the semi-Riemannian metric on $M$ induced by $\genmet$. Consider $\Sigma \hookrightarrow (M,g)$ a semi-Riemannian hypersurface, denote by $n$ a unit normal on $\Sigma$, and by $n_\pm$ the lifts of $n$  to the $\pm1$-eigenbundles $E_\pm$  of $\mathcal G^{\mathrm{End}}$. Then the generalised second fundamental form $\genext^{n_\pm} \in \Gamma(E_\Sigma^* \otimes E_\Sigma^*)$ and the the generalised shape tensor $\genextEnd^{n_\pm} \coloneqq Dn_\pm \in \Gamma(E_\Sigma^* \otimes E_\Sigma)$ are defined as
\begin{equation*}
    \genext^{n_\pm}(a,b) \coloneqq \genmet(D_a n_\pm, b), \qquad\qquad \genextEnd^{n_\pm}(a) \coloneqq D_a n_\pm.
\end{equation*}
The generalised second fundamental form does not capture all exterior curvature, we are missing the components of the \emph{conormal exterior curvature} $\gennormalext^\pm \in \Gamma(E_\Sigma^*)$, defined as
\begin{equation*}
    \gennormalext^\pm(a) \coloneqq \genmet(D_{n_\pm - n_\mp}n_\pm, a)
\end{equation*}
The conormal exterior curvature vanishes for a generalised LC connection if and only if $D \in \mathcal{D}^0(\genmet, \divergence; \Sigma)$. 

The main result characterizing the generalised second fundamental form is Lemma~\ref{genextcomputationlemma}, which states the following.
\begin{lemma}
    Let $\divergence = \divergence^\genmet - \scalbrack{e,\cdot}$ be a divergence operator, and let $D^\chi = D +\chi \in \mathcal{D}^0(\genmet,\divergence)$ be an arbitrary generalised LC connection, with $D$ the canonical choice and $\chi = \chi_+ + \chi_- \in \Sigma_0^+ \oplus \Sigma_0^-$, where
    \[ \Sigma_0^\pm := (\mathfrak{so}(E_\pm)^{\langle 1 \rangle })_0 := 
    \{ \chi \in \mathfrak{so}(E_\pm)^{\langle 1 \rangle } \mid \mathrm{tr}(\chi v) =0\text{ for all }v\in E_\pm\}. \]
    Denote $\chi_\pm^\perp(a_\pm,b_\pm) \coloneqq \chi_\pm(a_\pm,b_\pm,n_\pm)$. Then
    \begin{equation}
        \begin{split}
            \restr{\mathcal{K}^{n_\pm}}{E_\Sigma^\pm \times E_\Sigma^\pm} &= k - \chi_\pm^\perp - \frac{\scalbrack{e,n_\pm}}{\dim \Sigma} h \mp \frac{i_n H}{6} \\
            \restr{\mathcal{K}^{n_\pm}}{E_\Sigma^\mp \times E_\Sigma^\pm} &= k \mp \frac{i_n H}{2} \\
        \end{split}
    \end{equation}
    Herein, $h$ is the semi-Riemannian metric on $\Sigma$, $k$ the second fundamental form, $H$  the closed $3$-form associated with the pair $(E,\mathcal G)$, $E_\Sigma^\pm := E_\Sigma \cap E_\pm$  and we employed the canonical isometries $\sigma_\pm \colon (TM,g) \cong (E_\pm,\genmet)$ provided by the generalised metric $\genmet$ to identify $E_\pm$ with $TM$. 
\end{lemma}
In particular, the mixed-type part $\genext^\pm \coloneqq \restr{\genext^{n_\pm}}{E_\mp\times E_\pm}$ of the generalised second fundamental form is completely determined by the generalised metric $\genmet$ and the \emph{generalised mean curvature} $\genmean^\pm \coloneqq  \trwith{\genmet} \genext^{n_\pm} =  \trwith{h}k-\scalbrack{e,n_\pm}$ is entirely determined by the pair $(\genmet, \divergence)$ consisting of the generalised metric and the divergence operator. 

In section \ref{gengausseqsecion}, we compute the generalised Gauß and Codazzi equations. To that end, we take a generalised LC connection $D \in \mathcal{D}^0(\genmet,\divergence)$, and denote by $D^\Sigma \in \mathcal{D}^0(\genmetHS)$ the generalised connection obtained from restriction and projection, cf. (\ref{inheritedconndefeq}). So, in general, $\divergence_{D^\Sigma} \neq \divergence_\Sigma$. We obtain in Lemma \ref{finalgausslemma} the following 
\emph{generalised Gauß and Weingarten formulas}.
\begin{lemma}
    Let $a,b \in \Gamma(E_\Sigma)$, denote by $a_\pm,b_\pm$ their projections to $\Gamma(E_{\Sigma\pm})$, and let $f \in C^\infty(M)$. Then
    \begin{equation}
        \begin{split}
            (D_a b_\pm)^\parallel &= D^\Sigma_a b_\pm\\
            \genmet(D_a b_\pm,n_\pm) &= - \genext^{n_\pm}(a, b_\pm) \\
            D_a (f n_\pm) &= \pi(a)(f) n_\pm + f \genextEnd^{n_\pm}(a)
        \end{split}
    \end{equation}
    Herein, $(\cdot)^\parallel \colon \iota^*E \to E_\Sigma$ denotes the orthogonal projection.
\end{lemma}
The main results computed with this formula are the \emph{generalised Gauß and Codazzi equations}. The former are presented in Theorem \ref{riemgaußthm}, which we repeat here.  
\begin{theorem}
    The pure-type part of the generalised Riemann tensors obtained from $D$ and $D^\Sigma$ satisfies  the following Gauß equation, where $a,b,v,w \in \Gamma(E_{\Sigma}^\pm)$:
    \begin{equation*}
        \begin{split}
            &\pm 2 \varepsilon\left\{\genriem^D(a,b,v,w) - \genriem^{D^\Sigma}(a,b,v,w)\right\} \\
            &= \genext^{n_\pm}(w, a) \genext^{n_\pm}(v, b)  - \genext^{n_\pm}(v, a) \genext^{n_\pm}(w, b) \\
            & \quad+ \genext^{n_\pm}(a, w) \genext^{n_\pm}(b, v) - \genext^{n_\pm}(b, w) \genext^{n_\pm}(a, v) \\
            &\quad  +[\genext^{n_\pm}(v,w) - \genext^{n_\pm}(w,v)] [\genext^{n_\pm}(b, a) - \genext^{n_\pm}(a, b)] \\
        \end{split}
    \end{equation*}
    The mixed-type part satisfies with $a,v,w \in \Gamma(E_{\Sigma}^\pm), \bar{b} \in \Gamma(E_{\Sigma}^\mp)$
    \begin{equation*}
        \begin{split}
            &\pm2\varepsilon \left\{ \genriem^D(a, \bar{b}, v, w) - \genriem^{D^\Sigma}(a, \bar{b}, v, w)\right\} \\
            &=  \genext^{n_\pm}(a,w)\genext^{\pm}(\bar{b}, v)  -\genext^\pm(\bar{b}, w) \genext^{n_\pm}(a,v) \\
            &\quad +  [\genext^{n_\pm}(v, w)-\genext^{n_\pm}(w,v)] \genext^\pm(\bar{b}, a) 
        \end{split}
    \end{equation*}
\end{theorem}
The generalised Codazzi equations are derived in Theorem \ref{codazzithm}. We state them here.   
\begin{theorem}
    The normal components of the generalised Riemann curvature satisfy the Codazzi equations
    \begin{equation*}
        \pm2 \genriem^D(a, \bar{b}, n_\pm,w) = [D^\Sigma_{\bar{b}}\genext^{n_\pm}](a,w) - [D^\Sigma_a\genext^\pm](\bar{b},w)+ \varepsilon \genext^\pm(\bar{b},a) \gennormalext^\pm(w)
    \end{equation*}
    and
    \begin{equation*}
        \begin{split}
            &\pm2 \genriem^D(a, b, n_\pm- n_\mp,w) \\
            & = [D^\Sigma_w \genext^{n_\pm}](a,b) - [D^\Sigma_w\genext^{n_\pm}](b,a) + [D^\Sigma_{b}\genext^{n_\pm}](a,w) - [D^\Sigma_a\genext^{n_\pm}](b,w)\\
            &\quad + \varepsilon\:\big\{\gennormalext^\pm(w)[\genext^{n_\pm}(b,a)- \genext^{n_\pm}(a,b) ] +\gennormalext^\pm(b)\genext^{n_\pm}(w,a)-\gennormalext^\pm(a) \genext^{n_\pm}(w,b) \big\}
        \end{split}
    \end{equation*}
    and
    \begin{equation*}
        \begin{split}
            & \pm 2 \genriem^D(a, n_\pm - n_\mp, n_\pm - n_\mp, b) \\
            &= 2 [(\genext^{n_\pm})^2]^{\mathrm{sym}}(a,b) - 2 \genext^{n_\pm}(a, \genextEnd^{n_\pm} b) - \genext^\pm(\genextEnd^\mp a, b) \\
            &\quad + \trwith{\genmetHS} \left[\genext^{n_\pm}(\cdot,a) \genext^{n_\pm}(\cdot, b)\right]_{E_\pm \times E_\pm} - [D^\Sigma_a \gennormalext](b) - [D^\Sigma_b \gennormalext](a) 
        \end{split}
    \end{equation*}
    and
    \begin{equation*}
        \begin{split}
            &\pm 2 \genriem^D(a, n_\pm - n_\mp, n_\pm - n_\mp, \bar{b}) \\
            &= 4 \{(\genext^{n_\sigma})^2\}^{\sigmaantisym}(\bar{b},a) - 2 \{\genext^{n_\sigma}(\cdot,\genextEnd^\sigma \cdot) \}^\sigmaantisym +  [D_a^\Sigma \gennormalext^\mp](\bar{b}) - [D_{\bar{b}}^\Sigma \gennormalext^\pm](a)
        \end{split}
    \end{equation*}
    Herein, $a,b,w \in \Gamma(E_{\Sigma}^\pm)$ and $\bar{b} \in \Gamma(E_{\Sigma}^\mp)$.
\end{theorem}
From these equations, several Corollaries for important traces are derived. In particular, we find that the constraint equations for the generalised Einstein equations can be nicely  stated as follows, cf. Corollaries \ref{genenergyconstrcor} and \ref{genmomentumconstrcor}.
\begin{equation*}
    \begin{split}
        2 \genric^\pm(n_\mp, n_\pm) - \varepsilon \genscal &= -\varepsilon \genscal_\Sigma  - \absolute{\genext^\pm}^2 + \frac{(\genmean^+)^2+(\genmean^-)^2}{2} \\
        \genric^\pm(a_\mp, n_\pm) &= (\divergence_\Sigma \genextEnd^\pm)(a_\mp) - \pi a_\mp(\genmean^\pm) \\
    \end{split}
\end{equation*}
Herein, $\varepsilon = g(n,n)$.

In section \ref{flatCAsection}, we characterise the semi-Riemannian Courant algebroids which are flat with respect to the canonical generalised connection. We prove that in general, the metric induced on the base manifold is conformally flat, see Theorem \ref{genriemflatnessthm}. Inspecting the Riemannian and the Lorentzian case more closely in Corollaries \ref{riemannian_genriemflatnesscor} and \ref{lorentzian_genriemflatness_cor}, we find that in these signatures a flat exact CA is completely trivial: it is untwisted, has metric divergence (or trivial dilaton) and a flat Riemannian base manifold. In neutral signature, however, we establish the existence of a non-trivial flat exact CA, see Example \ref{genriemflatness_example}. Note that our result is more restrictive than \cite{cavalcanti-pedregal-rubio} which is based on a weaker notion of flatness.

 In section \ref{fundthmsection}, we employ the characterisation of flat exact CAs to bring a classical problem from Riemannian geometry to generalised geometry, as we establish the \enquote{generalised} version of the fundamental theorem for hypersurfaces. The classical result concerns data $(\Sigma,h,k)$ consisting of a Riemannian manifold $(\Sigma,h)$ and a symmetric two-tensor $k$ on $\Sigma$. It states that $(\Sigma,h)$ can be locally realised as a Riemannian hypersurface in Euclidean space such that $k$ corresponds to the second fundamental form if and only if the data $(\Sigma,h,k)$ satisfies the flat Gauß and Codazzi equations. We provide with Theorem~\ref{fth} an analogous result in the setting of exact Riemannian Courant algebroids.

As two appendices, we present in section \ref{genriemappendix} the expression for the generalised Riemann tensor of the canonical generalised LC connection with specified divergence, and in section \ref{mixedtypedivappendix}, we develop the formalism for taking the divergence of mixed-type tensors. In particular, this enables us to express the constraint equations for the generalised Einstein equations, derived in section \ref{gengausseqsecion} in the language of generalised geometry, in terms of classical geometrical objects.

\textbf{Acknowledgements.}   Research of VC is funded by the Deutsche\linebreak Forschungsgemeinschaft 
(DFG, German Research Foundation) under Germany's Excellence Strategy, EXC 2121 ``Quantum Universe,'' 390833306 and under -- SFB-Gesch\"aftszeichen 1624 -- Projektnummer 506632645.

As indicated in section \ref{genriemappendix} (respectively section \ref{fundthmsection}), the expression for the generalised Riemann tensor of the canonical generalised Levi-Civita connection is the result of work in progress \cite{vicentethomas,vicenteliana}.  We are very thankful to Liana David, Matas Mackevicius, and Thomas Mohaupt for their support of this project and for further discussions. We are indebted to Carlos Shahbazi for sharing ideas related to this project. We also thank him and Miguel Pino Carmona for support in calculating the constraint equations for the generalised Einstein equations. Finally, we thank Liana David for many helpful comments on a draft of this paper.  

\subsection{Notation and Conventions}
This section is intended to provide an overview of notation and conventions employed throughout the text.

We assume smoothness throughout. We mostly work on an exact Courant algebroid $\pi \colon E \to M$ over a manifold $M$, which has another manifold $N$ immersed. If $N$ is assumed to be of co-dimension one, we denote it as $\Sigma$. 

Generally, we denote vector fields on manifolds by capital Latin letters such as $V,W,X,Y,...$, the exception being unit normal vector fields on an immersed manifold $N$ or $\Sigma$, which we denote by $n$ to avoid confusion. Given a generalised metric on $E$, we denote the associated sections in the $\pm 1$ eigenbundle of $E$ by $n_\pm$, i.e.\ $\pi n_\pm = n$. We generally denote sections of an exact Courant algebroid, i.e. generalised vector fields, by lower case latin letters $a,b,v,w,...$, and try to match letters up if the sections correspond, so that we have e.g.\ $V = \pi v$.

Working with a representative $H \in \Omega^3(M)_{\mathrm{cl}}$ of the Ševera class of $E$, and  given a metric $g$ on $M$, we make frequent use of the following two contractions of $H^{\otimes 2}=H \otimes H$ with $g$: $H^{(2)}$, the $(0,4)$ tensor such that $H^{(2)}(X,Y,V,W) = \trwith{g} H(X,Y,\cdot)H(V,W,\cdot)$, and $H^2$, the $(0,2)$ tensor such that $H^2(X,Y) = \trwith{g} H^{(2)}(X, \cdot, Y, \cdot)$. Up to sign, these are the only possible non-trivial contractions of type $(0,4)$ and $(0,2)$, respectively. Note that, with $H_X$ the
$g$-skew-symmetric endomorphism such that $g(H_X Y, Z) = H(X, Y, Z)$ for all $X,Y,Z\in \Gamma (TM)$, we have $H^2(X,Y) = - \trwith{g}H_XH_Y$. Our convention for the scalar product of  alternating forms induced by the semi-Riemannian metric is 
that we consider forms as special tensors and use the 
usual scalar product on tensors. In particular,
\[ \absolute{H}_g^2 := H_{ijk}H^{ijk},\]
where indices were raised with the metric.  

Provided a generalised metric $\genmet$ and a divergence operator $\divergence$ on $E$, we often work with the generalised vector field $e \in \Gamma(E)$ defined via $\divergence = \divergence^\genmet - \scalbrack{e, \cdot}$. Due to the splitting induced by $\genmet$, we can then write $e = 2(X + \xi) \in \Gamma(\mathbb{T}M)$. Of particular importance will be the special case that $X= 0$ and $\xi = \extd \phi$ for some function $\phi \in C^\infty(M)$, which we sometimes refer to as \enquote{the dilaton}.

For the exterior curvature, our sign conventions are such that the second fundamental form is given by $k(X,Y)=g(\nabla^g_X n,Y)$, and the shape  tensor by $AX=\nabla^g_Xn$, $X,Y\in T\Sigma $, in terms of the Levi-Civita connection $\nabla^g$. Our sign convention for the Riemann tensor is such that $\riem(X,Y)Z = (\nabla^2_{X,Y} - \nabla^2_{Y,X}) Z$ and $\riem(W,Z,X,Y) = g(\riem(X,Y)Z,W)$. Then $\ric(X,Y) = \trwith{g} \riem(\cdot,X,\cdot,Y)$. We employ similar conventions in the \enquote{generalised} case.
\section{The Pullback of an Exact Courant Algebroid}\label{pullbackCAsection}

In this section, let $M$ be a smooth manifold, and $(E, \pi, \brackmap, \scalprodmap)$ an exact CA over $M$. Let furthermore $\iota \colon N \hookrightarrow M$ an immersion. We denote by $[H]\in H^3(M,\mathbb{R})$ the \v Severa class of $E$.

We recall the well-known definition of the pullback $\iota^! E$, cf. \cite[Lemma 3.7]{BursztynCavalcantiGualtieri2007}, \cite[Proposition 7.1.]{gualtieriPhDannals}.
\begin{definition}[Pullback CA]
    Denote by $K = \mathrm{Ann}(TN) \subset \iota^*(T^*M)$ the annihilator of $TN$. View it as a subbundle of $\iota^*E$ (by identifying $T^*M$ with $\pi^*T^*M \subset E$ via $\frac12 \pi^*$), and denote by $K^\perp \subset \iota^* E$ the orthogonal bundle with respect to the inner product on $E$. Then the \textit{pullback CA} $\iota^! E$ of $E$ is defined as the vector bundle
    \begin{equation*}
        \iota^! E \coloneqq \bigslant{K^\perp}{K}
    \end{equation*}
    equipped with the anchor $\iota^!\pi$, bracket $\iota^!\brackmap$ and inner product $\iota^!\scalprodmap$ that are naturally inherited from $E$.
\end{definition}
The bracket $\iota^![a,b]$ of $a, b \in \Gamma_{\mathrm{loc}}(\iota^!E)$ is defined as follows. Take $\tilde{a},\tilde{b} \in \Gamma_{\mathrm{loc}}(K^\perp)$ arbitrary lifts of $a,b$, and then $a_*,b_* \in \Gamma_{\mathrm{loc}}(E)$ smooth continuation of $\iota_*\tilde{a},\iota_*\tilde{b}$. Then $\iota^*[a_*,b_*] \in \Gamma_{\mathrm{loc}}(K^\perp)$ and $\iota^![a,b]$ is its projection to $\iota^!E$. 
One can easily show the following Lemma, see \cite[Lemma 3.7]{BursztynCavalcantiGualtieri2007}, \cite[Proposition 7.1.]{gualtieriPhDannals}
\begin{lemma}\label{pullbackCAlemma}
    The bracket is well-defined on $\Gamma (\iota^! E)$ and $(\iota^! E, \iota^!\pi, \iota^![\cdot ,\cdot ], \iota^!\scalprodmap)$ is an exact CA with \v Severa class $\iota^* [H]$.
\end{lemma}
It will at times be useful to work with a concrete realisation of the pullback CA, from which we will 
see why the pullback Courant algebroid $\iota^! E$ is exact. To give it, we introduce the following notion.

\begin{definition}[Courant Transversal Bundle]\label{Cou_trans:def}
    A \textit{Courant transversal bundle} $\mathcal{T} \subset \iota^*E$ for $N$ is a subbundle such that the sequence
    \begin{equation*}
        0 \longrightarrow K \stackrel{\pi^*}{\longrightarrow} \mathcal{T} \stackrel{\pi}{\longrightarrow} \nu \longrightarrow 0
    \end{equation*}
    is exact for some subbundle $\nu \subset \iota^* TM$ transversal to $TN$. A \textit{splitting} of an exact Courant algebroid $E$ over $M$ is an isomorphism $E \cong \mathbb{T}M$ of Courant algebroids, where the Dorfman bracket on the generalised tangent bundle is in general twisted. 
    We say that a splitting $F: E \cong \mathbb{T}M$ and a Courant transversal bundle are \textit{adapted}, if $F(\mathcal{T}) = \nu \oplus K$.
\end{definition}
\begin{remark}
    In particular, a Courant transversal bundle is of rank $\mathrm{rk}\, \mathcal{T} = 2\, \mathrm{codim}\; N= 2(\dim M -\dim N)$, and such that
    \begin{enumerate}[label = (\roman*)]
        \item $\nu \coloneqq \pi(\mathcal{T})$ is a transversal bundle for $N$, i.e. $\iota^*TM = TN \oplus \nu$,
        \item the annihilator $\mathrm{Ann}(TN) \subset \iota^*(T^*M)$ is contained in $\mathcal{T}$.
    \end{enumerate}
\end{remark}
\begin{lemma}\label{ad_splitting:lem}
    For every exact Courant algebroid $E$ over $M$ and any immersion $\iota : N \to M$ there exists a Courant transversal bundle $\mathcal T$. 
    \begin{enumerate}
        \item[(i)] For any splitting $F: E\cong \mathbb{T}M$ there is an adapted Courant transversal bundle $\mathcal T$. 
        \item[(ii)] For any Courant transversal bundle $\mathcal T$ there is an adapted splitting of $E$.
    \end{enumerate}
\end{lemma} 
\begin{proof} 
    Let $\nu\subset \iota^*TM$ be a subbundle transversal to $TN$ and $\tilde{\nu}$ any subbundle of $\iota^*E$ which projects isomorphically to $\nu$ under the anchor. Then $\mathcal T = K \oplus \tilde{\nu}$
    is a Courant transversal bundle. This proves the first claim. 
    
    Next we prove (i). 
    Given a splitting $F: E\cong \mathbb{T}M$ and a transversal bundle $\nu$, 
    $\mathcal T := K \oplus F^{-1}(\nu)$ is a Courant transversal bundle. 
    
    To prove (ii) let $F : E \cong \mathbb{T}M$ be a splitting and $\mathcal T$ a Courant transversal bundle. Then 
    $F(\mathcal T) = K\oplus \mathrm{graph} (\phi )$, where $\phi \in \Gamma (\mathrm{Hom} (\nu , \iota^*T^*M))$.  
    Decomposing $\phi = \phi_K + \psi$, according to 
    $\iota^*T^*M=K\oplus \mathrm{Ann} (\nu)$, we see that we can remove the component 
    $\psi\in \Gamma(\mathrm{Hom} (\nu ,\mathrm{Ann} (\nu)))\cong \Gamma (K\wedge \mathrm{Ann} (\nu))\subset \Gamma (\iota^*\bigwedge^2 T^*M)$ by a not necessarily closed B-field transformation on $M$. This yields a new splitting $F_1 : E\cong \mathbb{T}M$ such that 
    $F_1(\mathcal T) = K\oplus \mathrm{graph} (\phi_K ) = K\oplus \nu$. 
\end{proof}

\begin{proposition}\label{inducedCAprop}
    Let $\mathcal{T} \subset \iota^*E$ be a Courant transversal bundle. Then the orthogonal complement $E_N \coloneqq \mathcal{T}^\perp$ with respect to the inner product $\scalprodmap$ on $E$ is a concrete realisation of the quotient vector bundle $\iota^!E$ as a subbundle of $K^\perp$ complementary to $K$.

    In particular, we obtain the exact Courant algebroid $(E_N, \iota^*\pi, \brackmap_N, \iota^*\scalprodmap)$, which is canonically isomorphic to $\iota^! E$. Its bracket $\brackmap_N$ is defined by the formula
    \[ [u,v]_N \coloneqq \pi^\parallel ( \iota^* [\bar u, \bar v ]), \]
    where $\bar u, \bar v$ denote arbitrary extensions of $u, v\in \Gamma_{\mathrm{loc}} (E_N)\subset \Gamma_{\mathrm{loc}}  (\iota^* E)$
    to sections $\bar u ,\bar v\in \Gamma_{\mathrm{loc}}  (E)$ and $\pi^\parallel : \iota^*E = \mathcal T \oplus E_N\to E_N$ denotes the projection.
\end{proposition}
\begin{proof}
    To see that $E_N$ is a realisation of the quotient vector bundle $\iota^!E$, we choose an adapted splitting $F:E \to \mathbb{T}M$, which exists by Lemma~\ref{ad_splitting:lem}. Then 
    $F(\mathcal T)=K \oplus \nu$ and, hence, 
    \begin{equation*}
        \begin{split}
            K^\perp& = (F^{-1}(\mathrm{Ann}\:TN))^\perp = F^{-1}((\mathrm{Ann}\:TN)^\perp) = F^{-1}(TN \oplus \mathrm{Ann}\: \nu \oplus \mathrm{Ann}\:TN) \\
            &= F^{-1}((\nu \oplus \mathrm{Ann}\: TN)^\perp) \oplus K = \mathcal{T}^\perp \oplus K = E_N \oplus K
        \end{split}
    \end{equation*}
    Therefore, $E_N \cong \iota^!E$ as vector bundles. 
    
    To see that also $E_N \cong \iota^! E$ on the level of CAs, we note that, identifying $E_N$ and $\iota^!E$ as vector bundles, the inner product and the anchor defined on the two bundles are identical. It remains to see that also $\brackmap_N$ and $\iota^!\brackmap$ correspond to each other under the identification. This is obvious from the observation that 
    the identification $\iota^! E = K^\perp/K \cong E_N$ is precisely induced by the map 
    $\pi^\parallel|_{K^\perp}:K^\perp \to E_N$.
\end{proof}
\section{Inheritance of Geometrical Structures}\label{inheritedgeometrysection}
\subsection{Generalised Metrics}
In this chapter, we consider a semi-Riemannian immersion $N \to M$ into the base manifold of a semi-Riemannian exact CA $E\to M$. In this setup, we explain how $N$ naturally inherits the structure of a \textit{semi-Riemannian} exact CA. 

We recall the definition of a generalised (semi-Riemannian) metric \cite[Definition 2.16.]{krusche2022ricci}, compare   with the Riemannian case \cite[Definition 2.36.]{genricciflow}.
\begin{definition}[Generalised Metric] \label{genmet:def}
    A \textit{generalised metric} on on an exact Courant algebroid $E\to M$ is a non-degenerate section $\genmet \in \Gamma(\mathrm{Sym}^2\:E^*)$ such that 
    \begin{enumerate}[label = (\roman*)]
        \item $\restr{\genmet}{\Sym^2 \pi^* T^*M}$ is non-degenerate, and
        \item $(\genmet^\End)^2 = \idon{E}$.
    \end{enumerate}
    Herein, $\genmet^\End \in \Gamma(\End\: E)$ is the endomorphism obtained by composition of $\genmet \in \Gamma(\mathrm{Hom}(E,E^*))$  with the isomorphism $E^* \cong E$ determined by $\scalprodmap$.
\end{definition}
A generalised metric is equivalent to an orthogonal decomposition $E = E_+ \oplus E_-$ such that $\restr{\scalprodmap}{E_+}$ is non-degenerate. $E_+$ and $E_-$ are the Eigenbundle decomposition of $\genmet^\End$, corresponding by convention to the Eigenvalues $+1$ and $-1$ respectively. We denote by $\pi_\pm \coloneqq \frac{1}{2}(\id \pm \genmet^\End)\colon E \to E_\pm$ the orthogonal projections onto $E_\pm$.
\begin{example}\label{genmetricexample}
    Let $g$ be a semi-Riemannian metric on $M$. On the (twisted) generalised tangent bundle $\mathbb{T}M$, one can then define the endomorphism $\genmet_g^\End \in \Gamma(\End\: \mathbb{T}M)$ via
    \begin{equation*}
        \genmet_g^\End (X+ \xi) = \begin{pmatrix}
            0 & g^{-1} \\
            g & 0 \\
        \end{pmatrix} \begin{pmatrix}
            X \\ \xi
        \end{pmatrix}
        , \qquad\quad X+ \xi \in \Gamma(\mathbb{T}M)
    \end{equation*}
    Then $\genmet_g \coloneqq \scalbrack{\genmet_g^\End,\cdot}$ defines a semi-Riemannian generalised metric on $E$.
\end{example}

For our discussion, it will be essential to know the following decomposition of semi-Riemannian generalised metrics into splittings and semi-Riemannian metrics on the base, see e.g. \cite[§ 3.2]{krusche2022ricci}.
\begin{proposition}\label{genmetricprop}
    There is a one-to-one correspondence between semi-Riemannian generalised metrics $\genmet$ on $E\to M$ and tuples $(g,F)$ consisting of a semi-Riemannian metric $g$ on $M$, and a splitting $F\colon E \cong \mathbb{T}M$. Furthermore, $\genmet = F^*\genmet_g$ where $F^*\genmet_g$  denotes the pullback\footnote{Explicitly, $(F^*\genmet_g)(a,b) \coloneqq \genmet_g(Fa,Fb)$ for $a,b \in \Gamma(E)$.} of $\genmet_g$ with $F$, and $F(E_\pm) = \{X \pm g X \mid X \in \Gamma(TM)\}$.
\end{proposition} 
\begin{corollary}\label{sigmapmcor}
    Let $(g,F)$ be the tuple corresponding to a generalised metric $\genmet$. Then, the maps $\sigma_\pm\colon TM \to E_\pm, X \mapsto F^{-1}(X \pm g X)$ are isometries with respect to the inner products $g$ and $\genmet$.
\end{corollary}
\begin{remark}
    Note that, under the above correspondence, the signature of $\genmet$ is twice the signature of $g$. We will call $\genmet$ \emph{Riemannian} (respectively, \emph{Lorentzian}) if $g$ is Riemannian (respectively, Lorentzian).  
\end{remark}
From now on we view $M$ as being equipped with $g$ and consider a semi-Riemannian immersion $\iota : N \to M$, which means that $h=\iota^*g$ is nondegenerate.

\begin{remark}\label{courantnormalbundlerem}
    Note that a generalised metric $\genmet \cong (g, F)$ provides a canonical choice of Courant transversal bundle on $N \hookrightarrow M$, the \emph{Courant normal bundle} $\mathcal{N} \coloneqq K \oplus \genmet^\End K = F^{-1}(\nu \oplus g \nu)$. Herein, $\nu$ denotes the normal bundle induced by $g$. In particular, we obtain the canonical representation $E_N \coloneqq \mathcal{N}^\perp$ of the pullback CA $\iota^! E$.
\end{remark} 

\begin{proposition}\label{inducedsemiriemCAprop}
    The generalised metric $\genmet$ on $E$ induces a generalised metric $\genmetHS$ on $E_N$. Furthermore, over $N$, $\genmet$ is completely determined by $\genmetHS$ together with an orthonormal set of generalised vector fields $\{n_i\} \subset \Gamma(\iota^* E)$ such that $\mathcal N = \spanspace \{n_i\} \oplus K$.
\end{proposition}
\begin{proof}
    We claim that the symmetric bilinear form $\genmetHS := \restr{\iota^*\genmet}{E_N\times E_N}$ defines a generalised metric. Noting that $\mathcal N = \mathcal N \cap E_+ \oplus \mathcal N \cap E_-$, we immediately conclude that $\genmetHS$ is non-degenerate, as the decomposition $E = \mathcal{N} \oplus E_N$ is orthogonal also with respect to $\genmet$. 
    
    We now check that $\restr{\genmetHS}{\Sym\:T^*N}$ is non-degenerate. We remind ourselves that $T^*N$ is identified with the subbundle $\mathrm{Ann} (\nu )\subset \iota^*T^*M$. With this identification, the restriction $\pi^*|_{T^*N}: T^*N \to \iota^*E$ maps to $E_N$ such that the resulting map $T^*N\to E_N$  can be identified with the adjoint of the anchor $E_N\to TN$. Therefore, $\restr{\genmetHS}{\Sym^2\:T^*N} = \restr{\genmet}{\Sym^2\:(\mathrm{Ann}\:\nu)}$, and non-degeneracy of this restriction follows again from orthogonality of $E = \mathcal{N} \oplus E_N$, since $\mathrm{Ann}\:\nu\subset E_N$.
    
    To see that $\genmetHS^\End$ is involutive, note that the subbundles $\mathcal N, E_N = \mathcal N^\perp \subset \iota^*E$ are invariant under $\genmet^\End$ and non-degenerate with respect to $\genmet$. Therefore, $\genmetHS^\End = \iota^*\genmet^\End|_{E_N} : E_N \to E_N$, implying that $(\genmetHS^\End)^2 = \id$. We have concluded that $\genmetHS$ is indeed a generalised metric.

    Finally, note that $g$ is determined by the induced metric $h$ and a unit normal frame $\{U_i\}$. Thus, specifying the generalised metric $\genmet \cong (g,F)$ is equivalent to specifying a generalised metric $\genmetHS \cong (h,\restr{F}{E_N})$ and a set of orthonormal generalised vector fields $\{F^{-1}(U_i)\}$ that spans, together with $K$, $\mathcal{N}$.
\end{proof}

\begin{corollary} The generalised metric $\genmet$ on $E$ induces a generalised metric $\iota^!\genmet$ on $\iota^!E$. 
\end{corollary}
\begin{proof}
It suffices to define $\iota^!\genmet$ as the generalised metric on $\iota^!E$ which corresponds to the generalised metric $\genmetHS$ of 
Proposition \ref{inducedsemiriemCAprop} under the canonical
isomorphism $\iota^!E \cong E_N$. 
\end{proof}

\subsection{Divergence Operators}

Denote again by $(E,\genmet)$ a semi-Riemannian exact CA over a manifold $M$, naturally equipped with the semi-Riemannian metric $g$ associated to $\genmet$. In this section, we explain how, under a semi-Riemannian immersion $N\hookrightarrow (M,g)$, a divergence operator on $E$ induces a divergence operator $\divergence_N$ on $E_N$.

We start by recalling the following fundamental concept \cite{fernandez}.
\begin{definition}[Divergence Operator]
    A divergence operator on an exact CA $E \to M$ is a map $\divergence \colon \Gamma(E) \to C^\infty(M)$ such that $\divergence(f a) = \pi a(f) + f \divergence(a)$ for all $a \in \Gamma(E)$.
\end{definition}
\begin{example}\label{metridivexample}
    Let $E \to M$ be an exact CA, and $\nabla$ be a connection $M$. Then $\Gamma(E)\ni a \mapsto \tr (\nabla \pi a)$ defines a divergence operator on $E$. If $\nabla$ is the LC connection with respect to the metric $g$ associated to a generalised metric $\genmet$ on $E$, then this divergence operator is called the \textit{metric divergence} and denoted $\divergence^\genmet$.  
\end{example}
It is immediate from the definition that the space of divergence operators forms an affine space over $\Gamma(E^*) \cong \Gamma(E)$. It follows that, on semi-Riemannian exact CAs, we can identify a divergence operator $\divergence$ with the section $e\in \Gamma(E)$ such that $\divergence = \divergence^\genmet - \scalbrack{e, \cdot}$.
\begin{definition}[Induced Divergence Operator]
    Let $\divergence = \divergence^\genmet - \scalbrack{e, \cdot}$ be a divergence operator on an exact semi-Riemannian CA $(E, \genmet)$ with induced semi-Riemannian metric $g$ on $M$. Let $\iota\colon N \hookrightarrow (M,g)$ be a semi-Riemannian immersion with Courant normal bundle $\mathcal{N}$. Then we define on $E_N = \mathcal{N}^\perp$ the divergence operator
    \begin{equation*}
        \divergence_N \coloneqq \divergence^\genmetHS - \langle e^\parallel, \cdot \rangle
    \end{equation*}
    Herein, $e^\parallel$ is the orthogonal projection of $\iota^*e$ into $E_N$, and $\genmetHS$ is the generalised metric induced on $E_N$ by $\genmet$.
\end{definition}
We note that, in general, $\divergence_N$ depends on the choice of generalised metric $\genmet$. One can reduce this to a dependence only on the metric $g$ in the tuple $\genmet \cong (g, F)$, if one demands the expression $\scalbrack{e, \cdot}$ to uniquely define an element in $(\iota^!E)^*$. This amounts to $\iota^*e \in \Gamma(K^\perp)$ or, equivalently,  $\pi(\iota^*e) \in \Gamma(TN)$. 

Recall that a pair $(\genmet , \divergence )$ consisting of a generalised metric and a divergence operator is called \textit{compatible}  
if the generalised vector field $e =2(X+\xi) \in \Gamma(E)$ defined by the equation $\divergence = \divergence^\genmet - \scalbrack{e,\cdot}$ is generalised Killing, $[e,\genmet] = 0$  \cite[Lemma 3.32]{genricciflow}. 
Compatibility of $(\genmet, \divergence)$ can be seen to be equivalent to (see \cite[Proposition 2.53.]{genricciflow} for the Riemannian case, however the proof works in any signature)
    \begin{equation} \label{compa:eq}
        0 = L_X g, \qquad\qquad \extd \xi = H(X)
    \end{equation}
It is called \textit{closed} if $X=0$ and $\extd \xi = 0$, and \textit{exact} if $X=0$ and $\xi = \extd \phi$ for some $\phi \in C^\infty(M)$.

\begin{proposition}\label{induceddivergenceopsprop}
    Let $\divergence = \divergence^\genmet - \scalbrack{e, \cdot}$ be a divergence operator on the exact semi-Riemannian CA $(E, \genmet)$. Assuming $\iota^*e \in \Gamma(K^\perp)$,  the induced pair $(\genmetHS, \divergence_N)$ is compatible, closed, and exact if the pair $(\genmet, \divergence)$ has the respective attribute.
\end{proposition}
\begin{proof}
    We work in the splitting provided by the generalised metric, and denote $e = 2(X+\xi)$. The condition $\iota^*e \in \Gamma(K^\perp)$ translates to $\iota^*X = X^\parallel \in \Gamma(TN)$, where $X^\parallel$ denotes the component of $\iota^*X$ tangent to $N$. 
    
    Throughout this proof, we denote $\varphi^\parallel := \iota^*\varphi|_{\bigwedge^k TN}\in \Gamma(\bigwedge^kT^*N)$ for any $\varphi \in \Gamma(\bigwedge^kT^*M)$. In particular, $\xi^\parallel$ is given by $\iota^*\xi|_{TN}\in \Gamma (TN)\cong \Gamma (\mathrm{Ann}(\nu ))$.

    First we assume that $(\genmet , \divergence )$ is compatible and  deduce that $\genmetHS$ is compatible with $\divergence_N = \divergence^{\genmetHS} - 2\scalbrack{X^\parallel + \xi^\parallel, \cdot}$. In fact, from \eqref{compa:eq} see that
    \begin{equation*}
        L_{X^\parallel}h = \restr{L_X g}{TN} = 0
    \end{equation*}
    and
    \begin{equation*}
        \extd^N \xi^\parallel = (\extd \xi)^\parallel = (H(X))^\parallel = H^\parallel(X^\parallel)
    \end{equation*}
    wherein $\extd^N$ denotes the exterior derivative on $N$. This proves compatibility of $(\genmetHS, \divergence_N)$.

    Assume next $(\genmet, \divergence)$ to be closed. Then closedness of the induced pair immediately follows from $\extd^N \xi^\parallel = (\extd \xi)^\parallel = 0$.

    Finally, exactness of the induced pair is an immediate consequence of exactness of $(\genmet, \divergence)$, as with $\xi = \extd \phi$
    \begin{equation*}
        \xi^\parallel = (\extd \phi)^\parallel = \extd^N \iota^* \phi 
    \end{equation*}
    This finishes the proof.
\end{proof}
\begin{remark}
    By a slight generalisation of the above proof, we see that Proposition~\ref{induceddivergenceopsprop} holds under the weaker assumption that the normal component $X^\perp\in \Gamma (\nu)$ of $X=\pi e/2\in \Gamma (TM)$ extends to a Killing field 
    $X^\perp\in \Gamma (TM)$ satisfying $H(X^\perp)=0$. 
\end{remark}
\subsection{Generalised Connections}

Here we explain the inheritance of generalised connections. We assume as before that $(E,\genmet)$ is a semi-Riemannian exact CA over $M$, and that $\iota \colon N \hookrightarrow (M,g)$ is a semi-Riemannian immersion with respect to the semi-Riemannian metric $g$ associated with $\mathcal G$. We respectively denote by $\genmetHS$ and $h$ the induced generalised metric on $E_N$ and the usual metric $N$.

We recall the definition of generalised connections.
\begin{definition}[Generalised Connection]
    A \textit{generalised connection} on $E$ is a map $D\colon \Gamma(E) \times \Gamma(E) \to \Gamma(E), (u,v) \mapsto D_u v$ satisfying for all $f\in C^\infty(M)$ and $u,v,w \in \Gamma(E)$
    \begin{enumerate}[label = (\roman*)]
        \item $C^\infty$-linearity in the first entry, i.e. $D_{fu} v = f D_u $,
        \item the Leibniz rule $D_u (fv) = \pi u(f) v + f D_u v$, and
        \item compatibility with the inner product, i.e. $\pi u \scalbrack{v,w} = \scalbrack{D_u v, w} + \scalbrack{v, D_u w}$.
    \end{enumerate}
\end{definition}
Given a generalised metric $\genmet$, a generalised connection $D$ is called \textit{$\genmet$-metric}, if $D \genmet = 0$.

Every generalised connection $D$ is assigned a \textit{torsion-tensor} $T^D \in \Gamma({\bigwedge}^3E^*)$:
\begin{equation}\label{torsion_def:eq}
    T^D(u,v,w) \coloneqq \scalbrack{D_u v - D_v u - [u,v], w} + \scalbrack{D_wu,v}, \qquad\quad u,v,w \in \Gamma(E)
\end{equation}
The last term is needed to achieve tensoriality.

Furthermore, every generalised connection $D$ is assigned a divergence operator $\divergence_D$, given as $\Gamma(E)\ni a \mapsto \tr D a$. Given any divergence operator $\divergence$ on $E$, the pair $(D, \divergence)$ is called \textit{compatible} if $\divergence_D = \divergence$.

\begin{example}\label{genconnexample}
    On the generalised tangent bundle $\mathbb{T}M$ endowed with the Dorfman bracket associated with a closed $3$-form $H$, the map $(X+\xi,Y+\eta) \mapsto \nabla_X(Y+\eta)$ defines a generalised connection $D$. It is torsion-free if and only if $\nabla$ is torsion-free and $H=0$. In fact, $T^D(X+\xi, Y+\eta , Z+\zeta ) = \frac12 (\zeta (T^\nabla (X,Y)) -H(X,Y,Z))$. Given a semi-Riemannian metric $g$ on $M$ and the associated generalised metric $\genmet_g$ from Example \ref{genmetricexample}, it is $\genmet_g$-metric if and only if $\nabla$ is $g$-metric. Finally, it is compatible with the metric divergence $\divergence^\genmet$ from Example \ref{metridivexample}.
\end{example}

Given a tuple $(E,\genmet, \divergence)$ consisting of an exact CA $E$ equipped with a generalised metric $\genmet$ and a divergence operator $\divergence$, we call a generalised connection $D$ \textit{Levi-Civita with divergence $\divergence$}, if it torsion-free, metric, and divergence compatible. We denote the space of these generalised LC connections by $\mathcal{D}^0(\genmet, \divergence)$, which was described in \cite{fernandez}, compare \cite[§ 3]{genricciflow}.

For conceptual reasons and wider applicability we rephrase the result using the notion of generalised first prolongation. Recall \cite{cortesdavid} that  given a pseudo-Euclidean vector space $(V,\langle \cdot ,\cdot \rangle )$ and a Lie subalgebra $\mathfrak g \subset \mathfrak{so}(V)$, the \emph{generalised first prolongation} of $\mathfrak g$ is defined as 
\begin{equation} \label{1stprol:eq}\mathfrak g^{\langle 1\rangle} := \{ \chi \in V^*\otimes \mathfrak g \mid \partial_{\mathfrak g} \chi =0\},\end{equation}
where $\partial_\mathfrak g : V^*\otimes  \mathfrak g \to {\bigwedge}^3V$ is the restriction of the $\mathrm O(V)$-equivariant map 
\[ \partial=\partial_{\mathfrak{so}(V)}: V^*\otimes \mathfrak{so} (V)\cong V^* \otimes {\bigwedge}^2 V^* \to {\bigwedge}^3 V^*,\]
defined by $(\partial \chi )(u,v,w) = 
\sum_{\mathrm{cycl}}\langle \chi (u,v),w\rangle$ and 
the sum is over the cyclic permutations of $u,v,w\in V$. As a general fact \cite{cortesdavid}, the space of generalised connections on a Courant algebroid $E$ over $M$ compatible with a given reduction of the structure group $\mathrm O(k,\ell)$ to some Lie subgroup $G\subset \mathrm O(k,\ell )$ and with prescribed torsion $T\in \bigwedge^3 E^*$ is either empty or an affine space modeled on the space of sections of the vector bundle $\mathfrak g(E)^{\langle 1\rangle}\to M$. Here $\mathfrak g(E)\subset \mathfrak{so}(E)$
denotes the subbundle of which the fiber  over a point $p\in M$ is identified with $\mathfrak g$ by means of any $G$-frame of $E_p$. 

In particular, the affine space of generalised Levi-Civita  connections of 
a Courant algebroid endowed with a generalised metric $\mathcal G$ is 
modeled on the space of sections of 
\begin{equation}\label{so1:eq}
   (\mathfrak{so}(E)_{\mathcal G})^{\langle 1\rangle}=\mathfrak{so}(E_+)^{\langle 1\rangle}\oplus \mathfrak{so}(E_-)^{\langle 1\rangle},
\end{equation}
where $\mathfrak{so}(E)_{\mathcal G}=\mathfrak{so}(E_+)\oplus \mathfrak{so}(E_-)$ denotes the stabiliser 
of $\mathcal G$ in $\mathfrak{so}(E)$. 
To prove this it suffices to construct a 
single generalised Levi-Civita connection, 
which was done in \cite{fernandez}. As reviewed below, there is even a generalised Levi-Civita connection with prescribed divergence operator (under the assumption $\dim M>1$). 

Prescribing, in addition, a divergence operator therefore cuts out an affine subspace modeled on the space of sections of the subbundle $\mathfrak{so}(E_+)_0^{\langle 1\rangle}\oplus \mathfrak{so}(E_-)_0^{\langle 1\rangle}\subset \mathfrak{so}(E_+)^{\langle 1\rangle}\oplus \mathfrak{so}(E_-)^{\langle 1\rangle}$ of (totally) trace-free tensors. More explicitly, 
\begin{equation*} 
     \mathfrak{so}(E_\pm)_0^{\langle 1\rangle}:= \{\chi \in \mathfrak{so}(E_\pm )^{\langle 1\rangle}\mid \mathrm{tr}\, \chi =0\},
\end{equation*}  
where $\tr \chi \in \Gamma(E_\pm^*)$ denotes the non-trivial trace, that is 
\[ (\tr \chi)(v) := \tr (\chi v) =  \tr (u\mapsto \chi_u v),\]  
$v\in E_\pm$. Here we have used the standard notation $\chi_uv := \chi (u,v)$.
Summarising, we have the following result of \cite{fernandez}, which requires $\dim M>1$.
\begin{lemma}\label{genlcspacelemma}
    The set $\mathcal{D}^0(\genmet, \divergence)$ is an affine space over the space of sections of the vector bundle $\Sigma_0^+ \oplus \Sigma_0^-$, where
    $\Sigma_0^\pm \coloneqq \mathfrak{so}(E_\pm )_0^{\langle 1\rangle}$.
\end{lemma}

As announced,  we construct now a (canonical) generalised LC connection for any divergence operator, following \cite{fernandez} and \cite[§3]{genricciflow}.

Let $(M,g)$ be a semi-Riemannian manifold with LC connection $\nabla$ and $H \in \Omega^3(M)$. Then, we consider the connections $\nabla^\pm$ and $\nabla^{\pm 1/3}$ on $M$ such that 
\begin{equation}\label{nablapmdefeq}
    \nabla^\pm_X  \coloneqq \nabla_X  \pm \frac{1}{2} H_X, \qquad\quad \nabla^{\pm 1/3}_X \coloneqq \nabla_X  \pm \frac{1}{6} H_X, \quad  X \in \Gamma(TM),
\end{equation}
{\color{blue}}
where $H_X \in \Gamma (\mathfrak{so}(TM))$ is defined by $g(H_XY,Z) = H(X,Y,Z)$, $ Y,Z\in \Gamma(TM)$.

With this notation and taking for $H$ the  closed $3$-form associated with an exact CA $E$ endowed with a generalised metric $\genmet$ we state \cite[Lemma 3.17.]{genricciflow} 
\footnote{The result is stated in \cite{genricciflow}  for 
positive-definite $\genmet$. However, the construction and proof solely rely on the weaker properties  included in Definition \ref{genmet:def}, and thus work entirely independent of signature.}
\begin{lemma}\label{cangenconnlemma}
    Let $\divergence = \divergence^\genmet - \scalbrack{e, \cdot}, e \in \Gamma(E)$, be a divergence operator on  $(E, \genmet )$. Denote by $e_\pm \coloneqq \pi_\pm e \in \Gamma(E_\pm)$ the orthogonal projections, and by $\sigma_\pm \colon TM \to E_\pm$ the isometries from Corollary \ref{sigmapmcor}.  Then
    \begin{enumerate}[label = (\roman*)]
        \item the generalised connection $D^0$ defined by requiring for all $a_\pm,b_\pm \in \Gamma(E_\pm)$
            \begin{equation*}
                D^0_{a_\pm}b_\pm \coloneqq \sigma_\pm \nabla^{\pm 1/3}_{\pi a_\pm} \pi b_\pm, \qquad\quad D^0_{a_\mp}b_\pm \coloneqq \sigma_\pm \nabla^{\pm}_{\pi a_\mp} \pi b_\pm
            \end{equation*}
            is generalised LC with respect to $(\genmet, \divergence^\genmet)$, i.e.\  $D^0 \in \mathcal{D}^0(\genmet, \divergence^\genmet)$, and
        \item the generalised connection $D \coloneqq D^0 + \frac{\chi_+^{e_+} + \chi_-^{e_-}}{\dim M - 1}$, where $\chi_\pm \in \Gamma(\mathfrak{so}(E_\pm)^{\langle 1\rangle })$ is defined by
        \begin{equation*}
            \chi_\pm^{e_\pm}(a,b) \coloneqq \scalbrack{a,b} e_\pm - a \scalbrack{e_\pm, b}, \qquad\quad a,b \in \Gamma(E_\pm)
        \end{equation*}
        is generalised LC with respect to $(\genmet, \divergence)$, i.e.\ $D \in \mathcal{D}^0(\genmet, \divergence)$.
    \end{enumerate}
\end{lemma}
See \cite{vicentethomas} for another description of the above  generalised Levi-Civita connection $D$ that provides a geometric interpretation, which justifies to call it \emph{the} canonical generalised Levi-Civita connection.
The following lemma  establishes the basis of our discussion.
\begin{lemma}\label{gausslemma1}
    Denote by $\pi^\parallel \colon \iota^*E \to E_N$ the projection along the Courant normal bundle $\mathcal{N}$. Then restriction and projection $D \mapsto D^N = \pi^\parallel \restr{D}{E_N}$ of generalised connections on $E$  defines a map to the space of generalised connections on $E_N$. More precisely, given $u,v\in \Gamma_{\mathrm{loc}} (E_N)$ and extensions $\tilde{u}, \tilde{v}\in \Gamma_{\mathrm{loc}} (E)$ we define
    \[ D^N_uv := \pi^\parallel (D_{\tilde{u}}\tilde{v}).\] 
    The map  $D \mapsto D^N$ maps $\genmet$-metric to $\genmetHS$-metric generalised connections and torsion-free  to torsion-free generalised connections.  
    
    Furthermore, it maps the canonical connection $D^0 \in \mathcal{D}^0(\genmet, \divergence^\genmet)$ to the canonical connection $D^{N0}\in  \mathcal{D}^0(\genmetHS, \divergence^\genmetHS)$, thus $D^{N0} = (D^0)^N$.
\end{lemma}
\begin{proof} Using the fact that $\iota^*\pi u\in \Gamma (TN)$, it is easy to see that $\iota^*(D_{\tilde{u}}\tilde{v})$, and hence $\pi^\parallel(D_{\tilde{u}}\tilde{v})$, is independent of the extensions $\tilde{u}, \tilde{v}$ of $u,v\in \Gamma_{\mathrm{loc}} (E_N)$. This shows that the generalised connection $D^N= \pi^\parallel \restr{D}{E_N}$ is well-defined for any generalised connection $D$ on $E$. We recall that $\mathcal N = \mathcal N \cap E_+ \oplus \mathcal N \cap E_-$ and, hence,  the decomposition $\iota^*E = \mathcal N \oplus E_N$ is orthogonal for both, $\langle \cdot ,\cdot \rangle$ and $\genmet$. Therefore, $D^N$ is always a generalised connection on $E_N$ and is metric if $D$ is: 
    \begin{eqnarray*}
        \pi(u)\scalbrack{v,w} &=& \scalbrack{D_u v,w} + \scalbrack{v, D_u w} = \scalbrack{D^N_u v, w} + \scalbrack{v, D^N_u w} \\
        \pi(u)\genmet (v,w) &=& \genmet (D_u v,w) + \genmet (v, D_u w) = \genmet (D^N_u v, w) + \genmet (v, D^N_u w) 
    \end{eqnarray*}
    for all $u,v,w \in \Gamma(E_N)$. A similar calculation shows that the torsion $3$-forms of $D$ and $D^N$ are related by $T^{D^N}=(T^D)^\parallel$. In particular, $T^D=0$ implies $T^{D^N}=0$.

    For the last part of the statement, observe that the eigenbundles of $\mathcal H^\End$ are given by $(E_N)_\pm:= \pi^\parallel E_\pm$, and consider the explicit description of the pure-type and mixed-type parts of $D^0$ and $D^{N0}$ from Lemma \ref{cangenconnlemma}. Since the LC connection on $(M,g)$ gets mapped onto the LC connection on $(N, h)$ under restriction and projection, and since the twist on $E_N$ is given by $H^\parallel=\iota^*H|_{\bigwedge^3TN}$, the result follows.  
\end{proof}
Going on a brief tangent, we apply this result to show that, under natural conditions on the semi-Riemannian immersion $N\hookrightarrow M$, a generalised Kähler structure over $M$ induces one over $N$.
\begin{corollary}\label{GK_preserved:cor}
    Let $(E,\mathcal G, \mathcal J)$ be an exact Courant algebroid over $M$ endowed with a (possibly indefinite) generalised K\"ahler structure. Furthermore let $N\subset M$ be 
    a semi-Riemannian submanifold and $(E_N,\mathcal H)$ the induced exact semi-Riemannian Courant algebroid. Assume that  
    $\mathcal J E_N \subset  E_N$. Then $(E_N,\mathcal H, \mathcal J|_{E_N})$ is generalised K\"ahler. 
\end{corollary}
\begin{proof} According to \cite{cortesdavid} the 
generalised K\"ahler property for a generalised almost Hermitian structure $(\mathcal G, \mathcal J)$ on a Courant algebroid $E$ is equivalent to the existence of a LC generalised connection such that $D\mathcal J =0$. By virtue of Lemma \ref{gausslemma1} we know that the induced generalised connection $D^N$ is a Levi-Civita generalised connection for $\mathcal H$.
Therefore it suffices to show that $D^N\mathcal J_N=0$ for the generalised almost complex structure $\mathcal J_N := \mathcal J|_{E_N}$. 
Let $u,v\in \Gamma (E_N)$. Then 
\[ (D^N_u\mathcal J_N)v = D^N_u(\mathcal J_Nv) -\mathcal J_ND_u^Nv=\pi^{\parallel} (D_u\mathcal J) v=0\]
\end{proof}
Similarly, using the characterisation in \cite{cortesdavid} of generalised hyper-K\"ahler structures  as generalised almost hyper-Hermitian structures $(\mathcal G, \mathcal J_1, \mathcal J_2, \mathcal J_3)$ admitting a LC generalised connection $D$ such that $D\mathcal J_\alpha=0$ for all $\alpha \in \{ 1,2,3\}$  we obtain the following.
\begin{corollary}
    Let $(E,\mathcal G, \mathcal J_1, \mathcal J_2, \mathcal J_3)$ be an exact Courant algebroid over $M$ endowed with a (possibly indefinite) generalised hyper-K\"ahler structure. Furthermore let $N\subset M$ be 
    a semi-Riemannian submanifold and $(E_N,\mathcal H)$ the induced exact semi-Riemannian Courant algebroid. Assume that  
    $\mathcal J_\alpha E_N \subset  E_N$ for all $\alpha \in \{ 1,2,3\}$. Then $(E_N,\mathcal H, \mathcal J_\alpha|_{E_N}, \alpha = 1,2,3)$ is generalised hyper-K\"ahler. 
\end{corollary}

\begin{remark}
    The map $D \mapsto D^N$ does not, in general, preserve compatibility with a given divergence operator $\divergence$. Consider $D^0\in \mathcal{D}^0(\genmet, \divergence^\genmet)$ the canonical generalised connection from Lemma \ref{cangenconnlemma}. Then, by Lemma \ref{gausslemma1}, $(D^0)^N=\pi^\parallel D^0 = D^{N0} \in \mathcal{D}^0(\genmetHS, \divergence^\genmetHS)$ is the canonical generalised connection with divergence operator $\divergence_N= \divergence^\genmetHS$. However given an arbitrary divergence operator $\divergence = \divergence^\genmet - \scalbrack{e, \cdot}$, the reduced divergence operator is  
    $\divergence_N = \divergence^\genmetHS - \scalbrack{e^\parallel, \cdot}$. We claim that it does not in general coincide with the divergence
    operator of the generalised LC connection $D^N$ obtained from a generalised LC connection $D$ with divergence operator $\divergence$. Write  
    $D=D^0 +\chi$, where $\chi\in \Gamma (\mathfrak{so}(E)_{\mathcal G}^{\langle 1\rangle })$ is such that 
    $ \mathrm{tr}(\chi v) = -\langle e,v\rangle$ for all $v\in E$, compare  \eqref{1stprol:eq} and \eqref{so1:eq}. It follows for all $v\in \Gamma (E_N)$ that 
    \[ \divergence_{D^N}(v)=\divergence_{D^{N0}}(v) + \mathrm{tr}_{E_N} \pi^\parallel (\chi v) \]
    where $\mathrm{tr}_{E_N} \pi^\parallel (\chi v) = 
    \mathrm{tr}_{E_N} (\chi v) = \mathrm{tr} (\chi v)-
    \mathrm{tr}_{\mathcal N}(\chi v) = -\langle e^\parallel , v \rangle - \mathrm{tr}_{\mathcal N}(\chi v)$ and, in general, $\mathrm{tr}_{\mathcal N}(\chi v)\neq 0$.
    In fact, a particular choice, which corresponds to the     canonical LC connection in Lemma~\ref{cangenconnlemma} with divergence 
    operator $\divergence$ is $\chi= \frac{1}{\dim M - 1}(\chi_+^{e_+} + \chi_-^{e_-})$. For this choice and 
    $v\in (E_N)_+$ we obtain 
    \[ (\dim\: M - 1) \mathrm{tr}_{\mathcal N}(\chi v) = \sum n_i^*(\langle n_i,v\rangle e_+-\langle e_+,v\rangle n_i)= -(\dim M -\dim N)\langle e_+^\parallel ,v\rangle\]
    where $(n_i)$ is a (local) orthonormal frame of $\mathcal N_+:= \mathcal N \cap E_+$ and $(n_i^*)$ is the dual frame of $\mathcal N_+^*$. As a consequence we obtain for the canonical generalised LC connection $D$ with $\divergence_D=\divergence$,  
    \[ \divergence_{D^N}=\divergence_{D^{N0}} - 
    \frac{\dim N -1}{\dim M-1}\langle e^\parallel, \cdot \rangle = \divergence^{\genmetHS} -  
    \frac{\dim N -1}{\dim M-1}\langle e^\parallel, \cdot \rangle\]
    which is different from $\divergence_N = \divergence^\genmetHS - \scalbrack{e^\parallel, \cdot}$, as soon as $e^\parallel \neq 0$ (we do always assume $\dim N>1$ and hence $\dim M-1\neq0$).
\end{remark}
We have proven the following.
\begin{proposition}\label{divdiscompatibilityprop}
    Let $\divergence$ be any divergence operator on the semi-Riemannian CA $E$ and $D$ the canonical generalised LC connection with divergence operator $\divergence= \divergence^\genmet -\langle e,\cdot \rangle $. Then the reduced generalised LC connection $D^N$ has the divergence operator 
    \[ \divergence_{D^N}=\divergence^{\genmetHS} -  
        \frac{\dim N -1}{\dim M-1}\langle e^\parallel ,\cdot \rangle 
        \]
        while the reduction of $\divergence$ is 
        \[ \divergence_N = \divergence^{\genmetHS} -\langle e^\parallel ,\cdot \rangle\]
        More generally, for any generalised LC connection $D=D^0+\chi$, $D^0\in \mathcal D^0(\genmet, \divergence^{\mathcal G} )$ canonical,  with divergence operator $\divergence$ we have 
     \begin{equation} \label{div_deviation:eq}  \divergence_{D^N}(v) = \divergence_N(v) - \mathrm{tr}_{\mathcal N}(\chi v),\quad v\in \Gamma (E_N)\end{equation} 
\end{proposition}
It follows that $\mathrm{tr}_{\mathcal N}(\chi v) = - \frac{\dim M -\dim N}{\dim M -1}\langle e^\parallel, \cdot \rangle$, for the difference tensor $\chi =D-D^0$ associated with the canonical $D\in \mathcal D^0(\genmet, \divergence )$ but it can have other  values for other choices of $D\in \mathcal D^0(\genmet, \divergence )$. 
Thus we have the following. 
\begin{corollary}\label{cangenLC_divergence_submf_lemma}
Let $D$ be the canonical generalised Levi-Civita connection on $(E,\mathcal G)$ with divergence $\divergence=\divergence^{\genmet}-\langle e,\cdot \rangle $, $e\in \Gamma (E)$. Then 
\[ \divergence_{D^N} = \divergence_N + \frac{\dim M -\dim N}{\dim M -1}\langle e^\parallel ,\cdot \rangle.\]

\end{corollary}

We denote by $\mathcal{D}^0(\genmet, \divergence; N)$ the subspace of generalised connections in $\mathcal{D}^0(\genmet,\divergence)$ that project onto connections in $\mathcal{D}^0(\genmetHS, \divergence_N)$. 

\begin{proposition}
    $\mathcal{D}^0(\genmet, \divergence; N)$ is an affine space over the space of sections $\Gamma_N(\Sigma_0) :=  \Gamma_N(\Sigma_0^+) \oplus \Gamma_N(\Sigma^{-}_0)$, where
    \begin{equation*}
        \Gamma_N(\Sigma_0^\pm) = \big\{ \chi_\pm \in \Gamma(\Sigma_0^\pm) \mid \trwith{E_N} \iota^*\chi_\pm = 0 \big\}
    \end{equation*}
\end{proposition}

\begin{proof}
 The next lemma establishes that $\mathcal{D}^0(\genmet, \divergence; N)$ is non-empty. The claim is a consequence.  
\end{proof}

\begin{lemma}\label{HScanconnlemma}Let $D^0\in \mathcal{D}^0(\genmet, \divergence^{\mathcal G})$  be the canonical generalised LC connection with metric divergence. 
    Define the generalised connection $D = D^0 + \frac{\chi^N_+ + \chi^N_-}{\dim M -1}$, where $\chi_\pm^N \in \Gamma(\Sigma_0^\pm)$ is such that on $N$
    \begin{equation*}
        \begin{split}
            &\left[\chi_\pm^N - \chi^{e_\pm}_\pm\right](a_\pm,b_\pm,c_\pm) \\
            &\coloneqq \frac{\dim M - \dim N}{\dim N - 1} \chi_{\pm}^{e^\parallel_\pm}(a^\parallel_\pm,b^\parallel_\pm,c^\parallel_\pm) -  \chi_{\pm}^{e^\parallel_\pm}(a^\perp_\pm,b^\perp_\pm,c^\parallel_\pm) - \chi_{\pm}^{e^\parallel_\pm}(a^\perp_\pm,b^\parallel_\pm,c^\perp_\pm), 
        \end{split}
    \end{equation*}
    where we have abbreviated $\chi (a,b,c) := \langle \chi (a,b),c\rangle$.
    Then $D \in \mathcal{D}^0(\genmet, \divergence; N)$.
\end{lemma}
\begin{proof}
    It is immediate that $\chi^N_\pm - \chi^{e_\pm}_\pm \in \Gamma(\Sigma_0^\pm)$, as it is antisymmetric in the last two entries, its total antisymmetrisation vanishes, and its trace $\tr[\chi^N_\pm-\chi^{e_\pm}_\pm]$ over the first and third entry is computed as
    \begin{equation*}
        \begin{split}
            \trwith{E_N}\left[\chi_\pm^N - \chi^{e_\pm}_\pm\right] &= \frac{\dim M - \dim N}{\dim N -1} \tr \chi^{e^\parallel_\pm}_{N\pm}  -\sum_{i}\underbrace{\varepsilon_i^\pm\chi_\pm^{e_\pm^\parallel}(n_i^\pm,\cdot,n_i^\pm)}_{-\langle e_\pm^\parallel, \cdot \rangle} \\
            &= (\dim M - \dim N)\langle -  e^\parallel_\pm + e^\parallel_\pm,\cdot \rangle = 0,
        \end{split}
    \end{equation*}
    where $\chi^{e^\parallel_\pm}_{N\pm}\in \Gamma (E_N^*\otimes \mathfrak{so}(E_N))$ is the tensor obtained by restriction (and projection) of  $\chi^{e^\parallel_\pm}_{\pm}$ and $(n_i^\pm)$ is an orthonormal basis of $\mathcal N_\pm$ and $\varepsilon_i^\pm=\langle n_i^\pm,n_i^\pm \rangle=-\varepsilon_i^{\mp}$.
    Hence $D \in \iota^*\mathcal{D}^0(\genmet, \divergence)$. Furthermore, the restriction satisfies
    \begin{equation*}
        \restr{\chi^N_\pm}{E_N} = \frac{\dim M - 1}{\dim N - 1} \chi_{N\pm}^{e^\parallel_\pm}
    \end{equation*}
    which proves that $\pi^\parallel D \in \mathcal{D}^0(\genmetHS,\divergence)$ is the canonical connection.
\end{proof}

\begin{lemma}\label{gausslemma2}
    The restriction of the map $D \mapsto D^N$ to $\mathcal{D}^0(\genmet,\divergence; N)$ is a surjection onto $\mathcal{D}^0(\genmetHS,\divergence_N)$.
\end{lemma}
\begin{proof}
    We only have to check that restriction onto $E_N$ gives a surjective map $\Gamma_N(\Sigma_0^\pm) \to \Sigma_{0,N}^\pm$. This, however, is immediate, as every element of $\Sigma_{0,N}^\pm$ can be linearly extended to an element of $\Gamma_N(\Sigma_0^\pm)$ by requiring it to vanish on the normal bundle $\mathcal{N}$.
\end{proof}
\section{Generalised Exterior Geometry}\label{genextgeomsection}
Let $\pi\colon E \to M$ be an exact CA with generalised semi-Riemannian metric $\genmet \cong (g, F)$ and compatible divergence operator $\divergence = \divergence^\genmet - \scalbrack{e,\cdot}$, where $e \in \Gamma(E)$.  Let $\Sigma \subset M$ be a hypersurface\footnote{All results trivially generalise to $\iota\colon \Sigma \hookrightarrow M$ an immersion of co-dimension one.} which is nondegenerate with respect to $g$. Recall that the generalised metric $\mathcal G$ 
induces a generalised metric $\mathcal H$ on the exact CA $E_\Sigma$, see Proposition~\ref{inducedsemiriemCAprop}. Denote by $n \in \Gamma(\restr{TM}{\Sigma})$ a unit normal on $\Sigma$, and furthermore $n_\pm = \sigma_\pm n$.
\subsection{Generalised Exterior Curvature} \label{genextcurvsection}
\begin{definition}
    We define the \textit{generalised exterior curvature} or \textit{generalised second fundamental form} of $\Sigma$ with respect to $D \in \mathcal{D}^0(\genmet, \divergence)$ to be
    \begin{equation*}
        \mathcal{K}^{n_\pm} \in \Gamma(E_\Sigma^* \otimes E_\Sigma^*); \qquad\quad \mathcal{K}^{n_\pm}(a,b) \coloneqq \genmet(D_a n_\pm, b) = \pm \scalbrack{D_a n_\pm, b}
    \end{equation*}
    for $a,b \in \Gamma(E_\Sigma)$. The \emph{shape tensor} $\genextEnd^{n_\pm} \coloneqq D n_\pm \in \Gamma(\End E_\Sigma)$ is the corresponding endomorphism under the isomorphism $\mathcal G : E_\Sigma\to E_\Sigma^*$. 
\end{definition}
The \enquote{mixed-type} components of the generalised second fundamental form are of particular importance, hence we denote them by
\begin{equation*}
    \mathcal{K}^\pm \coloneqq \restr{\mathcal{K}^{n_\pm}}{E_\Sigma^\mp\times E_\Sigma^\pm}
\end{equation*} 
and by $\genextEnd^\pm = \mathcal G^{-1}\mathcal{K}^\pm: E_\Sigma^\mp \to E_\Sigma^\pm$ the corresponding components of the shape tensor.
We also introduce the \textit{conormal exterior curvature} $\gennormalext^\pm \in \Gamma(E_\Sigma^*)$, which is defined as
\begin{equation*}
    \gennormalext^\pm(a) \coloneqq \genmet(D_{n_\pm - n_\mp}n_\pm,a)
\end{equation*}
Note that $\pi(n_\pm - n_\mp) = 0$, so that the right hand side does not depend on any derivatives of $n$. The conormal exterior curvature can be understood as measuring how close a generalised LC connection with divergence $\divergence$ is to projecting to a generalised LC connection with divergence $\divergence_N$ on the hypersurface $\Sigma$:
\begin{lemma}\label{restr_conn:lem}
    Let $D\in \mathcal{D}^0(\genmet, \divergence)$. Then $D \in \mathcal{D}^0(\genmet, \divergence; \Sigma)$ if and only if both conormal exterior curvatures $\gennormalext^+$ and $\gennormalext^-$ vanish.
\end{lemma}
\begin{proof} 
    Let $v_\pm$ be a section of $E_\Sigma^\pm$. We can assume without loss of generality that $E$ is a twisted generalised tangent bundle and, thus, $n_\pm = n\pm gn$ and $v_\pm = V \pm gV$, where $V\in \Gamma (TN)$. We extend $n$ and $V$ such that $\nabla_n n = \nabla_n V = 0$. Recalling that  $[n_-,n_+]\in \ker \pi$, we see that $\langle v_+, [n_-,n_+]\rangle = \frac12 [n_-,n_+](\pi v_+)=0$. This proves that $\langle D_{n_-}v_+,n_+\rangle = -\langle v_+,D_{n_-}n_+\rangle =  -\langle v_+, [n_-,n_+]\rangle=0$. As a consequence, the condition $\mathcal L^+=0$ reduces to $\langle D_{n_+}v_+,n_+\rangle =0$. Similarly, $\mathcal L^-=0$ reduces to $\langle D_{n_-}v_-,n_-\rangle =0$. Therefore $\mathcal L^+=\mathcal L^-=0$ reduces to $0= \mathrm{tr}_{\mathcal N}(Dv_+)= \mathrm{tr}_{\mathcal N}(Dv_-)$. From $\nabla_n V=0$ it follows that $\langle D_{n_\pm}^0v_\pm,n_\pm\rangle =0$, in view of the formulas for the canonical connection $D^0$ in Lemma \ref{cangenconnlemma}. We have proven that the tensor $S=D-D^0$ has $\tr_{\mathcal N} (S v)=0$ for all $v\in \Gamma (E_\Sigma)$ if and only if the normal extrinsic curvatures $\mathcal L^\pm$ of $D$ are zero. From equation \eqref{div_deviation:eq} we see that this means that $D \in \mathcal{D}^0(\genmet, \divergence; \Sigma)$ if and only if $\mathcal L^+=\mathcal L^-=0$. 
\end{proof}
\begin{corollary} 
    The conormal exterior curvature of the canonical generalised LC connection $D^0$ with metric divergence vanishes. 
\end{corollary}
\begin{proof}
    This follows from the fact  that the canonical generalised LC connection $D^0$ with metric divergence projects to the canonical generalised LC connection $D^{N0}$ with metric divergence.
\end{proof}
Remember that the space of torsion-free, metric and divergence compatible generalised connections is an affine space over the space of sections of $\Sigma_0^+ \oplus \Sigma_0^-$, in virtue of Lemma~\ref{genlcspacelemma}. That is, taking $D = D^0 + \frac{1}{d-1} (\chi^{e_+}_+ + \chi^{e_-}_-)$ from Lemma \ref{cangenconnlemma}, we can describe an arbitrary connection in $\mathcal{D}^0(\genmet, \divergence)$ as
\begin{equation}\label{dchidef}
    D^\chi \coloneqq  D + \chi = D + \chi_+ + \chi_- 
\end{equation}
where $\chi_\pm \in \Gamma(\Sigma_0^\pm)$ and $d=\dim M$. 
We need the following computational result.  
\begin{lemma}\label{genextcomputationlemma}
    Denote $\chi_\pm^\perp(a_\pm,b_\pm) \coloneqq \chi_\pm(a_\pm, b_\pm, n_\pm)$. The pure-type and mixed-type parts of the generalised second fundamental form coming from $D^\chi$ is given by 
    \begin{equation}\label{genextcomputationlemmaeq}
        \begin{split}
            \restr{\genext^{n_\pm}}{E_\Sigma^\pm \times E_\Sigma^\pm} &= k - \chi_\pm^\perp - \frac{\scalbrack{e, n_\pm}}{\dim \Sigma} h \mp \frac{i_n H}{6} \\
            \restr{\genext^{n_\pm}}{E_\Sigma^\mp \times E_\Sigma^\pm} &= k \mp \frac{i_n H}{2} \\
        \end{split}
    \end{equation}
    Herein, we employed the isometries $\sigma_\pm\colon(TM,g) \cong (E_\pm,\genmet)$ to identify $TM$ with $E_\pm$.
\end{lemma}

\begin{proof}
    We calculate for $a=a_++a_-,b=b_++b_-\in \Gamma (E_\Sigma)$, $a_\pm , b_\pm \in \Gamma (E_\Sigma^\pm)$:
    \begin{equation*}
        \begin{split}
            & \genmet\left(D^\chi_a n_\pm, b\right) \\
            &= \genmet\left(D^\chi_{a_++a_-} n_\pm, b_+ + b_-\right) \\
            &= \genmet\left(D^\chi_{a_\pm} n_\pm, b_\pm\right) + \genmet\left(D^\chi_{a_\mp} n_\pm, b_\pm\right) \\
            &= \genmet\left(D^0_{a_\pm} n_\pm + \left(\frac{\chi^{e_\pm}_\pm}{\dim M -1} +\chi_\pm\right)(a_\pm, n_\pm), b_\pm\right) + g(\nabla^\pm_{\pi a_\mp} n, \pi b_\pm)  \\
            &= g(\nabla^{\pm 1/3}_{\pi a_\pm} n, \pi b_\pm) \pm \frac{1}{\dim \Sigma} \Big[\underbrace{\scalbrack{a_\pm,n_\pm}}_{=0}\scalbrack{e_\pm,b_\pm} - \scalbrack{e_\pm, n_\pm}\scalbrack{a_\pm, b_\pm}\Big]\\
            &\qquad \pm \chi_\pm(a_\pm, n_\pm, b_\pm) + g(\nabla^\pm_{\pi a_\mp} n, \pi b_\pm) \\
            &= g(\nabla^{\pm 1/3}_{\pi a_\pm} n, \pi b_\pm) - \frac{ \scalbrack{e, n_\pm}}{\dim \Sigma}h(\pi a_\pm, \pi b_\pm) \mp \chi_\pm^\perp(a_\pm, b_\pm)+ g(\nabla^\pm_{\pi a_\mp} n, \pi b_\pm) \\
            &= g(\nabla_{\pi a} n, \pi b_\pm) \pm \frac{1}{6} \big[ H(\pi a_\pm, n, b_\pm) + 3 H(\pi a_\mp, n, \pi b_\pm) \big] \\
            &\quad \mp \chi_\pm^\perp(a_\pm, b_\pm)- \frac{\scalbrack{e, n_\pm}}{\dim \Sigma} h(\pi a_\pm, \pi b_\pm)  \\
            &= k(\pi a, \pi b_\pm) \mp \chi_\pm^\perp(a_\pm, b_\pm) - \frac{\scalbrack{e,n_\pm}}{\dim \Sigma} h(\pi a_\pm, \pi b_\pm) \\
            &\quad \mp \frac{1}{6} \big[ H(n, \pi a_\pm, \pi b_\pm) + 3 H(n, \pi a_\mp, \pi b_\pm) \big] \\
        \end{split}
    \end{equation*}
    The result follows.
\end{proof}
\begin{corollary}\label{symmetryK_pm:cor}
    The mixed-type tensors $\genext^\pm$ only depend on the choice of generalised metric $\genmet$, and the sum $\genext^+ + \genext^-$ of the mixed-type exterior curvatures is symmetric: 
    \begin{equation}\label{mixedtypeextcursymmetryeq}
        \mathcal{K}^-(a_+, b_-) = \mathcal{K}^+(b_-, a_+)
    \end{equation}
\end{corollary}
Note that the pure-type tensors depend on the choice of $\chi$. In some settings, one may wish for the generalised metric $\genmet$ and divergence operator $\divergence$ to determine all relevant geometrical objects. Then, the quantities of interest are obtained by extracting the $\chi$-independent part of the pure-type tensors, i.e.\ by dividing out the action of $\Gamma(\Sigma_0)$. The result does not depend on whether we consider $\Gamma (\Sigma_0)$ acting on $\mathcal{D}^0(\genmet, \divergence)$ or 
$\Gamma_N(\Sigma_0)$ acting on 
$\mathcal{D}^0(\genmet, \divergence; \Sigma)$, respectively.
\begin{lemma}
     Denote by $\End^0(E_\Sigma^\pm)$ the space of traceless endomorphisms on $E_\Sigma^\pm$ and by $\mathrm{Bil}^0(E_\Sigma^\pm)$ the space of 
     trace-less bilinear forms with respect to $\langle \cdot ,\cdot \rangle$. Then the following gives a well-defined surjection:
     \begin{equation*}
         \Sigma_0^\pm \longrightarrow \mathrm{Bil}^0(E_\Sigma^\pm)\cong \End^0(E_\Sigma^\pm); \qquad\quad \chi_\pm \longrightarrow \chi_\pm^\perp
     \end{equation*}
     Furthermore, the map remains surjective if we restrict it to $\Gamma_N(\Sigma_0^\pm)$.
\end{lemma}
\begin{proof}
    It is clear that the bilinear form $\chi_\pm^\perp$ defines an endomorphism $\chi_\pm^\perp$ by the identification $\chi^\perp_\pm (a,b) = \langle \chi^\perp_\pm a,b\rangle$. Antisymmetry of $\chi_\pm$ in the last two entries  implies that $\chi_\pm(\alpha_\pm, n_\pm, n_\pm) = 0$ for all $\alpha_\pm \in \Gamma(E_\pm)$. Tracelessness of $\chi_\pm$ in the first and last entry thus implies
    \begin{equation*}
        0 = \left[\tr \chi_\pm\right](n_\pm) = -\left[\trwith{E_\Sigma} \chi^\perp_\pm\right]
    \end{equation*}
    so that indeed the map is well-defined. To prove surjectivity, we show, for instance, how to construct a pre-image
    $\chi_+\in \Gamma_N(\Sigma_0^+)$ of an arbitrary element $T \in \mathrm{Bil}^0(E_\Sigma^+)$. Let $S$ be the section of  the subbundle 
    \[ (E^+_\Sigma)^*\otimes (E_\Sigma^+)^*\wedge \langle n_+,\cdot \rangle \oplus \langle n_+,\cdot \rangle\otimes {\bigwedge}^2(E_\Sigma^+)^*\subset \iota^*(E_+^*\otimes {\bigwedge}^2E_+^*)\] which has the following values
    \[ S(a,b,n_+)=T(a,b), \quad S(n_+,a,b) = T(b,a)-T(a,b)\]
    where $a,b\in \Gamma (E_\Sigma^+)$. 
    Obviously $S$ is trace-free since $T$ is and 
    \begin{eqnarray*} (\partial S)(a,b,n_+) &=& S(a,b,n_+)+S(b,n_+,a)+S(n_+,a,b)\\
    &=& T(a,b)-T(b,a) + (T(b,a)-T(a,b))=0\end{eqnarray*}
    Therefore $S$ is a section of $\iota^*\Sigma_0^+$ and we extend it to a section $\chi_+$ of $\Sigma_0^+$.
    Then $\chi_+^\perp(a,b)=S(a,b,n_+)=T(a,b)$ for all $a,b\in E_\Sigma^+$ such that $\chi$ is a pre-image of $T$.
\end{proof}
It follows that the only component of the pure-type operators independent of the choice of generalised LC connection with specified divergence operator is the trace. We define the \textit{generalised mean curvature}
\begin{equation*}
        \genmean^\pm \coloneqq \trwith{\genmetHS} \genext^{n_\pm} 
\end{equation*}
From Lemma \ref{genextcomputationlemma}, we immediately get the following result.
\begin{corollary} \label{gen_mean:cor}
    The generalised mean curvature is given by 
    \begin{equation*}
        \genmean^\pm  = \trwith{h}k - \scalbrack{e,n_\pm}
    \end{equation*}
    In particular, it is independent of the choice of sign if and only if $g(\pi e,n) = 0$.
\end{corollary}
We now derive a counterpart  of the formula $k = \frac{1}{2} L_ng$ in generalised geometry. 
\begin{lemma}
     On $E_\Sigma$, it holds 
     \begin{equation*}
         \frac{1}{2}[n_\pm, \genmet] = \genext^+ + \genext^-
     \end{equation*}
\end{lemma}
\begin{proof}
    This follows from direct computation. We only consider $n_+$, the case of $n_-$ is completely analogous. Let us consider $a,b \in \Gamma(E)$, and write $a = a_+ + a_-$, $b =b_+ + b_-$. Then
    \begin{eqnarray*}   
          [n_+, \genmet](a_\pm,b_\pm)  &=& \pi(n_+)\genmet(a_\pm,b_\pm) - \genmet([n_+,a_\pm], b_\pm) - \genmet(a_\pm, [n_+, b_\pm]) \\
            &=& \pm (\pi( n_+) \scalbrack{a_\pm,b_\pm} - \scalbrack{[n_+,a_\pm], b_\pm} - \scalbrack{a_\pm, [n_+, b_\pm]})\\ 
            &=& 0   
    \end{eqnarray*}
    Furthermore
    \begin{eqnarray*}     
            [n_+, \genmet](a_+,b_-)
            &=& \pi(n_+) \genmet(a_+,b_-) - \genmet([n_+,a_+], b_-) - \genmet(a_+, [n_+, b_-]) \\
            &=&  \scalbrack{[n_+,a_+],b_-} - \scalbrack{a_+, [n_+, b_-]}\\
            &=& 2 \scalbrack{[b_-, n_+], a_+}\\
            &=& 2 \genmet(D_{b_-} n_+, a_+) \\
            &=& 2 \genext^+(b_-, a_+)
    \end{eqnarray*}
    Symmetry of $[n_\pm, \genmet]$ together with (\ref{mixedtypeextcursymmetryeq}) implies the claim.
\end{proof}

\subsection{Generalised Gauß and Codazzi Equations}\label{gengausseqsecion}
In this section we will determine the counterpart of the Gau{\ss} and Codazzi equations in generalised geometry. Let from now on $D \in \mathcal{D}^0(\genmet, \divergence)$ and $D^\Sigma \in \mathcal{D}^0(\genmetHS)$ the inherited connection. Denote $\varepsilon = g(n,n)= \genmet(n_\pm, n_\pm)$.
Remember that, in general, the inherited divergence operator $\divergence_\Sigma$ does not agree with the divergence operator $\divergence_{D^\Sigma}$ coming from the inherited connection $D^\Sigma$, i.e. $\divergence_\Sigma \neq \divergence_{D^\Sigma}$.

The starting point is the decomposition of the ambient generalised connection in tangent and normal parts.
\begin{lemma}[Generalised Gau{\ss} and Weingarten equations]\label{finalgausslemma}
    Let $a,b \in \Gamma(E_\Sigma)$ and denote by $a_\pm,b_\pm$ their projections to $\Gamma(E_\Sigma^\pm)$. Then
    \begin{equation}\label{gausslemmaeqs}
        \begin{split}
            (D_a b_\pm)^\parallel &= D^\Sigma_a b_\pm\\
            \genmet(D_a b_\pm,n_\pm) &= - \genext^{n_\pm}(a, b_\pm) \\
            D_a n_\pm &=   \genextEnd^{n_\pm}(a)
        \end{split}
    \end{equation}
    Herein, $(\cdot)^\parallel \colon \iota^*E \to E_\Sigma$ denotes the orthogonal projection.
\end{lemma}
\begin{proof}The last equation holds by definition of the generalised shape tensor
and the first by definition of $D^\Sigma$, see Lemma \ref{gausslemma1}. We recall that the values of $D_{\tilde{a}}\tilde{b}$ along $\Sigma$ do not depend on the chosen extensions $\tilde{a}, \tilde{b}$ of $a, b$.   For that reason the tildes are omitted. 
    Finally, 
    \begin{equation*}
        \genmet(D_a b_\pm, n_\pm) = -\genmet( b_\pm, D_an_\pm)= - \genext^{n_\pm}(a, b_\pm)
    \end{equation*}   
\end{proof}
\begin{lemma}\label{nflatlemma}
    For all $u_\pm,v_\pm \in \Gamma(E^\pm_\Sigma)$ and $a \in \Gamma(E_\Sigma)$, it holds
    \begin{enumerate}[label=(\roman*)]
        \item $\quad \allowdisplaybreaks{} [n_\pm-n_\mp, a] = 0$,
        \item $\quad \allowdisplaybreaks \genmet\left( D_{n_\pm-n_\mp} u_\pm, v_\pm\right)= \genext^{n_\pm}(u_\pm,v_\pm) - \genext^{n_\pm}(v_\pm,u_\pm)$, and
        \item $\quad \allowdisplaybreaks \genmet\left( D_{n_\pm-n_\mp} u_\pm, n_\pm\right)= -\gennormalext^\pm (u_\pm )$.
    \end{enumerate}
\end{lemma}
\begin{proof} To prove (i) it suffices to pick a splitting of $E$ adapted to the normal bundle, so that $n_\pm = n \pm n^\flat$, where $n^\flat :=gn$. Then
    \begin{equation*}
        [n_\pm- n_\mp, a] = \pm 2[n^\flat, a]_H = \mp 2 i_{\pi a}\extd n^\flat = 0
    \end{equation*}
    where the last equality follows from $\Sigma \subset M$ being a submanifold. Employing this result, we find (ii): 
    \begin{equation*}
        \begin{split}
            \genmet(D_{n_\pm-n_\mp} u_\pm,  v_\pm) &= \genmet(D_{u_\pm}(n_\pm-n_\mp) + [n_\pm-n_\mp,u_\pm],v_\pm)-\genmet(D_{v_\pm}(n_\pm- n_\mp),u_\pm) \\
            &= \genext^{n_\pm}(u_\pm,v_\pm) - \genext^{n_\pm}(v_\pm,u_\pm)
        \end{split}
    \end{equation*}
    Employing metricity of $D$, we immediately find (iii):
    \begin{equation*}
        \genmet(D_{n_\pm-n_\mp} u_\pm,  n_\pm) = - \genmet(u_\pm, D_{n_\pm-n_\mp} n_\pm) = - \gennormalext^\pm(u_\pm)
    \end{equation*}
\end{proof}
\begin{lemma}\label{secondderlemma}
    The second derivative $D^2$ is given by
    \begin{equation*}
        \begin{split}
            \genmet(a, D^2_{v,w} b_\pm) &= \genmet(a, (D^{\Sigma})^2_{v,w} b_\pm)-\varepsilon\genext^{n_\pm}(v,a)\genext^{n_\pm}(w,b_\pm) \\
            &\quad + \varepsilon \sum_{\sigma \in \{ +,-\}}\genext^{n_\sigma}(v, w) \genmet(a, D_{n_\sigma} b_\pm) 
        \end{split}
    \end{equation*}
    where $a,v,w \in\Gamma(E_{\Sigma})$ and $b_\pm \in \Gamma(E_\Sigma^\pm)$.
\end{lemma}
\begin{proof}
   First, a direct computation employing Lemma \ref{finalgausslemma} reveals
    \begin{equation*}
        \begin{split}
            \genmet(a, D_v D_w b_\pm) &= \genmet\left(a, D_v \left\{D^\Sigma_w b_\pm - \varepsilon n_\pm\genext^{n_\pm}(w, b_\pm)\right\}\right) \\
            &= \genmet(a, D_v^\Sigma D_w^\Sigma b_\pm) - \varepsilon \genext^{n_\pm}(v,a) \genext^{n_\pm}(w, b_\pm)\\
        \end{split}
    \end{equation*}
    Also
    \begin{equation*}
        \begin{split}
            \genmet(a, D_{D_v w} b_\pm) &= \genmet\left(a, D_{D_v^\Sigma w}b_\pm - \varepsilon\genext^{n_+}(v, w) D_{n_+} b_\pm - \varepsilon\genext^{n_-}(v, w) D_{n_-} b_\pm\right)\\
            &= \genmet\left(a,D^\Sigma_{D_v^\Sigma w} b_\pm - \sum_{\sigma \in \{ +,-\}}\varepsilon\genext^{n_\sigma}(v, w) D_{n_\sigma} b_\pm\right) 
        \end{split}
    \end{equation*}
    The result follows.
\end{proof}
We define the curvature tensors as in \cite{jurvco2017courant} and \cite{streets2024genricciflow}, which are mathematical accounts of the definition developed in \cite{siegel1993superspace} in the context of supergravity, and described in \cite{olaf2012riemann} in the context of double field theory. Let $D \in \mathcal{D}^0(\genmet, \divergence)$. Then, the generalised Riemann tensor is given as follows:
\begin{equation}\label{genriemdefeq}
    \begin{split}
        \genriem^D(a,b,v,w) &\coloneqq \frac{1}{2} \Big \{ \scalbrack{(D^2_{v,w}-D^2_{w,v})b, a} + \scalbrack{(D^2_{b,a}-D^2_{a,b})v,w} \\
        &\quad - \trwith{E} (\scalbrack{D v, w} \scalbrack{Db,a})\Big\}\\
    \end{split}
\end{equation}
The generalised Riemann tensor is an algebraic curvature tensor, i.e.\  it lives in $\Sym^2\,\Lambda^2\,E^*$, cf.\ e.g.\ \cite[Proposition 4.10.]{jurvco2017courant}, and satisfies the Bianchi identity\linebreak[4] 
$\sum_{\sigma(u,v,w)}\genriem^D(a,u,v,w) = 0$, 
cf.\  \cite[Lemma 1]{vicenteliana} (a result based on \cite[Theorem 4.13]{jurvco2017courant}).

\begin{remark}\label{genriem_domain_rem}
    It is crucial to note that $\genriem^D$ is only non-vanishing on the pure-type subbundle
    \begin{equation*}
        \Sym^2\,\Lambda^2\,E_+^*  \oplus \Sym^2\,\Lambda^2 E_-^*, 
    \end{equation*}
    and the mixed-type subbundle
    \begin{equation*}
        \left(\Lambda^2E_+^* \oplus \Lambda^2 E_-^*\right) \vee \left(E_+^* \wedge E_-^*\right)
    \end{equation*}
    To see this, one first notes that if the pairs $(a,b)$ and $(v,w)$ are both mixed-type, the metricity of $D$ implies that $\genriem^D(a,b,v,w) = 0$. From the Bianchi identity, it follows that also for $a,b \in \Gamma(E_+)$ and $v,w \in \Gamma(E_-)$ (or vice versa) one has $\genriem^D(a,b,v,w) = 0$.
\end{remark}
From the generalised Riemann tensor, one defines the full generalised Ricci tensor as
\begin{equation*}
    \genfullric^D(u,v) \coloneqq \trwith{E} \genriem^D(\cdot,u,\cdot, v)
\end{equation*}
where we recall that the trace is taken with $\langle \cdot , \cdot \rangle$. 

We note that, while the generalised Riemann tensor generally depends on the choice of generalised connection $D \in \mathcal{D}^0(\genmet,\divergence)$, the restrictions to mixed-type entries
\begin{equation*}
    \genric^\pm \coloneqq \restr{\genfullric^D}{E_\mp \times E_\pm}
\end{equation*}
do not \cite{streets2024genricciflow}.

Finally, one obtains the generalised scalar curvature as the trace with the metric:
\begin{equation*}
    \genscal = \frac{1}{2}\trwith{\genmet} \genfullric^D
\end{equation*}
\begin{remark}\label{genriccomparisonrem}
    In \cite[Theorem 8]{rubiogenricequivalence}, different definitions of the generalised Ricci curvature one encounters in the literature are related. In particular, it is found that the mixed-type tensors $\genric^\pm$ are related to the generalised Ricci tensors $\genric_\mathrm{GF}^\pm$ as defined in \cite{fernandez, genricciflow} by 
    \begin{eqnarray}\label{genriccomparisoneq}
        2\genric^+(a_-,b_+) = 2\genric^-(b_+,a_-) &=& \genric_{\mathrm{GF}}^+(a_-,b_+) + \genric_{\mathrm{GF}}^-(b_+,a_-).
    \end{eqnarray}
    That is, the sum $\genric^+ + \genric^- \in \Gamma (E_-^*\otimes E_+^* \oplus E_+^*\otimes E_-^*)$ is precisely the symmetrisation of $\genric_{\mathrm{GF}}^++\genric_{\mathrm{GF}}^-$. The latter being symmetric, i.e. $\genric_{\mathrm{GF}}^+(a_-,b_+)=\genric_{\mathrm{GF}}^-(b_+,a_-)$, is equivalent to the divergence operator $\divergence$ and the metric $\genmet$ being \textit{compatible}, which means that the generalised vector field $e \in \Gamma(E)$ associated to the pair $(\genmet, \divergence)$ via the equation $\divergence = \divergence^\genmet - \scalbrack{e,\cdot}$ is generalised Killing, $[e,\genmet] = 0$  \cite[Lemma 3.32]{genricciflow}. In particular, this is the case if $e = 2\xi$ for a closed one-form $\xi$. 
\end{remark}
\begin{theorem}[Generalised Gau{\ss} equations]\label{riemgaußthm}
    The pure-type part of the generalised Riemann tensors obtained from $D$ and $D^\Sigma$ are related by the following Gauß equation, where $a,b,v,w \in \Gamma( E_\Sigma^\pm )$:
    \begin{equation*}
        \begin{split}
            &\pm 2 \varepsilon\left\{\genriem^D(a,b,v,w) - \genriem^{D^\Sigma}(a,b,v,w)\right\} \\
            &= \genext^{n_\pm}(w, a) \genext^{n_\pm}(v, b)  - \genext^{n_\pm}(v, a) \genext^{n_\pm}(w, b) \\
            & \quad+ \genext^{n_\pm}(a, w) \genext^{n_\pm}(b, v) - \genext^{n_\pm}(b, w) \genext^{n_\pm}(a, v) \\
            &\quad  +[\genext^{n_\pm}(v,w) - \genext^{n_\pm}(w,v)] [\genext^{n_\pm}(b, a) - \genext^{n_\pm}(a, b)] \\
        \end{split}
    \end{equation*}
    The mixed-type part satisfies with $a,v,w \in \Gamma(E_{\Sigma}^\pm), \bar{b} \in \Gamma(E_{\Sigma}^\mp)$
    \begin{equation*}
        \begin{split}
            &\pm2\varepsilon \left\{ \genriem^D(a, \bar{b}, v, w) - \genriem^{D^\Sigma}(a, \bar{b}, v, w)\right\} \\
            &=  \genext^{n_\pm}(a,w)\genext^{\pm}(\bar{b}, v)  -\genext^\pm(\bar{b}, w) \genext^{n_\pm}(a,v) \\
            &\quad +  [\genext^{n_\pm}(v, w)-\genext^{n_\pm}(w,v)] \genext^\pm(\bar{b}, a) 
        \end{split}
    \end{equation*}
\end{theorem}
\begin{proof}
    We begin by computing for general sections $a,b,v,w \in \Gamma(E_\Sigma)$ with the help of Lemma \ref{secondderlemma}
    \begin{equation*}
        \begin{split}
            &2 \genriem^D(a,b,v,w) \\
            &= \scalbrack{D^2_{v,w} b - D^2_{w,v} b, a} + \scalbrack{D^2_{b, a} v - D^2_{a, b} v, w} - \trwith{E} (\scalbrack{D v, w} \scalbrack{Db,a})\\
            &= 2 \genriem^{D^\Sigma}(a,b,v,w)\\
            &\quad + \varepsilon\sum_\sigma \Big\{ -\sigma \genext^{n_\sigma}(v, a) \genext^{n_\sigma}(w, b) +\sigma \genext^{n_\sigma}(w, a) \genext^{n_\sigma}(v, b) \\
            & \quad -\sigma  \genext^{n_\sigma}(b, w) \genext^{n_\sigma}(a, v) +\sigma \genext^{n_\sigma}(a, w) \genext^{n_\sigma}(b, v) -\sigma \scalbrack{D_{n_\sigma} v, w} \scalbrack{D_{n_\sigma} b, a} \Big\}\\
            &\quad +\sum_{\sigma}
            \big\{ (\genext^{n_\sigma}(v, w)  - \genext^{n_\sigma}(w, v)) \scalbrack{D_{n_\sigma}b, a} +(\genext^{n_\sigma}(b, a)  - \genext^{n_\sigma}(a, b)) \scalbrack{D_{n_\sigma}v, w}\big\} 
        \end{split}
    \end{equation*}
    We restrict to the pure-type case, and then make use of Lemma \ref{nflatlemma} (ii):
    \begin{equation*}
        \begin{split}
            &\pm 2 \varepsilon\left\{\genriem^D(a,b,v,w) - \genriem^{D^\Sigma}(a,b,v,w)\right\} \\
            &= \genext^{n_\pm}(w, a) \genext^{n_\pm}(v, b)  - \genext^{n_\pm}(v, a) \genext^{n_\pm}(w, b) + \genmet(D_{n_\pm}b, a) [\genext^{n_\pm}(v, w) - \genext^{n_\pm}(w, v)] \\
            & \quad+ \genext^{n_\pm}(a, w) \genext^{n_\pm}(b, v) - \genext^{n_\pm}(b, w) \genext^{n_\pm}(a, v) + \genmet(D_{n_\pm}v, w) [\genext^{n_\pm}(b, a) - \genext^{n_\pm}(a, b)] \\
            &\quad - \genmet(D_{n_\pm} v, w) \genmet(D_{n_\pm} b, a) + \genmet(D_{n_\mp} v, w) \genmet(D_{n_\mp} b, a) \\
            &= \genext^{n_\pm}(w, a) \genext^{n_\pm}(v, b)  - \genext^{n_\pm}(v, a) \genext^{n_\pm}(w, b) \\
            & \quad+ \genext^{n_\pm}(a, w) \genext^{n_\pm}(b, v) - \genext^{n_\pm}(b, w) \genext^{n_\pm}(a, v) \\
            &\quad  + \genmet(D_{n_\pm - n_\mp}v, w) [\genext^{n_\pm}(b, a) - \genext^{n_\pm}(a, b)] \\
            &= \genext^{n_\pm}(w, a) \genext^{n_\pm}(v, b)  - \genext^{n_\pm}(v, a) \genext^{n_\pm}(w, b) \\
            & \quad+ \genext^{n_\pm}(a, w) \genext^{n_\pm}(b, v) - \genext^{n_\pm}(b, w) \genext^{n_\pm}(a, v) \\
            &\quad  +[\genext^{n_\pm}(v,w) - \genext^{n_\pm}(w,v)] [\genext^{n_\pm}(b, a) - \genext^{n_\pm}(a, b)] \\
        \end{split}
    \end{equation*}
    Finally, restricting to the mixed-type case, we get with $a,v,w \in \Gamma(E_\Sigma^\pm)$ and $\bar{b} \in \Gamma( E_\Sigma^\mp)$
    \begin{equation*}
        \begin{split}
            &\pm 2\varepsilon \left\{\genriem^D(a,\bar{b},v,w) -\genriem^{D^\Sigma}(a,\bar{b},v,w)\right\} \\
            &= \genext^{n_\pm}(a, w) \genext^{\pm}(\bar{b}, v) - \genext^{\pm}(\bar{b}, w) \genext^{n_\pm}(a, v) \\
            & \quad \pm \genext^{\pm}(\bar{b}, a) \scalbrack{D_{n_\pm}v, w} \mp \genext^{\mp}(a, \bar{b}) \scalbrack{D_{n_\mp}v, w} \\
        \end{split}
    \end{equation*}
    The result follows from symmetry of $\genext^\pm + \genext^\mp$ (see Corollary \ref{symmetryK_pm:cor}) and Lemma \ref{nflatlemma} (ii). 
\end{proof}
To present the next result succinctly, we introduce the new operations  $\{ \cdot \}^{\sigmasym}$ and $\{ \cdot \}^{\sigmaantisym}$ of \enquote{$\sigma$-symmetrisation} and \enquote{$\sigma$-antisymmetrisation}. They are defined for a pair $(T^\sigma)_{\sigma = \pm}$ of mixed-type tensors $T^\pm \in \Gamma(E_\mp^*\otimes E_\pm^*)$, as follows
\begin{equation*}
    \begin{split}
        2 \{T^\sigma\}^{\sigmasym}(\bar{a},b) &= 2 [T^\pm + T^\mp]^{\sym}(\bar{a},b) = T^\pm(\bar{a},b) + T^\mp(b, \bar{a}), \\
        2 \{T^\sigma\}^{\sigmaantisym}(\bar{a},b) &= 2 [T^\pm + T^\mp]^{\antisym}(\bar{a},b) = T^\pm(\bar{a},b) - T^\mp(b, \bar{a}). \\
    \end{split}
\end{equation*}
We use curly brackets for better distinction. Note that with this operation, one can express (\ref{genriccomparisoneq}) as $\genric^\pm(\bar{a},b) = \{\genric_{\mathrm{GF}}^\sigma\}^\sigmasym(\bar{a},b)$.
\begin{corollary}\label{riccigausscor}
    The pure-type part of the generalised Ricci tensors on $E$ and $E_N$ satisfy the following Gauß equation, where $a,b \in \Gamma(E_{\Sigma}^\pm)$
    \begin{equation*}
        \begin{split}
            & \varepsilon\left\{\genfullric^D(a,b) - \genfullric^{D^\Sigma}(a,b)\right\} \\
            &= \pm \genriem^D(n_\pm,a,n_\pm, b) + [(\genext^{n_\pm})^2- \genmean^\pm \genext^{n_\pm}]^{\mathrm{sym}}(a,b) + 2 [(\genext^{n_\pm})^{\mathrm{antisym}}]^2(a,b)
        \end{split}
    \end{equation*}
    The mixed-type parts satisfy with $\bar{a} \in \Gamma(E_{\Sigma}^\mp), b \in \Gamma(E_{\Sigma}^\pm)$
    \begin{equation*}
        \begin{split}
            &\varepsilon\left\{ \genric^\pm(\bar{a}, b) -\genric_\Sigma^\pm(\bar{a}, b) \right\}\\
            &= \pm2\{\genriem^D(n_\sigma, \cdot, n_\sigma, \cdot)\}^{\sigmaantisym}(\bar{a},b) + \left\{(\genext^{n_\sigma})^2 - \genmean^\sigma \genext^\sigma\right\}^{\sigmasym}(\bar{a},b)\\
            &\quad  + 2 \left\{(\genext^{n_\sigma})^{\mathrm{antisym}}(\genextEnd^{\sigma}\cdot,\cdot)\right\}^{\sigmasym}(\bar{a},b)
        \end{split}
    \end{equation*}
where $(\genext^{n_\pm})^2(\bar{a}, b) := \genext^{n_\pm}( \genextEnd^{\pm} \bar{a}, b)$ and $\genextEnd^{\pm}$ is defined by $\mathcal G(\genextEnd^{\pm}\bar a,b)=\genext^{\pm}(\bar a, b)$. The expression $(\genext^{n_\sigma})^{\mathrm{antisym}}$ stands for the antisymmetric part of the tensor $\genext^{n_\sigma}$, $\sigma \in \{ \pm\}$.
\end{corollary}
\begin{proof}
    For the pure-type equation, we calculate from contraction of the Gauß equation for the generalised Riemann tensors, Theorem \ref{riemgaußthm},
    \begin{equation*}
        \begin{split}
            & 2 \varepsilon\left\{\genfullric^D(a,b) - \genfullric^{D^\Sigma}(a,b)\right\} \mp 2\genriem^D(n_\pm,a,n_\pm,b) \\
            &= \genext^{n_\pm}(\genextEnd^{n_\pm} b, a) - \genext^{n_\pm}(b, a) \genmean^\pm + \genext^{n_\pm}(\genextEnd^{n_\pm}a, b) - \genmean^\pm \genext^{n_\pm}(a, b)  \\
            &\quad + \sum_i \varepsilon_i[\genext^{n_\pm}(b,e_i^\pm) - \genext^{n_\pm}(e_i^\pm,b)][\genext^{n_\pm}(e_i^\pm,a) - \genext^{n_\pm}(a,e_i^\pm)] \\
            & = 2[(\genext^{n_\pm})^2- \genmean^\pm \genext^{n_\pm}]^{\mathrm{sym}}(a,b) + 4 [(\genext^{n_\pm})^{\mathrm{antisym}}]^2(a,b)
        \end{split}
    \end{equation*}
    where $\varepsilon_i=\pm \langle e_i^\pm,e_i^\pm\rangle$.
    The statement regarding the mixed-type parts follows from a similar calculation employing Theorem \ref{riemgaußthm}:
    \begin{equation}\label{riccigausscoreq1}
        \begin{split}
            &2 \varepsilon\left\{ \genric^\pm(\bar{a}, b) -\genric_\Sigma^\pm(\bar{a}, b) \right\}\\
            &= \pm  2\varepsilon \left\{ \sum_i\varepsilon_i\genriem^D(e_i^\pm,\bar{a},e_i^\pm,b) - \sum_i\varepsilon_i\genriem^{D^\Sigma} (e_i^\pm,\bar{a},e_i^\pm,b) + \varepsilon \genriem^D(n_\pm, \bar{a}, n_\pm, b) \right\} \\
            & \quad \mp 2\varepsilon \left\{ \sum_i\varepsilon_i\genriem^D(e_i^\mp,\bar{a},e_i^\mp,b) - \sum_i\varepsilon_i\genriem^{D^\Sigma} (e_i^\mp,\bar{a},e_i^\mp,b) + \varepsilon \genriem^D(n_\mp,\bar{a}, n_\mp, b) \right\} \\
            &= \pm 2 \genriem^D(n_\pm,\bar{a}, n_\pm, b) + 2 (\genext^{n_\pm})^2(\bar{a}, b) - \genmean^\pm \genext^\pm(\bar{a}, b) - \genext^{n_\pm}(b, \genextEnd^\pm \bar{a})\\
            &\quad \mp 2 \genriem^D(n_\mp,\bar{a}, n_\mp, b) + 2 (\genext^{n_\mp})^2(b, \bar{a}) - \genmean^\mp \genext^\mp(b, \bar{a}) - \genext^{n_\mp}(\bar{a}, \genextEnd^\mp b) \\
        \end{split}
    \end{equation}
    We are done. 
\end{proof}
\begin{remark}
    If we assume the pair $(\genmet, \divergence)$ to be compatible, we can simplify the result of Corollary \ref{riccigausscor} for the mixed-type case via a comparison to the generalised Ricci tensors $\genric^\pm_{\mathrm{GF}}$ as introduced in \cite{fernandez}, cf. Remark \ref{genriccomparisonrem}. Employing (\ref{genriccomparisoneq}) and $\genric^\pm_{\mathrm{GF}}(\bar{a},b) = \genric^\mp_{\mathrm{GF}}(b,\bar{a})$, we obtain
    \begin{equation*}
        \begin{split}
            2\genric^\pm(\bar{a},b) &= \genric^\pm_{\mathrm{GF}}(\bar{a},b) + \genric^\mp_{\mathrm{GF}}(b,\bar{a}) = 2 \genric^\pm_{\mathrm{GF}}(\bar{a},b) \\
            & = \trwith{E_{\pm}} \genriem^D(\cdot,\bar{a}, \cdot, b)
        \end{split}
    \end{equation*}
    so that
    \begin{equation*}
        \begin{split}
            &\varepsilon\left\{ \genric^\pm(\bar{a}, b) -\genric_\Sigma^\pm(\bar{a}, b) \right\}\\
            &= \pm2\genriem^D(n_\pm, \bar{a}, n_\pm, b)] + \left[(\genext^{n_\pm})^2 - \genmean^\pm \genext^\pm\right](\bar{a},b)  + 2 (\genext^{n_\pm})^{\mathrm{antisym}}(\genextEnd^\pm\bar{a},b).
        \end{split}
    \end{equation*}
\end{remark}
For the next Proposition, we need the following
\begin{lemma}\label{riccipuremixedcomparisonlemma}
    Let $D \in \mathcal{D}^0(\genmet, \divergence; \Sigma)$. Then
    \begin{equation*}
        \begin{split}
            &\genfullric^D(n_+,n_+) + \genfullric^D(n_-,n_-) - \genfullric^D(n_+,n_-)- \genfullric^D(n_-,n_+) \\
            &= \sum_\pm \left\{- \trwith{\genmetHS} (\genext^{n_\pm})^2 + \frac{1}{2}\absolute{\restr{\genext^{n_\pm}}{E_\pm \oplus E_\pm}}^2_\genmetHS + \frac{1}{2} \absolute{\genext^\pm}_\genmetHS^2 \right\} \\
        \end{split}
    \end{equation*}
\end{lemma}
\begin{proof}
    For this calculation, we assume $\nabla_n n = 0$. As $D \in \mathcal{D}^0(\genmet, \divergence; \Sigma)$, this implies $D_{n_\pm}n_\pm = D_{n_\mp} n_\pm = 0$. The second equation follows from the proof of Lemma \ref{restr_conn:lem}, and the first is then a consequence of its statement. 
    \begin{equation*}
        \begin{split}
            &\genfullric^D(n_+,n_+) + \genfullric^D(n_-,n_-) - \genfullric^D(n_+,n_-)- \genfullric^D(n_-,n_+) \\
            &= \sum_{i=1}^n \varepsilon_i\sum_\pm \pm \big\{\genriem^D(e_i^\pm, n_+, e_i^\pm, n_+) + \genriem^D(e_i^\pm, n_-, e_i^\pm, n_-)\\
            &\qquad - 2\genriem^D(e_i^\pm, n_+, e_i^\pm, n_-) \big\} \\
            &=\sum_{i=1}^n \varepsilon_i\sum_\pm \pm \frac{1}{2} \left\{ 2\scalbrack{D_{e_i^\pm}D_{n_\pm}n_\pm - D_{n_\pm} D_{e_i^\pm} n_\pm - D_{D_{e_i^\pm}n_\pm}n_\pm + D_{D_{n_\pm} e_i^\pm}n_\pm, e_i^\pm}\right. \\
            &\qquad - \trwith{E}\left(\scalbrack{De_i^\pm, n_\pm} \scalbrack{Dn_\pm, e_i^\pm}\right) \\
            &\qquad -  2\scalbrack{D_{e_i^\pm}D_{n_\mp}n_\pm - D_{n_\mp} D_{e_i^\pm} n_\pm - D_{D_{e_i^\pm}n_\mp}n_\pm + D_{D_{n_\mp} e_i^\pm}n_\pm, e_i^\pm} \\
            & \qquad \left.+2\: \trwith{E}\left(\scalbrack{De_i^\pm, n_-} \scalbrack{Dn_+, e_i^\pm}\right) \right\}\\
            &\stackrel{*}{=}\sum_{i=1}^n \varepsilon_i\sum_\pm \pm  \left\{\scalbrack{ - D_{n_\pm} D_{e_i^\pm} n_\pm - D_{D_{e_i^\pm}n_\pm}n_\pm + D_{D_{n_\pm} e_i^\pm}n_\pm, e_i^\pm}\right. \\
            &\qquad - \frac{1}{2}\trwith{E}\left(\scalbrack{De_i^\pm, n_\pm} \scalbrack{Dn_\pm, e_i^\pm}\right) \\
            &\qquad \left. -  \scalbrack{- D_{n_\mp} D_{e_i^\pm} n_\pm - D_{D_{e_i^\pm}n_\mp}n_\pm + D_{D_{n_\mp} e_i^\pm}n_\pm, e_i^\pm}\right\} \\
            &= \sum_{i=1}^n \varepsilon_i\sum_\pm \pm \left\{\mp n( \genext^{n_\pm}(e_i^\pm, e_i^\pm)) \pm \genext^{n_\pm}(e_i^\pm, D_{n_\pm}e_i^\pm) \mp \genext^{n_\pm}(\genextEnd^{n_\pm}(e_i^\pm),e_i^\pm) \right. \\
            &\qquad \pm \genext^{n_\pm}(D_{n_\pm}e_i^\pm, e_i^\pm) + \frac{1}{2}\trwith{E}\left( \genext^{n_\pm}(\cdot, e_i^\pm)\genext^{n_\pm}(\cdot,e_i^\pm) \right) \\
            &\qquad \left. \pm n(\genext^{n_\pm}(e_i^\pm, e_i^\pm)) \mp \genext^{n_\pm}(e_i^\pm, D_{n_\mp}e_i^\pm) \pm \genext^{n_\pm}(\genextEnd^{n_\mp}(e_i^\pm),e_i^\pm) \pm \genext^{n_\pm}(D_{n_\mp} e_i^\pm, e_i^\pm) \right\} \\
            &\stackrel{\#}{=} \sum_\pm \left\{- \trwith{\genmetHS} (\genext^{n_\pm})^2 + \frac{1}{2}\absolute{\genext^{n_\pm}}_E^2 + \tr(\genextEnd^{n_\pm} \genextEnd^{n_\mp}) \right\} \\
            &=  \sum_\pm \left\{- \trwith{\genmetHS} (\genext^{n_\pm})^2 + \frac{1}{2}\absolute{\restr{\genext^{n_\pm}}{E_\pm \oplus E_\pm}}^2_\genmetHS + \frac{1}{2} \absolute{\genext^\pm}_\genmetHS^2 \right\} \\
        \end{split}
    \end{equation*}
    In $*$, we used that in the preceding line both terms vanish, and that $D_{n_\pm}n_\pm = D_{n_\mp} n_\pm = 0$. In $\#$, we used that the first terms in the first and third line cancel, and that 
    \[ \sum_i \varepsilon_i \left\{\genext^{n_\pm}(D_{n_\pm-n_\mp}e_i^\pm, e_i^\pm) +\genext^{n_\pm}(e_i^\pm, D_{n_\pm-n_\mp}e_i^\pm)  \right\} = \pi (n_\pm -n_\mp)\mathcal{T}^{\pm}= 0 \] 
    \end{proof}
\begin{corollary}\label{genenergyconstrcor}
    The generalised scalar curvatures on $E$ and $E_\Sigma$ are related by the Gauß equation
    \begin{equation*}
        2 \genric^\pm(n_\mp, n_\pm) - \varepsilon \genscal = -\varepsilon \genscal_\Sigma  - \absolute{\genext^\pm}^2 + \frac{(\genmean^+)^2+(\genmean^-)^2}{2}
    \end{equation*}
    Herein, $\varepsilon = g(n,n)$ where $n$ is the unit normal on $\Sigma$ and $\genmet \cong (g, F)$.
\end{corollary}
\begin{proof}
    This is obtained by taking the trace of the pure-type generalised Gauß equation, cf.\ Corollary \ref{riccigausscor}:
    \begin{equation*}
        \begin{split}
            & 4 \varepsilon \left\{\genscal - \genscal_\Sigma\right\} \\
            &= \sum_\pm \left\{4 \genfullric^D(n_\pm,n_\pm) - 2 (\genmean^\pm)^2 + 4 \trwith{\genmetHS} (\genext^{n_\pm})^2 -2 \absolute{\restr{\genext^{n_\pm}}{E_\pm\oplus E_\pm}}^2_\genmetHS \right\} \\
            &= 8 \genric^\pm(n_\mp,n_\pm)- 2(\genmean^+)^2 - 2(\genmean^-)^2 + 2 \sum_\pm \absolute{\genext^{\pm}}_\genmetHS^2 \\ 
        \end{split}
    \end{equation*}
    For the last equality, we assumed $D \in \mathcal{D}^0(\genmet, \divergence; \Sigma)$ and employed Lemma \ref{riccipuremixedcomparisonlemma}.
\end{proof}
The above has, in the spirit of the ADM formulation of general relativity, the interpretation of a generalised energy constraint, if $\Sigma$ is a spacelike hypersurface in the base manifold of a generalised Lorentzian CA.  In order to compute the generalised momentum constraint, we derive the generalised Codazzi equations.
\begin{theorem}[Generalised Codazzi equations]\label{codazzithm}
    The normal components of the generalised Riemann curvature satisfy the Codazzi equations
    \begin{equation*}
        \pm2 \genriem^D(a, \bar{b}, n_\pm,w) = [D^\Sigma_{\bar{b}}\genext^{n_\pm}](a,w) - [D^\Sigma_a\genext^\pm](\bar{b},w)+ \varepsilon \genext^\pm(\bar{b},a) \gennormalext^\pm(w)
    \end{equation*}
    and
    \begin{equation*}
        \begin{split}
            &\pm2 \genriem^D(a, b, n_\pm- n_\mp,w) \\
            & = [D^\Sigma_w \genext^{n_\pm}](a,b) - [D^\Sigma_w\genext^{n_\pm}](b,a) + [D^\Sigma_{b}\genext^{n_\pm}](a,w) - [D^\Sigma_a\genext^{n_\pm}](b,w)\\
            &\quad + \varepsilon\:\big\{\gennormalext^\pm(w)[\genext^{n_\pm}(b,a)- \genext^{n_\pm}(a,b) ] +\gennormalext^\pm(b)\genext^{n_\pm}(w,a)-\gennormalext^\pm(a) \genext^{n_\pm}(w,b) \big\}
        \end{split}
    \end{equation*}
    and
    \begin{equation*}
        \begin{split}
            & \pm 2 \genriem^D(a, n_\pm - n_\mp, n_\pm - n_\mp, b) \\
            &= 2 [(\genext^{n_\pm})^2]^{\mathrm{sym}}(a,b) - 2 \genext^{n_\pm}(a, \genextEnd^{n_\pm} b) - \genext^\pm(\genextEnd^\mp a, b) \\
            &\quad + \trwith{\genmetHS} \left[\genext^{n_\pm}(\cdot,a) \genext^{n_\pm}(\cdot, b)\right]_{E_\pm \times E_\pm} - [D^\Sigma_a \gennormalext^{\pm}](b) - [D^\Sigma_b \gennormalext^{\pm}](a) 
        \end{split}
    \end{equation*}
    and
    \begin{equation*}
        \begin{split}
            &\pm 2 \genriem^D(a, n_\pm - n_\mp, n_\pm - n_\mp, \bar{b}) \\
            &= 4 \{(\genext^{n_\sigma})^2\}^{\sigmaantisym}(\bar{b},a) - 2 \{\genext^{n_\sigma}(\cdot,\genextEnd^\sigma \cdot) \}^\sigmaantisym +  [D_a^\Sigma \gennormalext^\mp](\bar{b}) - [D_{\bar{b}}^\Sigma \gennormalext^\pm](a)
        \end{split}
    \end{equation*}
    Herein, $a,b,w \in \Gamma(E_{\Sigma}^\pm)$ and $\bar{b} \in \Gamma(E_{\Sigma}^\mp)$.
\end{theorem}
\begin{proof}
    The result follows from straightforward calculation. We check that
    \begin{equation*}
        \begin{split}
            & \pm2 \genriem^D(a, \bar{b}, n_\pm,w) \\
            & = \genmet(D_{\bar{b}} D_{a} n_\pm - D_{a} D_{\bar{b}} n_\pm - D_{D_{\bar{b}} a} n_\pm + D_{D_{a} \bar{b}} n_\pm, w) \\
            & = + \genmet\left(D^\Sigma_{\bar{b}}[\genextEnd^{n_\pm}a] - D^\Sigma_a[\genextEnd^\pm \bar{b}],w\right) - \genext^{n_\pm}(D^\Sigma_{\bar{b}}a,w) + \genext^{\pm}(D^\Sigma_a \bar{b}, w) \\
            &\quad - \varepsilon\genmet(D_{\bar{b}}a,n_\pm) \genmet(D_{n_\pm}n_\pm,w) + \varepsilon\genmet(D_{a}\bar{b},n_\mp) \genmet(D_{n_\mp}n_\pm,w) \\
            & = [D^\Sigma_{\bar{b}}\genext^{n_\pm}](a,w) - [D^\Sigma_a\genext^\pm](\bar{b},w) + \varepsilon \genext^\pm(\bar{b},a) \gennormalext^\pm(w)
        \end{split}
    \end{equation*}
    where we used the symmetry of $\gennormalext^++\gennormalext^-$, see Corollary \ref{symmetryK_pm:cor}.
    In the following calculations, we assume for simplicity $\nabla_n n = 0$ and hence $D_{n_\mp}n_\pm = [n_\mp,n_\pm]= 0$.
    With the help of Lemma \ref{nflatlemma}  we compute the following.
    \begin{equation*}
        \begin{split}
            & \pm2 \genriem^D(a, b, n_\pm - n_\mp,w) \\
            & = \genmet ( D_{n_\pm-n_\mp} D_w b - D_w D_{n_\pm-n_\mp}b - D_{D_{n_\pm-n_\mp}w}b + D_{D_w (n_\pm-n_\mp)}b,a )\\
            & \quad +\genmet\left(D_{b} D_{a} n_\pm - D_{a} D_{b} n_\pm - D_{D_{b} a} n_\pm + D_{D_{a} b} n_\pm, w\right)\\
            &\quad \mp \trwith{E}(\genmet(D n_\pm, w) \genmet(Db,a)) \\
            & = [D^\Sigma_w \genext^{n_\pm}](a,b) - [D^\Sigma_w\genext^{n_\pm}](b,a) \pm \genmet(D^\Sigma_{\genmet(D_{(\cdot)}n_\pm,w)}b,a) \\
            &\quad + \varepsilon\big[\gennormalext^\pm(b)\genext^{n_\pm}(w,a)-\gennormalext^\pm(a)\genext^{n_\pm}(w,b)\big]\\
            & \quad + [D^\Sigma_{b}\genext^{n_\pm}](a,w) - [D^\Sigma_a\genext^{n_\pm}](b,w) + \varepsilon\: \gennormalext^\pm(w)[\genext^{n_\pm}(b,a) - \genext^{n_\pm}(a,b) ]\\
            &\quad\mp \trwith{E}(\genmet(D n_\pm, w) \genmet(Db,a)) \\
            & = [D^\Sigma_w \genext^{n_\pm}](a,b) - [D^\Sigma_w\genext^{n_\pm}](b,a) + [D^\Sigma_{b}\genext^{n_\pm}](a,w) - [D^\Sigma_a\genext^{n_\pm}](b,w)\\
            &\quad + \varepsilon\:\big\{\gennormalext^\pm(w)[\genext^{n_\pm}(b,a)- \genext^{n_\pm}(a,b) ] +\gennormalext^\pm(b)\genext^{n_\pm}(w,a)-\gennormalext^\pm(a) \genext^{n_\pm}(w,b) \big\}
        \end{split}
    \end{equation*}
    where $\genmet(D_{(\cdot)}n_\pm,w)\in \Gamma (E_\Sigma^*)$ is identified with an element to $\Gamma (E_\Sigma)$ by means of the scalar product $\langle \cdot , \cdot \rangle$. 
    For the next Codazzi equation, one has
    \begin{equation*}
        \begin{split}
            & \pm 2 \genriem^D(a, n_\pm - n_\mp, n_\pm - n_\mp, b) \\
            & = \genmet\left(D_{n_\pm - n_\mp} D_b n_\pm - D_b D_{n_\pm - n_\mp} n_\pm - D_{D_{n_\pm - n_\mp}b} n_\pm+ D_{D_b (n_\pm - n_\mp)}n_\pm, a \right) \\
            &\quad + \genmet\left(D_{n_\pm - n_\mp} D_a n_\pm - D_a D_{n_\pm - n_\mp} n_\pm - D_{D_{n_\pm - n_\mp}a} n_\pm+ D_{D_a (n_\pm - n_\mp)}n_\pm, b \right) \\
            & \quad \mp\trwith{E}(\genmet(Dn_\pm,b) \genmet(Dn_\pm,a)) 
        \end{split}
    \end{equation*}
    We calculate
    \begin{equation*}
        \begin{split}
            & \genmet\left(D_{n_\pm - n_\mp} D_b n_\pm - D_b D_{n_\pm - n_\mp} n_\pm - D_{D_{n_\pm - n_\mp}b} n_\pm+ D_{D_b (n_\pm - n_\mp)}n_\pm, a \right) \\
            & = \genext^{n_\pm}(\genextEnd^{n_\pm} b,a) - \genext^{n_\pm}(a, \genextEnd^{n_\pm} b) - [D_b^\Sigma \gennormalext](a) \pm \trwith{E} \genmet(D n_\pm,b) \genmet(Dn_\pm,a)
        \end{split}
    \end{equation*}
    Therefore
    \begin{equation*}
        \begin{split}
            & \pm 2 \genriem^D(a, n_\pm - n_\mp, n_\pm - n_\mp, b) \\
            &= \genext^{n_\pm}(\genextEnd^{n_\pm}a,b) + \genext^{n_\pm}(\genextEnd^{n_\pm}b,a) - \genext^{n_\pm}(a, \genextEnd^{n_\pm} b) - \genext^{n_\pm}(b, \genextEnd^{n_\pm} a) \\
            &\quad - [D_a^\Sigma \gennormalext^\pm](b) - [D_b^\Sigma \gennormalext^\pm](a) \pm \trwith{E} \genext^{n_\pm}(\cdot,a) \genext^{n_\pm}(\cdot, b)\\ 
            &= 2 [(\genext^{n_\pm})^2]^{\mathrm{sym}}(a,b) - 2 \genext^{n_\pm}(a, \genextEnd^{n_\pm} b)\\  
            &\quad - [D_a^\Sigma \gennormalext^\pm](b) - [D_b^\Sigma \gennormalext^\pm](a) \pm \left[  \trwith{E_\pm} \genext^{n_\pm}(\cdot,a) \genext^{n_\pm}(\cdot, b)+\trwith{E_\mp} \genext^{n_\pm}(\cdot,a) \genext^{n_\pm}(\cdot, b)\right]\\
        \end{split}
    \end{equation*}
    where, using Corollary \ref{symmetryK_pm:cor}, 
    \[ \pm \trwith{E_\mp} \genext^{n_\pm}(\cdot,a) \genext^{n_\pm}(\cdot, b)=-\trwith{\mathcal H} \genext^{\pm}(\cdot,a) \genext^{\pm}(\cdot, b)= -\trwith{\mathcal H} \genext^{\mp}(a,\cdot ) \genext^{\pm}(\cdot, b)\]
    Analogously, we do the last Codazzi equation by calculating
     \begin{equation*}
        \begin{split}
            & \pm 2 \genriem^D(a, n_\pm - n_\mp, n_\pm - n_\mp, \bar{b}) \\
            & = \genmet\left(D_{n_\pm - n_\mp} D_{\bar{b}} n_\pm - D_{\bar{b}} D_{n_\pm - n_\mp} n_\pm - D_{D_{n_\pm - n_\mp}\bar{b}} n_\pm+ D_{D_{\bar{b}} (n_\pm - n_\mp)}n_\pm, a \right) \\
            &\quad +\genmet\left(D_{n_\pm - n_\mp} D_a n_\mp - D_a D_{n_\pm - n_\mp} n_\mp - D_{D_{n_\pm - n_\mp}a} n_\mp+ D_{D_a (n_\pm - n_\mp)}n_\mp, \bar{b} \right) \\
            & \quad \mp\trwith{E}(\genmet(Dn_\mp,\bar{b}) \genmet(Dn_\pm,a)) 
        \end{split}
    \end{equation*}
    We note that
    \begin{equation*}
        \begin{split}
            & \genmet\left(D_{n_\pm - n_\mp} D_{\bar{b}} n_\pm - D_{\bar{b}} D_{n_\pm - n_\mp} n_\pm - D_{D_{n_\pm - n_\mp}\bar{b}} n_\pm+ D_{D_{\bar{b}} (n_\pm - n_\mp)}n_\pm, a \right) \\
            & = \genext^{n_\pm}(\genextEnd^{\pm} \bar{b},a) - \genext^{n_\pm}(a, \genextEnd^{\pm} \bar{b}) - [D_{\bar{b}}^\Sigma \gennormalext^\pm](a) \pm \trwith{E} \genmet(D n_\mp,\bar{b}) \genmet(Dn_\pm,a)
        \end{split}
    \end{equation*}
    (cf.\ Lemma \ref{nflatlemma}) and obtain 
    \begin{equation*}
        \begin{split}
            & \pm 2 \genriem^D(a, n_\pm - n_\mp, n_\pm - n_\mp, \bar{b}) \\
            &= \genext^{n_\pm}(\genextEnd^{\pm}\bar{b},a) - \genext^{n_\pm}(a, \genextEnd^{\pm} \bar{b}) -\genext^{n_\mp}(\genextEnd^{\mp}a,\bar{b}) + \genext^{n_\mp}(\bar{b}, \genextEnd^{\mp} a) \\
            &\quad + [D_a^\Sigma \gennormalext^\mp](\bar{b}) - [D_{\bar{b}}^\Sigma \gennormalext^\pm](a) \pm \trwith{E} \genext^{n_\pm}(\cdot,a) \genext^{n_\mp}(\cdot, \bar{b})\\ 
            &=4 \{(\genext^{n_\sigma})^2\}^{\sigmaantisym}(\bar{b},a) - 2 \{\genext^{n_\sigma}(\cdot,\genextEnd^\sigma \cdot) \}^\sigmaantisym +  [D_a^\Sigma \gennormalext^\mp](\bar{b}) - [D_{\bar{b}}^\Sigma \gennormalext^\pm](a)\\ 
        \end{split}
    \end{equation*}
    where we used that
    \begin{equation*}
        \begin{split}
            & \pm \trwith{E} \genext^{n_\pm}(\cdot,a) \genext^{n_\mp}(\cdot, \bar{b}) \\
            & = \pm \left[  \trwith{E_\pm} \genext^{n_\pm}(\cdot,a) \genext^{n_\mp}(\cdot, \bar{b})+\trwith{E_\mp} \genext^{n_\pm}(\cdot,a) \genext^{n_\mp}(\cdot, \bar{b})\right] \\
            &= \genext^{n_\pm}(\genextEnd^\pm \bar{b}, a) - \genext^{n_\mp}(\genextEnd^\mp a,\bar{b})
        \end{split}
    \end{equation*}
    This proves the claim.
\end{proof}
As a Corollary, we obtain the generalised momentum constraint.
\begin{corollary}\label{genmomentumconstrcor}
    The mixed components of the generalised Ricci tensor on $E$ are given as follows:
    \begin{equation*}
        \begin{split}
            \genric^\pm(a_\mp, n_\pm) &= (\divergence_\Sigma \genextEnd^\pm)(a_\mp) - \pi a_\mp(\genmean^\pm) \\
        \end{split}
    \end{equation*}
    Herein, $a_\mp \in \Gamma(E_{\Sigma}^\mp)$ arbitrary, $\divergence_\Sigma \genextEnd^\pm$ is the divergence of the mixed-type tensor $\genextEnd^\pm$ with respect to $(\genmet, \divergence)$, as in Lemma \ref{divgentensorlemma}.
\end{corollary}
\begin{proof}
    Take the trace over $a$ and $w$ in the mixed-type equation of Theorem \ref{codazzithm}, choosing $D \in \mathcal{D}^0(\genmet, \divergence; \Sigma)$ so that $\divergence_\Sigma = \divergence_{D^\Sigma}$. Note that for such $D$ the conormal exterior curvature $\gennormalext^\pm$ vanishes.
\end{proof}
Translating this equation into the non-generalised language, one can employ the following Lemma to re-express the constraint equations for the generalised Einstein tensor $\genric^+ + \genric^--\frac12\genscal\cdot\mathcal  G$ in terms of ordinary (rather than ``generalised'') geometric objects on $\Sigma$.
\begin{lemma}\label{shapetensordivergencelemma}
    It holds
    \begin{equation*}
        (\divergence_\Sigma \genextEnd^\pm)(a_\mp) = \divergence^{e^\pm}\left(k \mp \frac{\iota_n H}{2}\right)(\pi a_\mp) - \frac14 H^2(n, \pi a_\mp)
    \end{equation*}
     where $a_\mp \in  \Gamma(E_{\Sigma}^\mp)$ is arbitrary, and $\divergence^{e^\pm} \colon \Gamma(T\Sigma)\stackrel{h}{\cong} \Gamma(T^*\Sigma) \to C^\infty(\Sigma)$ is defined by $\divergence^{e^\pm}(X) = \divergence_h(X) \mp g(X, \pi e^\pm)$.
\end{lemma}
\begin{proof}
    We denote by $D^\Sigma \in \mathcal{D}^0(\genmetHS, \divergence_\Sigma)$ the generalised LC connection from Lemma~\ref{cangenconnlemma}. We furthermore denote by $\nabla^{\Sigma\pm}$ the connections on $\Sigma$ with torsion $\pm \iota^*H$, cf. (\ref{nablapmdefeq}). We choose with $e_i$ an orthonormal frame for $T\Sigma$ and put $\varepsilon_i = g(e_i,e_i)$. Then, we calculate with Lemma~\ref{divgentensorlemma}
    \begin{equation*}
        \begin{split}
            &(\divergence_\Sigma \genextEnd^\pm)(a_\mp) \\
            & = (\divergence^{\genmetHS,\pm} \genextEnd^\pm)(a_\mp) - \scalbrack{e,\genextEnd^\pm}(a_\mp) \\
            &= \sum_i  \varepsilon_i\left(\nabla^{\Sigma\mp}_{e_i}\left[k \mp \frac{\iota_n H}{2}\right]\right)(\pi a_\mp,e_i) \mp \left[k \mp \frac{\iota_n H}{2}\right](\pi a_\mp,\pi e^\pm) \\
            &= \divergence^{e^\pm} \left[k \mp \frac{\iota_n H}{2}\right](\pi a_\mp) \pm \frac{1}{2} \sum_i \varepsilon_i \left[k\mp \frac{\iota_n H}{2}\right](H(e_i,\pi a_\mp),e_i)\\
            & = \divergence^{e^\pm}\left[k \mp \frac{\iota_n H}{2} \right](\pi a_\mp)-\frac14 H^2(n, \pi a_\mp)
        \end{split}
    \end{equation*}
\end{proof}

\begin{corollary}
    Assume that on the ambient exact CA $E\to M$, the generalised Einstein equations are satisfied, i.e. $\genric = 0$ and $\genscal = 0$. If $e = 2\xi$, $\extd\xi = 0$, the constraint equations for the generalised Einstein equations are given by
    \begin{equation*}
        \begin{split}
            -\varepsilon\rscal_\Sigma + (\tr k)^2 - \absolute{k}^2 &=-\varepsilon \frac{\absolute{H^\parallel}^2}{12} + \frac{\absolute{H^\perp}^2}{4}  -2\varepsilon (\extd^\Sigma)^*\xi^\parallel + 2 (\tr k) x  -\varepsilon \absolute{\xi^\parallel}^2-x^2 \\
            \divergence_h\: k - \extd^\Sigma \tr k &= \frac{1}{4}\scalbrack{H^\perp, H^\parallel} - \extd^\Sigma x + i_{\xi^\parallel} k \\
            0  &=\left((\extd^\Sigma)^* +i_{\xi^\parallel}\right) H^\perp \\
        \end{split}
    \end{equation*}
    Herein, we decomposed on $\Sigma$ 
    \begin{equation*}
        H = H^\parallel + \varepsilon\: n^\flat \wedge H^\perp, \qquad\quad \xi = \xi^\parallel +  \varepsilon\: x\: n^\flat
    \end{equation*}
     and denoted by $\scalbrack{H^\perp, H^\parallel}$ the one-form obtained from contracting $H^\perp$ with $H^\parallel$ (matching corresponding entries) by means of the induced metric $h$. The first equation is equivalent to the generalised energy constraint, and the other two are equivalent to the generalised momentum constraint. (Note that the case of a space-like hypersurface $(\Sigma , h=g|_\Sigma)$ in a Lorentzian manifold $(M,g)$ corresponds to $\varepsilon=-1$.)
\end{corollary}
\begin{proof}
    Recall from Corollary \ref{genenergyconstrcor} that the generalised energy constraint is given by 
    \begin{equation*}
        0 = -\varepsilon \genscal_\Sigma  - \absolute{\genext^\pm}^2 + \frac{(\genmean^+)^2+(\genmean^-)^2}{2}
    \end{equation*}
    The generalised scalar curvature for $\genscal_\Sigma$ is given by, cf. (\ref{genscaleq}),
    \begin{equation*}
        \genscal_\Sigma = \rscal_\Sigma - \frac{\absolute{H^\parallel}^2}{12} - 2 (\extd^\Sigma)^*\xi^\parallel - \absolute{\xi^\parallel}^2
    \end{equation*}
    where we have used that 
    \begin{equation*}
        \frac{1}{2}\absolute{2\xi^\parallel}_\genmetHS^2 = 2 \genmetHS(\xi^\parallel, \xi^\parallel) = 2 \scalbrack{\genmetHS \xi^\parallel, \xi^\parallel} = h^{-1}(\xi^\parallel ,\xi^\parallel) = \absolute{\xi^\parallel}_h^2
    \end{equation*} 
    Furthermore, by Lemma \ref{genextcomputationlemma} and Corollary \ref{gen_mean:cor}, we have
    \begin{equation*}
        \absolute{\genext^\pm}_\genmetHS^2 = \absolute{k}^2 + \frac{\absolute{H^\perp}^2}{4}, \qquad\qquad (\genmean^\pm)^2 = (\tr k - x)^2
    \end{equation*}
    Equivalence of the first equation to the generalised energy constraint follows from insertion. 
    
    Now recall the generalised momentum constraint
    \begin{equation*}
        0 = (\divergence_\Sigma \genextEnd^\pm)(a_\mp) - \pi a_\mp(\genmean^\pm) 
    \end{equation*}
    By Lemma \ref{shapetensordivergencelemma}, we have
    \begin{equation*}
        (\divergence_\Sigma \genextEnd^\pm)(a_\mp) = (\divergence_h - i_{\xi^\parallel})\left(k \mp \frac{H^\perp}{2}\right)(\pi a_\mp) -\frac{1}{4} \scalbrack{H^\perp, H^\parallel}(\pi a_\mp)
    \end{equation*}
    where we have used that $\pi (e_\pm) = \pm h^{-1}\xi^\parallel$ and $i_{\xi^\parallel}:=i_{h^{-1}\xi^\parallel}$. 
    The result follows from considering the sum and the difference of the cases \enquote{$+$} and \enquote{$-$}.
\end{proof}

\section{Flat Semi-Riemannian Courant Algebroids}\label{flatCAsection}

In this section we study exact semi-Riemannian  Courant algebroids for which the canonical Levi-Civita generalised connection is flat. We may call them canonically flat. For Riemannian and Lorentzian signature we show that canonical flatness implies complete triviality, by which we mean that the Courant algebroid is untwisted ($H=0$), the divergence operator coincides with the metric divergence (constant dilaton in physics terminology) and the underlying Riemannian or Lorentzian metric is flat. We begin with two simple consequences of the generalised Einstein  equations.

In this section, we always assume the exact Courant algebroid to be equipped with the semi-Riemannian generalised metric $\genmet$ and a divergence operator $\divergence$, which we express as $\divergence = \divergence^\genmet - \scalbrack{e, \cdot}$ for some $e \in \Gamma(E)$. Employing the splitting induced by the generalised metric, we also write $e = 2(X + \xi) \in \Gamma(\mathbb{T}M)$. 
\begin{lemma}\label{dilatoneomlemma}
    The generalised Einstein equations $\genric = 0$ and $\genscal = 0$ imply
    \begin{equation}\label{dilatoneomeq}
        0 = \frac{\absolute{H}^2}{6} - \extd^* \xi - \frac{1}{2}\absolute{e}_\genmet^2
    \end{equation}
    We will sometimes refer to this as the \enquote{dilaton's equation of motion}.
\end{lemma}
\begin{proof}
    We compute the trace of the expression for the generalised Ricci curvature (\ref{genriceq}) as
    \begin{equation}\label{mixedgenrictraceeq}
        \trwith{g} [\genric \circ (\sigma_\pm, \sigma_\mp)] = \rscal - \frac{\absolute{H}^2}{4} - \extd^* \xi
    \end{equation}
    The dilaton's equation of motion is obtained by subtracting the trace (\ref{mixedgenrictraceeq}) from the generalised scalar curvature (\ref{genscaleq}).
\end{proof}
\begin{corollary}\label{geneinsteinriemtrivialitycor}
    Let $E \to M$ be an exact CA over a closed and orientable  manifold $M$. Let $(\genmet, \divergence)$ be a pair consisting of a Riemannian generalised metric $\genmet$ and an exact divergence operator $\divergence = \divergence^\genmet - 2\scalbrack{\xi, \cdot}$, $\xi \in \Gamma(T^*M)$, $\xi =\extd \phi$. 
    
    Then, generalised Ricci and scalar flatness, i.e. $\genric = 0$ and $\genscal = 0$, imply $\ric = 0$, $H= 0$ and $\xi = 0$. In particular, generalised Riemann flatness implies complete triviality, i.e. $\riem = 0$, $H= 0$, and $\xi = 0$.
\end{corollary}
\begin{proof}
    We follow the argument as given in \cite[Theorem  3.50.]{genricciflow}. The dilaton's equation of motion is for $e = 2 \xi$ given as
    \begin{equation}\label{geneinsteinriemtrivialityeq1}
        0 = \frac{\absolute{H}^2}{6} - \extd^*\xi - \absolute{\xi}^2
    \end{equation}
    Integrating (\ref{geneinsteinriemtrivialityeq1}) with the volume form $\mu = e^{-2\phi}\vol_g$, one gets
    \begin{equation*}
        0 = \int_M \left(\frac{\absolute{H}^2}{6} - \extd^*\extd \phi - \absolute{\extd \phi}^2 \right) \mu = \int_M \left(\frac{\absolute{H}^2}{6} +\absolute{\extd \phi}^2 \right) \mu
    \end{equation*}
    The result follows.
\end{proof}
We want to show that generalised Riemann flatness implies complete triviality without assuming  exactness of the divergence operator nor compactness and orientability. The following result shows conformal flatness, independent of signature and topology.
\begin{theorem}\label{genriemflatnessthm}
    Let $D \in \mathcal{D}^0(\genmet, \divergence)$ be the canonical connection, denote $\divergence = \divergence^\genmet - \scalbrack{e,\cdot}$. Then $\genriem^D = 0$ implies that the Weyl tensor of $(M,g)$ vanishes. In particular, if $d = \dim M \geq 4$, $(M,g)$ is conformally flat. 
    
    More precisely, the quadratic components of the Riemann tensor are given for perpendicular vector fields $A, B \in \Gamma(TM)$, $g(A,B) = 0$, as
    \begin{equation}\label{genriemflatnessthm_riem_eq}
         \begin{split}
             &\riem(A,B,A,B) = \frac{\pm \divergence^\genmet(e_\pm)}{2d(d-1)}  \absolute{A}^2\absolute{B}^2  \\
             &\qquad+ \frac{3}{4(d-1)^2} \left\{2 \absolute{e_\pm}^2 \absolute{A}^2 \absolute{B}^2 - \absolute{B}^2 g(A,\pi e_\pm)^2 - \absolute{A}^2 g(B,\pi e_\pm)^2 \right\} \\
         \end{split}
    \end{equation}
    Furthermore, it holds $\divergence^\genmet(e_+) = - \divergence^\genmet(e_-)$, $\pi e_+ = \varepsilon \: \pi e_-$ for some $\varepsilon \in \{+1,-1\}$, and
    \begin{equation}\label{genriemflatnessthmeq}
        \nabla H = 0, \qquad\quad  g \nabla^\pm e_\pm = \frac{\divergence^\genmet(e_\pm)}{d} g, \qquad\quad 0 = \frac{\absolute{H}^2}{6} + \divergence^\genmet(e_+) - \absolute{e_+}^2.
    \end{equation}
\end{theorem}

\begin{remark}\label{genriemflatness_jet_rem}
    The proof of the above statement does not rely on working on an exact Courant algebroid. Instead, one can consider the following setting.

    Let $p\in M$, $\nabla$ a torsion-free connection\footnote{A connection at $p$ is a bilinear map 
    $\nabla : T_pM \times \Gamma (TM) \to T_pM, (v, X) \mapsto \nabla_vX$ satisfying the Leibniz rule: $\nabla_v (fX)=
    v(f)X_p + f(p)\nabla_vX$, $f\in C^\infty (M)$. Note that 
    $\nabla_vX$ depends only on the $1$-jet of $X$ at $p$. We will write $\nabla_vX|_p$ to emphasize that $\nabla_vX$ is defined at $p$.} at $p$, and $g\in J^1_p(M, \Sym^2(M))$ such that $g_p$ is non-degenerate and $\nabla g|_p = 0$. Let furthermore $\riem \in \Sym^2(\Lambda^2 T_p^*M)$, $H \in J^1_p(M, \Lambda^3 T^*M)$, and $e = X + \xi \in J^1_p(M, TM \oplus T^*M)$. Assume that at $p$ the equations in Propositions \ref{puretyperiemprop} and \ref{mixedtyperiemprop} hold with the generalised Riemann tensor set to zero.\footnote{Note that these equations immediately imply that $\riem$ satisfies the Bianchi identity and hence is an algebraic curvature tensor.} Then, also equations (\ref{genriemflatnessthm_riem_eq}) and (\ref{genriemflatnessthmeq}) hold at $p$.
\end{remark}

\begin{proof}
    Let us start by considering the following components of the mixed-type part of the generalised Riemann tensor. Take $A,V,W \in \Gamma(TM)$, denote $a = \sigma_\pm A, \Bar{a} = \sigma_\mp A, v = \sigma_\pm V$, and $w = \sigma_\pm W$. By Theorem \ref{mixedtyperiemprop} we obtain from $\genriem^D = 0$:
    \begin{equation}\label{genriemflatnesseq1}
        \begin{split}
            0 &= 2 \genriem^D(a, \Bar{a},v, w) \\
            &= - \frac{1}{3} [\nabla_A H](A,V,W) + \frac{1}{d-1} [D^0_{\Bar{a}}\chi_\pm^{e_\pm}](a,v,w)
        \end{split}
    \end{equation}
    We calculate
    \begin{equation*}
        \begin{split}
            [D^0_{\Bar{a}}\chi_\pm^{e_\pm}](a, v,  w) = g(A, V) g(\nabla^\pm_A e_\pm, W) - g(A, W) g(\nabla^\pm_A e_\pm, V)
        \end{split}
    \end{equation*}
    Considering (\ref{genriemflatnesseq1}) with $V = A \perp W$, we can conclude $0 = [D^0_{\Bar{a}}\chi_\pm^{e_\pm}](a, a,  w)$ and thus $g(\nabla^\pm_A e_\pm, W) = 0$. In particular the tensor $g\nabla^\pm e_\pm$ is symmetric and  also $0 = [D^0_{\Bar{a}}\chi_\pm^{e_\pm}](a, v,  w)$ for $V \perp A$. It follows that 
    \begin{equation*}
        0 = [D^0_{\Bar{a}}\chi_\pm^{e_\pm}](a, v,  w) \qquad\quad \text{ for all } A,V,W \in \Gamma(TM)
    \end{equation*}
    Inserting this back into (\ref{genriemflatnesseq1}), we obtain that $[\nabla_A H](A,V,W) = 0$. By a polarisation argument, it follows that $[\nabla_A H](B,V,W) = -[\nabla_B H](A,V,W)$ for all $A,B,V,W \in \Gamma(TM)$. We conclude that for all $A,B,V,W \in \Gamma(TM)$
    \begin{equation*}
        \begin{split}
            0 &= [\extd H](A,B,V,W) \\
            &= [\nabla_A H](B,V,W) - [\nabla_B H](A,V,W) + [\nabla_V H](A,B,W) - [\nabla_W H](A,B,V) \\
            &= 4 [\nabla_A H](B,V,W)
        \end{split}
    \end{equation*}
    i.e.\ $\nabla H = 0$. With this, we see that
    \begin{equation}\label{genriemflatnesseq2}
        \begin{split}
            0 &= 2 (d-1) \genriem^D(a,\Bar{b},v,w) - 2 (d-1)\genriem^D(w, \Bar{v}, b, a) \\
            &= [D^0_{\Bar{b}} \chi_\pm^{e_\pm}](a,v,w) - [D^0_{\Bar{v}} \chi_\pm^{e_\pm}](w,b,a)
        \end{split}
    \end{equation}
    since $\riem, H^{(2)}\in \Gamma (\mathrm{Sym}^2{\bigwedge}^2T^*M)$. Taking the trace of this equation over $A$,$V$ and assuming $B= W$, we obtain $d\cdot g(\nabla_B e_\pm, B) = \absolute{B}^2 \divergence^\genmet(e_\pm)$, implying
    \begin{equation*}
        g \nabla^\pm e_\pm =  \frac{\divergence^\genmet(e_\pm)}{d} g
    \end{equation*}
since the tensor $g\nabla^\pm e_\pm$ is symmetric, as shown above. Therefore
    \begin{equation*}
        [D^0_{\Bar{b}} \chi_\pm^{e_\pm}](a,v,w) = \frac{\divergence^\genmet(e_\pm)}{d}\left[g(a,v)g(b,w) - g(a,w)g(b,v) \right]  
    \end{equation*}
    To conclude the investigation of the mixed-type generalised Riemann tensor, we compute its quadratic components as
    \begin{equation}\label{genriemflatnesseq3}
        \begin{split}
            \pm 2 \genriem^D(a,\Bar{b},a,b) &= \riem(A,B,A,B) - \frac{1}{4} H^{(2)}(A,B,A,B) \\
            &\quad \pm \frac{\divergence^\genmet(e_\pm)}{d(d-1)} \mathrm{CS}(A,B)
        \end{split}
    \end{equation}
    where we introduced
    \begin{equation*}
        \mathrm{CS}(A,B) \coloneqq \absolute{A}^2\absolute{B}^2 - g(A,B)^2.
    \end{equation*}
    We can conclude from comparing the cases \enquote{$+$} and \enquote{$-$} in (\ref{genriemflatnesseq3}) that $\genriem^D = 0$ implies
    \begin{equation}\label{divpmv2:eq}
        \divergence^\genmet(e_+) = - \divergence^\genmet(e_-)
    \end{equation}
    Now, we move our attention to the pure-type tensor. To that end, note that
    \begin{equation*}
        \begin{split}
            &[D^0_b \chi_\pm^{e_\pm} - D^0_{\Bar{b}} \chi_\pm^{e_\pm}](a,v,w) \\
            &= g(a,v) g([\bismutthird - \bismut]_B e_\pm, w) - g(a,w) g([\bismutthird - \bismut]_B e_\pm, v) \\
            &= \pm \frac{1}{3}[g(a,w)H(b,e_\pm,v) - g(a,v)H(b,e_\pm,w) ]
        \end{split}
    \end{equation*}
    and hence
    \begin{equation*}
        \begin{split}
            [D^0_b \chi_\pm^{e_\pm}](a,v,w) &= \frac{\divergence^\genmet(e_\pm)}{d}\left[g(a,v)g(b,w) - g(a,w)g(b,v) \right]\\
            &\quad \pm \frac{1}{3}[g(a,w)H(b,e_\pm,v) - g(a,v)H(b,e_\pm,w) ]
        \end{split}
    \end{equation*}
    so that
    \begin{equation*}
        \begin{split}
            &[D^0_v \chi_\pm^{e_\pm}](w,b,a) - [D^0_w \chi_\pm^{e_\pm}](v,b,a)  + [D^0_b \chi_\pm^{e_\pm}](a,v,w) - [D^0_a \chi_\pm^{e_\pm}](b,v,w) \\
            &= \frac{\divergence^\genmet(e_\pm)}{d} \big\{g(w,b)g(v,a) - g(w,a)g(v,b) - g(v,b)g(w,a) + g(v,a)g(w,b)\\
            &\qquad\quad  + g(a,v)g(b,w) - g(a,w)g(b,v) - g(b,v)g(a,w) + g(b,w)g(a,v) \big\}\\
            &\quad \pm \frac{1}{3}\big\{g(w,a)H(v,e_\pm,b) - g(w,b)H(v,e_\pm,a) \\
            &\qquad\quad - g(v,a)H(w,e_\pm,b) + g(v,b)H(w,e_\pm,a)  \\
            &\qquad\quad + g(a,w)H(b,e_\pm,v) - g(a,v)H(b,e_\pm,w) \\
            &\qquad\quad- g(b,w)H(a,e_\pm,v) + g(b,v)H(a,e_\pm,w)\big\} \\
            &=  \frac{4\:\divergence^\genmet(e_\pm)}{d} \big[g(w,b)g(v,a) - g(w,a)g(v,b) \big]\\
        \end{split}
    \end{equation*}
    Now, with this and Proposition \ref{puretyperiemprop}, we compute the quadratic terms of the pure-type tensor
    \begin{equation}\label{genriemflatnesseq4}
        \begin{split}
            \pm \genriem^D(a,b,a,b) &= \riem(A,B,A,B) - \frac{1}{12} H^{(2)}(A,B,A,B) \\
            &\quad \pm \frac{2\:\divergence^\genmet(e_\pm)}{d(d-1)} \mathrm{CS}(A,B) - Q(A,B) \\
        \end{split}
    \end{equation}
    where we defined
    \begin{equation}\label{Qdefeq}
        \begin{split}
            Q(A,B) &\coloneqq \frac{1}{2(d-1)^2} \left\{2 \absolute{e_\pm}^2 \mathrm{CS}(A,B)- \absolute{B}^2 g(A,e_\pm)^2 - \absolute{A}^2 g(B,e_\pm)^2\right.\\ 
            & \qquad\qquad\qquad \left. + 2 g(A,B) g(A,e_\pm) g(B,e_\pm)\right\}.
        \end{split}
    \end{equation}
    Note that, by (\ref{divpmv2:eq}) and $\genriem^D=0$, $Q$ has to be independent of the choice of sign, and thus $\pi e_+ = \varepsilon\: \pi e_-$ for some  $\varepsilon \in \{+1,-1\}$. 
    More precisely, 
    \begin{equation} \label{epm:eq}
        \begin{cases}
            e=2X \;\; \mbox{and}\;\; e_\pm = X\pm gX, &\mbox{if}\;\;  \varepsilon=+1\\
            e=2\xi\;\;\mbox{and}\;\; e_\pm = \pm g^{-1}\xi \pm\xi, &\mbox{if}\;\;   \varepsilon=-1
        \end{cases}
    \end{equation}
    Comparing the quadratic expressions (\ref{genriemflatnesseq3}) and (\ref{genriemflatnesseq4}), we obtain that
    \begin{equation}\label{genriemflatnesseq5}
        \begin{split}
            &\frac{1}{6}H^{(2)}(A,B,A,B) = \mp \frac{\:\divergence^\genmet(e_\pm)}{d(d-1)} \mathrm{CS}(A,B) + Q(A,B)\\
        \end{split}
    \end{equation}
    so that $\genriem^D = 0$ implies
    \begin{equation*}
        \begin{split}
            &2 \riem(A,B,A,B) = \frac{1}{2} H^{(2)}(A,B,A,B) \mp 2\frac{\divergence^\genmet(e_\pm)}{d(d-1)} \mathrm{CS}(A,B) \\
            &= \mp5\: \frac{\divergence^\genmet(e_\pm)}{d(d-1)} \mathrm{CS}(A,B)  +3 Q(A,B)
        \end{split}
    \end{equation*}
    From this, we derive the quadratic components of the Ricci tensor as
    \begin{equation*}
        \begin{split}
            &2 \ric(B,B) = \mp5\: \frac{\divergence^\genmet(e_\pm)}{d} \absolute{B}^2 + \frac{3}{2(d-1)^2} \left\{(2d-3) \absolute{e_\pm}^2\absolute{B}^2 - (d-2) g(B,e_\pm)^2 \right\} 
        \end{split}
    \end{equation*}
    Therefore
    \begin{equation}\label{genriemflatnesseq6}
        \begin{split}
            2 \ric &= \mp5\: \frac{\divergence^\genmet(e_\pm)}{d} g + \frac{3}{2(d-1)^2} \left\{(2d-3) \absolute{e_\pm}^2g - (d-2) e_\pm \otimes e_\pm \right\}  \\
            &= \frac{3 \absolute{e_\pm}^2 \mp 5\: \divergence^\genmet(e_\pm)}{d} g + \frac{3(d-2)}{2(d-1)^2} \left\{\frac{ \absolute{e_\pm}^2}{d} g - e_\pm \otimes e_\pm \right\}  \\
        \end{split}
    \end{equation}
    For the Ricci scalar, we obtain
    \begin{equation*}
        2 \rscal = 3 \absolute{e_\pm}^2 \mp 5\: \divergence^\genmet(e_\pm) 
    \end{equation*}
    and for the trace-free part of the Ricci tensor 
    \begin{equation*}
        \begin{split}
            2 Z &\coloneqq 2\ric - \frac{2 \rscal}{d} \\
            &= \frac{3(d-2)}{2(d-1)^2} \left[  \frac{\absolute{e_\pm}^2}{d}g - e_\pm \otimes e_\pm \right]
        \end{split} 
    \end{equation*}
    We obtain, in the decomposition of the Riemann tensor $\riem  = S + E + W$ into scalar curvature part, trace-free Ricci part and Weyl part, assuming $A\perp B$
    \begin{equation*}
        \begin{split}
            S(A,B,A,B) \coloneqq \frac{\rscal}{d(d-1)} \absolute{A}^2\absolute{B}^2  = \frac{3\absolute{e_\pm}^2 \mp 5\: \divergence^\genmet(e_\pm)}{2d(d-1)}  \absolute{A}^2\absolute{B}^2 
        \end{split}
    \end{equation*}
    and, again assuming $A\perp B$
    \begin{equation*}
        \begin{split}
            &E(A,B,A,B) \\
            &\coloneqq \frac{1}{d-2} \left\{Z(A,A)\absolute{B}^2 + Z(B,B) \absolute{A}^2 \right\}  \\
            &= \frac{\absolute{B}^2}{2(d-2)} \left\{\frac{3(d-2)}{2(d-1)^2}\left[\frac{1}{d}\absolute{e_\pm}^2\absolute{A}^2 - g(A,e_\pm)^2 \right] \right\}\\
            &\quad + \frac{\absolute{A}^2}{2(d-2)} \left\{\frac{3(d-2)}{2(d-1)^2}\left[\frac{1}{d}\absolute{e_\pm}^2\absolute{B}^2 - g(B,e_\pm)^2 \right] \right\}\\
            & = \frac{3}{4(d-1)^2} \left\{\frac{2}{d}\absolute{e_\pm}^2 \absolute{A}^2 \absolute{B}^2 - \absolute{A}^2 g(B,e_\pm)^2 - \absolute{B}^2 g(A,e_\pm)^2 \right\}
        \end{split}
    \end{equation*}
    Therefore with $A\perp B$
    \begin{equation*}
        \begin{split}
            &2[E+S](A,B,A,B) \\
            & =  \frac{\mp5\: \divergence^\genmet(e_\pm)}{d(d-1)}  \mathrm{CS}(A,B)   + 3 Q(A,B) \\
            &= 2\: \riem(A,B,A,B)
        \end{split}
    \end{equation*}
    implying that the Weyl tensor vanishes! 
    
    Finally, we note that employing $\pi e_+ = \varepsilon \pi e_-$, the dilaton's equation of motion (\ref{dilatoneomeq}) becomes 
    \begin{equation*}
        0 = \frac{\absolute{H}^2}{6} - \extd^* \xi - \absolute{e_+}^2
    \end{equation*}
    In view of \eqref{epm:eq}, this coincides with the last equation in \eqref{genriemflatnessthmeq}: if $\varepsilon = +1$, then $\xi = 0$, but also $\divergence^\genmet(e_+) = 0$, since $\divergence^{\genmet} (e_+-e_-) = \divergence^\genmet(0) = 0$ and $\divergence^{\genmet} (e_++e_-)=0$ by \eqref{divpmv2:eq}.  And if $\varepsilon = -1$, then $\divergence^\genmet(e_+) = - \extd^* \xi$.
\end{proof}
The following corollary shows that, in the Riemannian case, we obtain complete triviality.
\begin{corollary} \label{riemannian_genriemflatnesscor}
     $\genriem^D = 0$ and $\dim M>2$ implies that $\absolute{H}^2_g = \absolute{e_\pm}^2_\genmet = 0$,  $\nabla^\pm \pi e_\pm = \nabla \pi e_\pm = 0$ and $\iota_{\pi e_\pm }H= 0$.
     In particular, if $\genmet$ is Riemannian, we obtain complete triviality, i.e.\  $\riem =0, H=0$, and $e= 0$.
\end{corollary}
\begin{remark}\label{genriemflatness_jet_rem2}
    Note that this Corollary holds under the assumptions described in Remark \ref{genriemflatness_jet_rem}.
\end{remark}
\begin{proof}
    Recall from \eqref{epm:eq} that $e = 2(X+\xi) \in \Gamma(T \oplus T^*)$ has $\xi = 0$ if $\varepsilon = +1$ and $X = 0$ if $\varepsilon = -1$. Consider the equation
    \begin{equation} \label{genriemflatnesseq10}
        g \nabla^\pm \pi e_\pm = \frac{\divergence^\genmet(e_\pm)}{d} g
    \end{equation}
    Its antisymmetric parts  
    \begin{equation*}
        \extd (g\pi e_\pm) \mp \iota_{\pi e_\pm }H = 0
    \end{equation*}
decouples in any case and yields, in particular, $\iota_{\pi e_\pm }H= 0$ and, hence, $\iota_{\pi e_\pm }H^{(2)}= 0$. It follows that 
    \begin{equation*}
        H^{(2)}(\pi e_+,B,\pi e_+,B)=0 \qquad \text{ for all $B \in \Gamma(TM)$ such that $g(e_+, B) = 0$,}
    \end{equation*}
    which is by (\ref{genriemflatnesseq5}) equivalent to
    \begin{equation}\label{genriemflatnesseq11} 
        \pm\frac{\divergence^\genmet(e_\pm)}{d(d-1)} = \frac{\absolute{e_\pm}^2}{2(d-1)^2} 
    \end{equation}
    If $\varepsilon = +1$, the left hand side is only independent of sign if it vanishes, and so $\absolute{e_\pm}^2 = \absolute{X}^2 = 0$. In that case, from the dilaton's equation of motion in (\ref{genriemflatnessthmeq}), we obtain $\absolute{H}^2 = 0$ and from (\ref{genriemflatnesseq10}) with $\iota_{\pi e_\pm }H= 0$ that $\nabla^\pm \pi e_\pm = \nabla \pi e_\pm = 0$.
    
    From now on, we can assume $\varepsilon = -1$, equivalently $e = 2\xi \in \Gamma(T^*M)$. Equations (\ref{genriemflatnesseq10}) and (\ref{genriemflatnesseq11}) become
    \begin{equation}\label{genriemflatnesseq12}
        \nabla \xi =- \frac{\extd^*\xi}{d} g, \qquad\qquad \frac{\absolute{\xi}^2}{2(d-1)} = -\frac{\extd^* \xi}{d}
    \end{equation}
    From the dilaton's equation of motion in (\ref{genriemflatnessthmeq}), we obtain that
    \begin{equation*}
        0 = \frac{\absolute{H}^2}{6} - \frac{d-2}{2(d-1)} \absolute{\xi}^2
    \end{equation*}
    and then with $\nabla H = 0$ also $\nabla \absolute{\xi}^2 = 0$ (since $d>2$). This, however, implies
    \begin{equation*}
        0 = \nabla_\xi \absolute{\xi}^2 = 2 g(\nabla_\xi \xi, \xi) = \frac{\absolute{\xi}^4}{2(d-1)}
    \end{equation*}
    Inserting this back into (\ref{genriemflatnesseq12}) and the dilaton's equation of motion, we respectively obtain $\nabla \xi = 0$ and $\absolute{H}^2 = 0$. With $\iota_{\pi e_\pm }H= 0$, we obtain that also in this case $\nabla^\pm \pi e_\pm = \nabla \pi e_\pm = 0$.
\end{proof}

\begin{corollary}\label{lorentzian_genriemflatness_cor}
    If $n+1 = \dim M > 2$ and $\genmet$ is Lorentzian, $\genriem^D = 0$ implies complete triviality.
\end{corollary}
\begin{proof}
    Let $p \in M$. We know that $\pi e_+|_p$ is a null vector. We can complete $\pi e_+|_p$ to a basis $\{e_0, ..., e_n\}, U = e_0, V = e_1,$ such that $\pi e_+|_p = V$, $\absolute{U}^2 = 0$, $g(U,V) = -2$ and for $i,j = 2,...,n$
    \begin{equation*}
        g(e_i, e_j) = \delta_{ij}, \qquad\quad g(e_i, U) = g(e_i, V) = 0.
    \end{equation*}
    It follows that
    \begin{equation}\label{lorentzian_genriemflatness_eq1}
        H^{(2)}(U, e_i, U, e_i) = \sum_{\mu,\nu} g^{\mu\nu} H(e_\mu, U, e_i)H(e_\nu, U, e_i) = \sum_{j=2}^n g^{jj} [H(e_j, U, e_i)]^2 \geq 0
    \end{equation}
    At the same time
    \begin{equation*}
        Q(U, e_i) = \frac{-1}{2(d-1)^2} \absolute{e_i}^2 g(U,V)^2 \leq 0
    \end{equation*}
    Noting that all possibly non-vanishing components of $Q$ are of the form $Q(U, e_i)$, we conclude  that $H^{(2)} (e_\mu, e_\nu, e_\mu, e_\nu)= 6Q(e_\mu, e_\nu) = 0$, cf.\  (\ref{genriemflatnesseq5}), \eqref{genriemflatnessthmeq}, and Corollary~\ref{riemannian_genriemflatnesscor}.  This implies $H= 0$: from (\ref{lorentzian_genriemflatness_eq1}), we already know that $H(U,e_i,e_j) = 0$ for all $i,j\in\{2,...,n\}$. Furthermore, from Corollary \ref{riemannian_genriemflatnesscor}, we also know that $\iota_VH = 0$. Now $\absolute{H}^2 = 0$ implies $H(e_i,e_j,e_k) = 0$ for all $i,j,k \in \{2,...n\}$.
    
    Finally, we can use e.g.\ (\ref{genriemflatnesseq4}) to conclude that $\riem = 0$. The formula for the Ricci tensor (\ref{genriemflatnesseq6}) then implies that $e_\pm = 0$.
\end{proof}

\begin{proposition}\label{genriemflatness_cfflatness_prop}
    Assume that $\genriem^D = 0$. Then, the factor of conformal flatness is given by
    \begin{equation}\label{genriemflatness_cfflatness_prop_eq}
        \grad_g \varphi = \frac{\sigma\sqrt{3}}{2(d-1)}\: \pi e_+
    \end{equation}
    where $\sigma \in \{+1,-1\}$. In other words, if $\varphi$ is a (local) solution of (\ref{genriemflatness_cfflatness_prop_eq}), then $\Tilde{g} = e^{2\varphi} g$ is flat.
\end{proposition}

\begin{proof}
    For $\varphi \in C^\infty(M)$ to be such that $\Tilde{g} = e^{2\varphi} g $ is flat, it has to solve
    \begin{equation}\label{genriemflatness_cfflatness_eq1}
        \ric = (d-2) \left(\nabla \extd \varphi - \extd \varphi \otimes \extd \varphi\right) + \left(\Delta \varphi + (d-2) \absolute{d\varphi}^2\right) g 
    \end{equation}
    We compare this to the Ricci tensor (\ref{genriemflatnesseq6})
    \begin{equation*}
        \ric = -\frac{3(d-2)}{4(d-1)^2} g\pi e_\pm \otimes g \pi e_\pm 
    \end{equation*}
    where we made use of Corollary \ref{riemannian_genriemflatnesscor}. Since $\absolute{e_\pm}^2= 0$ and $\nabla \pi e_\pm = 0$ by Corollary \ref{riemannian_genriemflatnesscor}, we see that indeed (\ref{genriemflatness_cfflatness_prop_eq}) provides a solution to (\ref{genriemflatness_cfflatness_eq1}).
\end{proof}

\begin{remark}
    That we get one solution $\varphi_\pm$ for the factor of conformal flatness in (\ref{genriemflatness_cfflatness_prop_eq}) for each of the two values $\sigma = \pm1$ implies 
    that we have a pair of conformally equivalent flat metrics $e^{2\varphi_\pm}g$.

    To see this more explicitly, we focus on $\varphi = \varphi_+$. We denote $\Tilde{g} = e^{2 \varphi} g$ (which is a flat metric) and $\Tilde{\nabla}$ for the Levi-Civita connection of $\Tilde{g}$. Then 
    \begin{equation}\label{genriemflatness_cffactor_minkowski_eq}
        0 = \nabla \extd \varphi = \Tilde{\nabla} \extd \varphi + 2 \extd\varphi \otimes \extd\varphi - \absolute{\extd \varphi}^2 g = \Tilde{\nabla} \extd \varphi + 2 \extd\varphi \otimes \extd\varphi
    \end{equation}
    implying that $e^{-2\varphi}g=  e^{-4\varphi}\Tilde{g}$ is again flat since its Ricci curvature vanishes:
    \begin{equation*}
        \ric[e^{-4 \varphi}\Tilde{g}] = -(d-2) \left( -2\Tilde{\nabla} \extd \varphi - 4\extd \varphi \otimes \extd \varphi\right) - \left(-2\Tilde{\Delta} \varphi + 4(d-2) \absolute{d\varphi}^2_{\Tilde{g}}\right) \Tilde{g} = 0
    \end{equation*}
    Herein, we used that $\Tilde{\Delta}\varphi  = -2 \absolute{\extd\varphi}^2_{\Tilde{g}} = 0$.
\end{remark}

\begin{example}\label{genriemflatness_example}
    Consider $M=\mathbb{R}^{2m}$, $d = 2m > 2$, endowed with the flat metric 
    \[ \Tilde{g}=\sum_{i=1}^m\left( (dx^i)^2 -(dy^i)^2\right) \] 
    of signature $(m,m)$ and its LC connection $\Tilde{\nabla}$, where  $(x^\mu)_{1\le\mu\le 2m}=(x^1,...,x^m,y^1,...,y^m)$ are standard coordinates. Replacing $(x^1,y^1)$ by $(u,v) = (x^1 + y^1, x^1-y^1)$, we obtain a coordinate system for which $\partial_u$ and $\partial_v$ are null. We consider only the region $u> 0$. Set $\varphi = \frac{\log(u)}{2}$. Then $\extd \varphi = \frac{\extd u}{2u} = \frac{e^{-2 \varphi}}{2}\extd u $, and hence (compare this to (\ref{genriemflatness_cffactor_minkowski_eq}))
    \begin{equation*}
        \Tilde{\nabla} \extd \varphi = -2 \extd \varphi \otimes \extd \varphi
    \end{equation*}
    It follows that $\nabla d\varphi = 0$, for the LC connection $\nabla$ of $g = e^{-2\varphi}\Tilde{g}$.  Hence the Ricci curvature of $g$ is given by
    \begin{equation*}
        \ric = -(d-2) \extd \varphi \otimes \extd \varphi 
    \end{equation*}
   
    To relate this to the dilaton, we set
    \begin{equation*}
        \pi e_+ = \frac{2(d-1)}{\sqrt{3}}  \grad_g \varphi = \frac{2(d-1)}{\sqrt{3}} \partial_v
    \end{equation*}
    where we used that $g(\partial_v) = d\varphi$. Note that we remain agnostic as to the value of $\varepsilon$ in $\pi e_+ = \varepsilon \: \pi e_-$. 
    We make the ansatz
    \begin{equation*}
        H = \sum_{i=2}^m f(u) \: g\pi e_+ \wedge \extd x^i \wedge \extd y^i,\quad f>0.
    \end{equation*}
    A small computation reveals that
    \begin{equation*}
        \begin{split}
            H^{(2)} &= \sum_{\mu,\nu} g^{\mu\nu} H(\partial_\mu) \otimes H(\partial_\nu) \\
            &= \sum_{i=2}^m f^2 e^{2\varphi} \left[(g\pi e_+ \wedge \extd y^i)^{\otimes 2} - (g\pi e_+ \wedge \extd x^i)^{\otimes 2} \right]
        \end{split}
    \end{equation*}
    Therefore, for $A, B \in \Gamma(TM)$
    \begin{equation*}
        \begin{split}
            &H^{(2)}(A,B,A,B) \\
            &= f^2 e^{2\varphi} \sum_{i=2}^m \left\{[g(\pi e_+, A) B_{y^i} -  g(\pi e_+, B) A_{y^i}]^2 - [g(\pi e_+, A) B_{x^i} -  g(\pi e_+, B) A_{x^i}]^2 \right\} \\
            &= f^2 e^{4\varphi} \left\{-g(\pi e_+, A)^2 \absolute{\rho B}^2_g - g(\pi e_+, B)^2 \absolute{\rho A}^2_g + 2 g(\pi e_+, A) g(\pi e_+, B) g(\rho A,\rho B) \right\}
        \end{split}
    \end{equation*}
    where $\rho$ denotes the orthogonal projection onto the distribution  
    $\{dx^1= d y^1 =0\}\subset TM$. Denoting $A_1 = A - \rho A$ and $B_1 = B - \rho B$, we compute $Q$ as (cf. (\ref{Qdefeq}))
    \begin{equation*}
        \begin{split}
            &2(d-1)^2 Q(A,B) \\
            &= - \absolute{B}^2 g(A,e_+)^2 - \absolute{A}^2 g(B,e_+)^2 + 2 g(A,B) g(A,e_+) g(B,e_+) \\
            &=  - \absolute{\rho B }^2 g(A,e_+)^2 - \absolute{\rho A}^2 g(B,e_+)^2 + 2 g(\rho A, \rho B) g(A,e_+) g(B,e_+)  \\
            &\quad- \absolute{B_1}^2 g(A,e_+)^2 - \absolute{A_1}^2 g(B,e_+)^2 + 2 g(A_1, B_1) g(A,e_+) g(B,e_+) \\
        \end{split}
    \end{equation*}
    Note that the last line vanishes, since up to a sign it is the norm-squared of 
    \begin{equation*}
        \begin{split}
            g(A_1,e_+)B_1 - g(B_1,e_+)A_1 &= g(\partial_u, e_+) \{A_u (B_u \partial_u + B_v \partial_v) - B_u(A_u \partial_u + A_v \partial_v) \} \\
            &=  g(\partial_u, e_+) (A_u B_v - A_v B_u) \partial_v 
        \end{split}
    \end{equation*}
    which is a null vector. We conclude that $H^{(2)}(A,B,A,B) = 6 Q(A,B)$ if and only if $f = \frac{\sqrt{3}}{d-1} e^{-2\varphi}$. Note that for $\mu\neq 1,m+1$
    \begin{equation*}
        \begin{split}
            \nabla (e^{-\varphi} \extd x^\mu) &= \Tilde{\nabla}(e^{-\varphi} \extd x^\mu) + e^{-\varphi} \extd \varphi \otimes \extd x^\mu + e^{-\varphi} \extd x^\mu \otimes \extd \varphi - e^{-\varphi} \Tilde{g}(\extd \varphi, \extd x^\mu) \Tilde{g} \\
            &= e^{-\varphi} \extd x^\mu \otimes \extd \varphi = \frac{\sqrt{3}}{2}(d-1) e^{-\varphi} \extd x^\mu \otimes g\pi e_+
        \end{split}
    \end{equation*}
    so that
    \begin{equation*}
        \begin{split}
            \nabla H &= \sum_{i=2}^m\frac{\sqrt{3}}{d-1} \nabla( g \pi e_+ \wedge e^{-\varphi} \extd x^i \wedge e^{-\varphi} \extd y^i) \\
            &= \sum_{i=2}^m\frac{3}{2} e^{-2\varphi} \left[\extd x^i \otimes g \pi e_+ \wedge g\pi e_+ \wedge \extd y^i + \extd y^i \otimes g \pi e_+ \wedge \extd x^i \wedge g\pi e_+  \right] \\
            & = 0
        \end{split}
    \end{equation*}
    Because $g$ is conformally flat by construction and $H^{(2)}(A,B,A,B) = 6 Q(A,B)$, we find that generalised Riemann flatness is equivalent to (cf.\ also (\ref{genriemflatnesseq3}, \ref{genriemflatnesseq4}))
    \begin{equation*}
        \ric = \frac{1}{4} H^2 = -\frac{f^2(d-2)}{4} e^{4\varphi} g \pi e_+ \otimes e \pi e_+ = -\frac{3(d-2)}{4(d-1)^2} g \pi e_+ \otimes g \pi e_+
    \end{equation*}
    which is true by construction.
\end{example}

\section{The Fundamental Theorem for Generalised Hypersurfaces}\label{fundthmsection}
    
     The fundamental theorem for hypersurfaces is a classical result \cite[Chapter VII, Theorem 7.1]{KN2} establishing that every tuple $(\Sigma, h, k)$ consisting of a Riemannian manifold $(\Sigma, h)$ and a symmetric two-tensor $k$ satisfying the flat Gauß and Codazzi equations is locally equivalent to a hypersurface in flat space equipped with the induced metric and exterior curvature tensor. In this section, we establish the analogous result for the generalised case.
    
    The idea in our proof is to make use of the result that generalised Riemann flatness implies complete triviality (Corollary \ref{riemannian_genriemflatnesscor}) to reduce to the classical situation, where we already have the theorem. The problem is that we do not have access to the generalised Riemann tensor of an ambient space to which we could apply Corollary~\ref{riemannian_genriemflatnesscor}. Our solution to this  is inspired by the classical case.

    One approach in the classical case (cf., though in  slightly different settings, the proofs of \cite[Theorem 8.1]{nomizu_sasaki} or \cite[Theorem 2.3]{cortes2014completeness}) is to define a synthetical version $\tilde{T}$ of the tangent space to the ambient manifold by defining the \enquote{normal bundle} $\nu = \Sigma \times \reals$ with generator $n$ and setting $\tilde{T} = T\Sigma \oplus \nu$. One then defines a connection $\tilde{\nabla}$ on $\tilde{T}$ by setting $\tilde{\nabla}_X Y = \nabla^\Sigma_X Y - k(X,Y)n$ and $\tilde{\nabla}_X n = AX$, where $k = h A$. One finds that $\tilde{\nabla}$ is flat and hence that $(\tilde T, \tilde\nabla )$ can be locally identified with $M \times \mathbb{R}^{d}$, where $d-1= \dim \Sigma$.  After some further steps,\footnote{We actually go through the needed steps in our proof of Theorem \ref{fth} as we do not refer to but reprove the classical result.} one obtains from this the local immersion into $\mathbb{R}^{d}$.    Note that $\tilde{\nabla}$ is not a $\tilde{T}$-connection, i.e.\ $n$ is not allowed as a direction entry, owing to the fact that one cannot know the normal derivative of an object only defined on $\Sigma$. 
    
    Transferring the approach from the classical to the generalised case, one synthetically defines the generalised normal bundle $\mathcal{N} = \langle n, n^\flat \rangle$ and then the generalised tangent bundle $\tilde{E}_\Sigma = E_\Sigma \oplus \mathcal{N}$. Now, one might naively expect that one should define an $E_\Sigma$-connection $\tilde{D}$ on $\tilde{E}_\Sigma$, in the hope that it is flat so that one can somehow apply Corollary \ref{riemannian_genriemflatnesscor} to its curvature tensor, reduce to the classical case and proceed as there. However, the key insight is that such a connection would not encode sufficient information.  It turns out that the needed extra input are the covariant derivatives in the conormal direction $n^\flat$. In contrast to those in the normal direction, these can be meaningfully defined, as their computation does not involve taking an actual derivative. In particular, one can compute $\tilde{D}_{n^\flat} T$ for a generalised tensor $T$ only defined over $\Sigma$. Thus, introducing $L = \langle n^\flat \rangle$, $\tilde{D}$ should be defined as an $E_\Sigma \oplus L$-connection on $\tilde{E}_\Sigma$. Indeed, we find in Lemma \ref{gft_extendedconn_lemma} that one can associate to such a connection a well-defined torsion tensor $T^{\tilde{D}} \in \Gamma(\Lambda^3(E_\Sigma \oplus L)^*)$, as well as a generalised Riemann tensor $\genriem^{\tilde{D}}$ defined on the bundles (\ref{genriem_mixedtypesubbundle_def} - \ref{quad_L:eq}). It is to this curvature tensor that we will apply Corollary  \ref{riemannian_genriemflatnesscor} and obtain the fundamental theorem for generalised hypersurfaces.
    
    We begin now by describing the general assumptions of this section and presenting some preparatory constructions.
    \paragraph{General assumptions.}
    Let $(E_\Sigma, \genmetHS, \divergence_\Sigma, D^\Sigma, \genext, \mathcal L)$ be a tuple consisting of an exact CA $E_\Sigma \to \Sigma$ over a $(d-1)$-dimensional manifold $\Sigma$, equipped with a generalised Riemannian metric $\genmetHS$, a divergence operator $\divergence_\Sigma$, the canonical generalised LC connection $D^\Sigma$ on $(E_\Sigma , \genmetHS)$ such that (cf.\ Corollary \ref{cangenLC_divergence_submf_lemma})
    \begin{equation}\label{fth_div_assummption_eq}
        \divergence_{D^\Sigma} = \divergence_\Sigma + \frac{1}{d-1}\scalbrack{e_\Sigma,\cdot},
    \end{equation}
    two sections $\mathcal L^\pm \in \Gamma ((E_\Sigma^\pm)^*)$, and two bilinear forms $\genext^{n_\pm} \in \Gamma((E_\Sigma)^* \otimes (E_\Sigma^\pm)^*)$. The latter we ask to decompose as in (\ref{genextcomputationlemmaeq}) with $\chi_\pm^\perp = 0$, i.e.
    \begin{equation}
        \begin{split}
            \restr{\genext^{n_\pm}}{E_\Sigma^\pm \times E_\Sigma^\pm} &= k - \frac{e^\perp_\pm}{2 (d-1)} h \mp \frac{H^\perp}{6} \\
            \restr{\genext^{n_\pm}}{E_\Sigma^\mp \times E_\Sigma^\pm} &= k \mp \frac{H^\perp}{2} \\
        \end{split}
    \end{equation}
    for the Riemannian metric $h$  associated with $\mathcal H$, a symmetric two-tensor $k$, two functions $e_\pm^\perp \in C^\infty(\Sigma)$ and a two-form $H^\perp$. 
    \begin{remark}
        Note that the decomposition of $\genext^{n_\pm}$ into $k$, $e_\pm^\perp$ and $H^\perp$ is unique.
    \end{remark}

    \paragraph{Extending the Courant algebroid and the generalised connection.}
     Define the trivial vector bundles $\nu \coloneqq \Sigma \times \reals$  and $\mathcal{N} \coloneqq \Sigma \times \reals^2$ and denote their trivialising frames by $n$ and $\{n_+,n_-\}$. Define then $\tilde{T} \coloneqq T\Sigma \oplus \nu$ and $\Tilde{E}_\Sigma \coloneqq E_\Sigma \oplus \mathcal{N}$. We extend the inner product onto $\Tilde{E}_\Sigma$ by requiring it to make the direct sum orthogonal and setting $\scalbrack{n_\pm,n_\pm} 
    =\pm 1$, $\langle n_+,n_-\rangle =0$. Similarly, we define a second   symmetric tensor $\tilde{\genmet}$ on $\Tilde{E}_\Sigma$ by the conditions that the subbundles  $\tilde{E}_{\Sigma}^\pm \coloneqq E_{\Sigma}^\pm \oplus \mathbb{R} n_\pm$ are $\tilde{\genmet}$-orthogonal and that $\tilde{\genmet}= \pm \langle \cdot ,\cdot \rangle >0$ on those subbundles. We define a positive definite symmetric tensor $\tilde{g}$ on $\tilde{T}$ such that it agrees with $h$ on $T\Sigma$ and makes $n$ a unit normal to that subspace.
    
    Now we define an $E_\Sigma$-connection $\Tilde{D}$ on $\Tilde{E}_\Sigma$ by
    \begin{equation*}
        \Tilde{D}_a b_\pm \coloneqq D^\Sigma_a b_\pm - \genext^{n_\pm}(a,b_\pm) n_\pm, \quad \Tilde{D}_a n_\pm \coloneqq \genextEnd^{n_\pm}(a),\quad a\in E_\Sigma,\; b_\pm \in \Gamma (E_\Sigma^\pm ),
    \end{equation*}
    where $\genextEnd^{n_\pm}:={\mathcal H}^{-1}\genext^{n_\pm}$. Then one can check that $\Tilde{D}\tilde{\mathcal G}= \Tilde{D}\langle \cdot , \cdot \rangle = 0$. Next we 
    consider the line bundle $L\to \Sigma$ spanned by $n_+-n_-$ and     define an anchor map $\tilde{\pi}$ for $E_\Sigma\oplus L\to \Sigma$ as the composition of the natural projection  $E_\Sigma \oplus L\to E_\Sigma$ with the anchor of $E_\Sigma$. We also extend the Dorfman bracket 
    of $E_\Sigma$ trivially to $\Gamma (E_\Sigma \oplus L)$ such that $[\Gamma (L), \Gamma (L) ]=0$ and $[v, f(n_+-n_-) ]=-[f(n_+-n_-),v ]=\pi (v)(f)(n_+-n_-)$ 
    for all $v\in \Gamma (E_\Sigma)$, $f\in C^\infty (\Sigma )$. 
    We note that $(E_\Sigma \oplus L, \scalprodmap|_{E_\Sigma \oplus L},\brackmap )$ satisfies the axioms of a Courant algebroid with exception of the non-degeneracy of the scalar product. 
    \begin{lemma}\label{gft_extendedconn_lemma}
        There is a unique extension of the $E_\Sigma$-connection $\tilde D$ to an $(E_\Sigma \oplus L , \tilde{\pi})$-connection on $\tilde E_\Sigma = E_\Sigma \oplus \mathcal N$ (which we denote by the same symbol $\tilde D$) which is still metric for $\tilde{\mathcal{G}}$ and $\langle \cdot , \cdot \rangle$ and satisfies the following formula (inspired by Lemma~\ref{nflatlemma}~(ii) and the definition of the conormal exterior curvature $\mathcal L^\pm$ in the case of an actual hypersurface) 
        \begin{equation}\label{genfundhypersurfaceeq1} 
            \begin{split}
                \Tilde{D}_{n_\pm - n_\mp} n_\pm &= \gennormalshape^\pm := \mathcal H^{-1}\mathcal L^\pm,\\
                \genmet(\Tilde{D}_{n_\pm-n_\mp} a_\pm, b_\pm) &= \genext^{n_\pm}(a_\pm,b_\pm) - \genext^{n_\pm}(b_\pm,a_\pm),\\
                \genmet( \Tilde{D}_{n_\pm-n_\mp} a_\pm, n_\pm) &= -\gennormalext^\pm (a_\pm ).
            \end{split}
        \end{equation}
    \end{lemma}
    \begin{proof}
        To check that $\tilde D \scalprodmap =0$ one can easily evaluate the covariant derivative $\langle \tilde D_u v, w\rangle + \langle  v,\tilde D_u w\rangle -\pi (u) \langle v,w\rangle$ on elements $u\in E_\Sigma \oplus L$, 
        $v,w\in \Gamma (\tilde E_\Sigma )$, where it is sufficient to consider sections of $E_\Sigma^\pm$ and the sections $n_+-n_-$ and $n_\pm$
        of $L$ and $\mathcal N$, respectively. Since the decomposition $\tilde E_\Sigma = \tilde E_\Sigma^+ \oplus \tilde E_\Sigma^-$ is invariant under $\tilde D$, we see that also $\tilde D \mathcal G=0$. 
    \end{proof}
    We now construct an auxiliary ambient exact CA $E_M \to M$ restricting to  $\tilde{E}_\Sigma$ over $\Sigma$  and reproducing the structure on it. On $E_M$ we construct \textit{an} extension $D^M$ of $\tilde{D}$ to a 
    (generally non-flat) canonical Levi-Civita connection realising the given exterior curvature tensors.

    \begin{lemma}\label{tildeD_extension_lemma}
        There exists an exact Courant algebroid $\pi_M\colon E_M \to M$ over a $d$-dimensional manifold $M$ containing $\Sigma$ as a submanifold and on $E_M$ a generalised metric $\genmet_M$, a divergence operator $\divergence_M$, and finally a vector bundle isomorphism $\Phi \colon \tilde{E}_\Sigma \cong E_M|_\Sigma$ such that 
        \begin{enumerate}[label = (\arabic*)]
            \item $\Phi^*\genmet_M = \tilde{\genmet}$ and $\Phi^*\scalprodmap_M = \scalprodmap$,
            \item $\Phi(n)$ becomes, under the $\genmet_M$-induced identification $E_M \cong \mathbb{T}M$, the unit normal $\mathbf{n} \in \Gamma(TM)$, and
            \item denoting by $H_M$ the preferred representative of the \v Severa class and writing $\divergence_M = \divergence^{\genmet_M} - \scalbrack{e_M, \cdot}_M$,
            \begin{equation*}
                H_M|_\Sigma = H_\Sigma + \mathbf{n}^\flat \wedge H^\perp, \qquad e_M|_\Sigma = \Phi(e_\Sigma) + e_+^\perp \mathbf{n}_+ - e_-^\perp \mathbf{n}_-.
            \end{equation*}
            where $\mathbf{n}_\pm = \mathbf{n} \pm \mathbf{n}^\flat$.
        \end{enumerate}
        In particular, with $D^M$ the canonical generalised Levi Civita connection on $E_M$ with divergence $\divergence_M$, it holds for all $a,b,c \in \Gamma(E_\Sigma \oplus L)$ and all $\beta \in \Gamma(\tilde{E}_\Sigma)$
        \begin{equation}\label{DM_tildeD_eq}
            {\scalbrack{[a^*,b^*]_{M}, c^*}_M|}_\Sigma = \scalbrack{[a,b],c}, \qquad\quad D^M_{a^*} \beta^*|_\Sigma = \Phi (\tilde{D}_{a} \beta).
        \end{equation}
        Herein, $a^*,b^*,c^*$, and $\beta^*$ denote arbitrary extensions of respectively $\Phi(a),\Phi(b), \Phi(c)$, and $\Phi(\beta)$ to sections of $E_M$. 
    \end{lemma}
    \begin{proof}
        Consider $M = \Sigma \times \reals$. Identify $E_\Sigma$ with $\mathbb{T}\Sigma$  via the isomorphism induced by the generalised metric $\genmetHS$, denote the corresponding closed three-form by $H_\Sigma\in \Omega^3(\Sigma)$. Denote $\mathbf{n} = \partial_t$ and $\mathbf{n}^\flat = \extd t$, and then $\mathbf{n}_\pm = \mathbf{n} \pm \mathbf{n}^\flat$. Arbitrarily extend $H = H_\Sigma + \mathbf{n}^\flat \wedge H^\perp$ to a closed three-form on $M$ and $e = e_\Sigma + e_+^\perp \mathbf{n}_+ - e_-^\perp \mathbf{n}_-$ to a section of $\mathbb{T}M$. Define a generalised metric $\genmet_M$ on $M$ such that over $\Sigma$, it is determined by $\genmetHS$ and $\mathbf{n}$ as indicated in Proposition \ref{inducedsemiriemCAprop}. We consider now $E_M = \mathbb{T}M$ as the generalised tangent bundle with twist $H$, standard scalar product $\langle \cdot ,\cdot \rangle_M = \langle \cdot ,\cdot \rangle$,  generalised metric $\genmet_M$, and divergence operator $\divergence_M = \divergence^\genmet - \scalbrack{e,\cdot}$. Define $\Phi\colon \tilde{E}_\Sigma \to E_M|_\Sigma$ on $E_\Sigma = \mathbb{T}\Sigma$ as the canonical identification of $\mathbb{T}\Sigma$ with the orthogonal subbundle $\langle \mathbf{n}_+, \mathbf{n}_- \rangle^\perp \subset E_M$, and on $\mathcal{N}$ by demanding $n_\pm \mapsto \mathbf{n}_\pm$. 
        In this way, $\Phi (E_\Sigma)$ is identified with the semi-Riemannian Courant algebroid over $\Sigma$ induced from $E_M$.

        Now consider the canonical generalised Levi-Civita connection $D^M$. By the assumption (\ref{fth_div_assummption_eq}) made on $\divergence_{D^\Sigma}$ and Corollary \ref{cangenLC_divergence_submf_lemma}, $D^\Sigma$ is the $\Phi$-pullback of the  generalised Levi-Civita connection that $D^M$ induces on $\Phi(\mathbb{T}\Sigma) \subset \mathbb{T}M$. The statement in (\ref{DM_tildeD_eq}) regarding the generalised connections is then due to Lemmas \ref{finalgausslemma} and \ref{nflatlemma}.  Regarding the statement for the bracket, note that $L$ is the radical of $(E_\Sigma\oplus L, \scalprodmap)$ and
    \[\pi_M([a^*,b^*]_M)|_\Sigma = \mathcal{L}_{\pi_M a^*}\pi_M b^*|_\Sigma = \mathcal{L}_{\pi a} \pi b \]
        is tangent to the submanifold  $\Sigma \subset M$ ($\mathcal{L}_XY$ stands for the Lie derivative of vector fields) and thus  $[a^*,b^*]_M|_\Sigma$ is a section of $\Phi (E_\Sigma \oplus L)$. Hence, 
        for $a,b,c \in \Gamma(E_\Sigma)$, the statement for the brackets reduces to Proposition \ref{inducedCAprop}. For more general $a,b,$ and $c$, it follows from the observation that the bracket of a pair of sections of $E_M$ extending a section of $L$ and a section of $E_\Sigma \oplus L$ satisfies on $\Sigma$ the defining equations for the bracket on $E_\Sigma \oplus L$ up to terms in $L$. Using the projection \[ \rho \colon E_M|_\Sigma = \Phi( \tilde{E}_\Sigma)= \Phi (E_\Sigma\oplus L \oplus \langle n\rangle )\to \Phi (E_\Sigma \oplus \langle n\rangle )\] we can write these relations as follows. 
        For all $\xi,\eta \in \Gamma(E_M)$ such that $\xi|_\Sigma = f \mathbf{n}^\flat$ and $\eta|_\Sigma = g \mathbf{n}^\flat$ for some functions $f,g\in C^\infty(\Sigma)$ and all $a \in \Gamma(E_M)$ such that $a|_\Sigma \in \Gamma(\Phi(E_\Sigma))$ the following equations hold on $\Sigma$:
        \begin{equation*}
            [\xi, \eta]_M=0 \qquad \text{and} \qquad \rho[a, \xi ]_M=-\rho[\xi,a ]_M  =0  
        \end{equation*}
    \end{proof}

    \begin{corollary}
        We can define a torsion-tensor $T^{\tilde{D}}\in \Gamma (\Lambda^3 (E_\Sigma \oplus L)^*)$ for $\tilde D$ via the formula \eqref{torsion_def:eq} and this tensor vanishes. 
        Via the formula (\ref{genriemdefeq}), we can associate to $\tilde{D}$ the generalised Riemann tensor $\genriem^{\tilde{D}}$, which is a section of the direct sum of the pure-type subbundle
        \begin{equation}\label{genriem_puretypesubbundle_def}
            \Sym^2\,\Lambda^2\,\big(E_\Sigma^+\big)^*  \oplus \Sym^2\,\Lambda^2 \big(E_\Sigma^-\big)^*, 
        \end{equation}
        the mixed-type subbundle 
        \begin{equation}\label{genriem_mixedtypesubbundle_def}
            \left(\Lambda^2(\tilde{E}_\Sigma^+)^* \oplus \Lambda^2 (\tilde{E}_\Sigma^-)^*\right) \vee \left((E_\Sigma^+)^* \wedge (E_\Sigma^-)^*\right)
        \end{equation}
        and the conormal subbundles
        \begin{eqnarray}
             && (L^*\wedge (E_\Sigma^+)^* \oplus L^*\wedge (E_\Sigma^-)^*)\vee (\Lambda^2(E_\Sigma^+)^* \oplus \Lambda^2(E_\Sigma^-)^*)\quad\mbox{and}\label{lin_L:eq}\\
             && 
             \Sym^2(L^*\wedge (E_\Sigma^+)^*\oplus L^*\wedge (E_\Sigma^-)^*),\label{quad_L:eq}
        \end{eqnarray}
        which are respectively linear and quadratic in $L^*$, where $L:=\langle n_+ - n_-\rangle$.
    \end{corollary}
    \begin{remark}
        Technically, on the mixed-type subbundle, the formula (\ref{genriemdefeq}) cannot be employed, since in the direction entry of $\tilde{D}$ there could feature contributions proportional to $n = (n_++n_-)/2$, which would be ill-defined. For instance, $n_+ = \frac12 (n_+-n_-) + n$ contains such a contribution.  However, we can deal with this by formally assuming that also the derivative $\tilde{D}_n$ respects the decomposition $\tilde{E}_\Sigma = \tilde{E}_\Sigma^+ \oplus \tilde{E}_\Sigma^-$, as then one obtains immediately that the terms featuring such a problematic derivative vanish. Alternatively, if one wishes to avoid such formal evaluations, one can adapt the definition of $\genriem^{\tilde{D}}$ on the mixed-type subbundle to simply exclude these undesired terms, which yields a definition more in line with that of the generalised Riemann tensors $\genriem_{\mathrm{GF}}^\pm$ found in \cite{genricciflow}.
    \end{remark}
    \begin{proof}
        Employing the generalised connection $D^M$ from Lemma~\ref{tildeD_extension_lemma}, which extends  $\tilde{D}$, the formulas for $T^{\tilde{D}}$ and $\genriem^{\tilde{D}}$ respectively equate by (\ref{DM_tildeD_eq}) to restrictions of the formulas for the torsion $T^{D^M}$ and the curvature $\genriem^{D^M}$ of $D^M$, which are well-defined tensors. 
    \end{proof}

\paragraph{The Generalised Fundamental Theorem for Hypersurfaces.} Assume now that the geometric data $(E_\Sigma, \genmetHS, \divergence_\Sigma, D^\Sigma, \genext, \mathcal L)$,  specified 
under the general assumptions of section \ref {fundthmsection},  satisfy the equations of Gauß and Codazzi for a flat ambient Riemannian CA, i.e.\ the equations from Theorems~\ref{riemgaußthm} and \ref{codazzithm} with an ambient space of vanishing generalised Riemann curvature. Then, we obtain the announced  theorem.

\begin{theorem}\label{fth} 
    Under the above assumptions, for every simply connected open set $U \subset \Sigma$ there exists a Riemannian (hypersurface) immersion $U \to \reals^{d}$ such that the canonical structure on the untwisted generalised tangent bundle $\mathbb{T}\reals^{d}$ induces $(\genmetHS,\divergence_\Sigma, D^\Sigma, \genext, \mathcal L)$. 
    
    In particular, $e=0$,  $H = 0$,   $\mathcal L=0$, and $(\Sigma ,h,k)$ satisfies the well-known flat Gauß and Codazzi equations from Riemannian geometry.
\end{theorem}
\begin{remark}
    The Riemannian immersion $U \to \mathbb{R}^{d}$ is unique up to isometries of $\mathbb{R}^{d}$. This follows from the fundamental theorem for hypersurfaces, which includes this uniqueness statement in the version \cite[Chapter VII, Theorem 7.2]{KN2} and applies since $(\Sigma,h,k)$ satisfies the classical Gauß and Codazzi equations.
\end{remark}
\begin{proof}
   Denote by $\tilde{D}$ the ($E_\Sigma \oplus L)$-connection from Lemma \ref{gft_extendedconn_lemma}, and by $\genriem^{\tilde{D}}$ its generalised Riemann curvature. Note that the components of the generalised Riemann tensor involved in the Gauß and Codazzi equations (cf. Theorems \ref{riemgaußthm} and \ref{codazzithm}) precisely determine the tensor on the subbundles (\ref{genriem_puretypesubbundle_def}), (\ref{genriem_mixedtypesubbundle_def}),  \eqref{lin_L:eq} and \eqref{quad_L:eq}.

    Next we relate the ($E_\Sigma \oplus L)$-connection $\tilde{D}$ to a triple $(\tilde \nabla , H, e)$ consisting of:
    \begin{enumerate}
    \item the metric connection $\tilde{\nabla}$ on $\tilde{T}$ obtained from setting $\Tilde{\nabla}_X Y = \nabla^\Sigma_X Y - k(X,Y) n$, $\Tilde{\nabla}_Xn=AX$, where $k=hA$ and $\nabla^\Sigma$ is the Levi-Civita connection of $\Sigma$,
    \item $H = H_\Sigma + n^\flat \wedge H^\perp\in \Gamma (\Lambda^3 \tilde T^*)$, and 
    \item $e = e_\Sigma + e_+^\perp\: n_+ - e_-^\perp\: n_- \in \Gamma (\tilde E_\Sigma )$.
    \end{enumerate}
    To that end, we introduce the related connections $\tilde{\nabla}^{\pm}$ and $\tilde{\nabla}^{\pm1/3}$ on $\tilde{T}$ by setting, cf.\ (\ref{nablapmdefeq}),
    \begin{equation*}
        \tilde{\nabla}^\pm_X = \tilde{\nabla}_X  \pm \frac{1}{2} H_X, \qquad\quad \tilde{\nabla}^{\pm1/3}_X = \tilde{\nabla}_X  \pm \frac{1}{6} H_X
    \end{equation*}
    \begin{proposition}
   Employing the obvious identifications $\tilde{E}_\Sigma^\pm \cong \tilde{T}$ mapping $n_\pm$ to $n$ and restricting to the anchor on $E_\Sigma^\pm$, the relationship between $\tilde D$ and the triple $(\tilde \nabla , H, e)$ is given by, compare Lemma \ref{cangenconnlemma},
\begin{equation}\label{genfundhypersurfaceeq2}
        \begin{split}
            \Tilde{D}_a \beta &= \Tilde{\nabla}^{\pm1/3}_a \beta + \frac{1}{d-1}\chi_\pm^{e_\pm}(a,\beta), \qquad\quad \Tilde{D}_{\bar{a}} \beta = \Tilde{\nabla}^{\pm}_{\bar{a}} \beta, \\
            \Tilde{D}_{n_\pm - n_\mp} \beta &= \mp\frac{1}{3}H^\perp(\beta)+ \frac{1}{d-1}\chi_\pm^{e_\pm}(n_\pm,\beta)
        \end{split}
    \end{equation}
    where $\chi_\pm^{e_\pm}$ is formally defined as in Lemma \ref{cangenconnlemma}, $a \in \Gamma(E_\Sigma^\pm), \bar{a} \in \Gamma(E_\Sigma^\mp)$, and $\beta \in \Gamma(\tilde{E}_\Sigma^\pm)$. In particular, the restrictions of 
    $\tilde D$ to pure-type entries $\Gamma (E_\Sigma^\pm) \times \Gamma (\tilde E_\Sigma^\pm )$
    and mixed-type entries $\Gamma (E_\Sigma^\pm) \times \Gamma (\tilde E_\Sigma^\mp )$ respectively define 
    two connections in $\tilde T$. 
    \end{proposition}
    \begin{proof}
        This follows from applying Lemma \ref{cangenconnlemma} to the extension of $\tilde{D}$ to the canonical generalised LC connection $D^M$ from Lemma \ref{tildeD_extension_lemma}.
    \end{proof}
    Since $\tilde{\nabla}$ is a metric connection in $\tilde T \to \Sigma$, it  has a curvature tensor $\riem^{\tilde{\nabla}}$ ,  
    defined via the classical formula \begin{equation}\label{genfundhypersurfaceeq3}
        \riem^{\Tilde{\nabla}}(\tilde{V},\tilde{W},X,Y) = \tilde{g}\left(\tilde{V}, \tilde{\nabla}_{X} \tilde{\nabla}_{Y} \tilde{W} - \tilde{\nabla}_{Y} \tilde{\nabla}_{X} \tilde{W} - \tilde{\nabla}_{[X, Y]} \tilde{W} \right), 
    \end{equation}
    where $ X, Y \in \Gamma(T\Sigma)$ and $\tilde{V}, \tilde{W} \in \Gamma(\tilde{T})$. For better comparison with generalised curvatures, we symmetrically extend it to a section 
    $\riem^{\Tilde{\nabla}}$ of 
    \begin{equation*}
        (\Lambda^2 \tilde{T}^*) \vee (\Lambda^2 T^*\Sigma).
    \end{equation*}
    Note that $\riem^{\tilde{\nabla}}$ is related to $\riem^{\nabla^\Sigma}$ by the classical Gauß and Codazzi equations. From these, together with the Bianchi identity for $\riem^{\nabla^\Sigma}$ and the symmetry of $k$ it follows that $\riem^{\tilde{\nabla}}$ satisfies the Bianchi identities 
    \begin{equation}\label{riemnablatilde_bianchi_eq}
        \begin{split}
            0 &= \sum_{\sigma(X,Y,Z)}\riem^{\tilde{\nabla}}(X,Y,Z,V) \\
            0 &= \sum_{\sigma(X,Y,Z)}\riem^{\tilde{\nabla}}(X,Y,Z,n) = \sum_{\sigma(n,X,Y)} \riem^{\tilde{\nabla}}(Z,n,X,Y) \\
        \end{split}
    \end{equation}
    where $X,Y,Z,V \in \Gamma(T\Sigma)$. The second equality in the second line follows from the symmetries of $\riem^{\tilde{\nabla}}$.
    \begin{lemma} \label{gft_genriem_expression_lemma}
        The generalised Riemann tensor $\genriem^{\tilde{D}}$ and the Riemann tensor $\riem^{\tilde{\nabla}}$ are related via the formulas given in Propositions \ref{puretyperiemprop} and \ref{mixedtyperiemprop}, where $\riem^{\Tilde{\nabla}}$ takes the place of $\riem$, $\genriem^{\tilde{D}}$ takes the place of $\genriem^D$, $H = H_\Sigma + n^\flat  \wedge H^\perp$ and $e = e_\Sigma + e_+^\perp n_+ - e_-^\perp n_-$. 
    \end{lemma}
    \begin{remark}
        We briefly clarify which formulas Lemma \ref{gft_genriem_expression_lemma} claims to hold for the conormal components of the generalised Riemann tensor, given that these components are absent from Propositions \ref{puretyperiemprop} and \ref{mixedtyperiemprop}, which only feature formulas for pure- and mixed-type components. 

        To start with, we note that the curvature components specified in Propositions \ref{puretyperiemprop} and \ref{mixedtyperiemprop} fully determine any algebraic generalised Riemann tensor, as defined below equation \eqref{genriemdefeq}, cf.\ Remark \ref{genriem_domain_rem}. The assertion in Lemma \ref{gft_genriem_expression_lemma} is that, formally employing the multi-linearity of the tensor and then applying the replacement prescription detailed in its statement to the extracted equations, we obtain formulas that hold in our setting. In particular, Lemma \ref{gft_genriem_expression_lemma} claims that for every component of $\genriem^{\tilde{D}}$ that is well-defined, we obtain in this way an equation for it whose right-hand side is also well-defined.
    \end{remark}
    \begin{proof}
        Consider the extension $D^M$ of $\tilde{D}$ from Lemma \ref{tildeD_extension_lemma} and their relationship (\ref{DM_tildeD_eq}). As $D^M$ is a canonical generalised LC connection, Propositions \ref{puretyperiemprop} and \ref{mixedtyperiemprop} apply to it. Since $E_M$ is equipped with structure inducing $H$ and $e$ on $\Sigma$, and since $\genriem^{D^M} = \genriem^{\tilde{D}}$ on the domain of the latter, the result follows.
    \end{proof}
    Let us consider $\Sigma \subset M=\Sigma \times \reals$ via the canonical embedding $p\mapsto (p,0)$. We denote $n = \partial_t$, $n_\pm = \partial_t \pm \extd t$, and define at every point $p \in \Sigma$ the following data. First, we extend $\tilde{\nabla}|_p$ to a connection on $M$ at $p$ by setting $\tilde{\nabla}_n X|_p = \partial_tX|_p +   A(X)|_p$ for any vector field $X$ tangent to $\Sigma$ and $\tilde{\nabla}_n (fn)|_p =(\partial_tf) n|_p$, $f\in C^\infty (M)$. We then extend $\riem^{\tilde{\nabla}}|_p$ to a  tensor\footnote{In fact, from the Bianchi identities (\ref{riemnablatilde_bianchi_eq}) for $\riem^{\tilde{\nabla}}$, it follows that $\riem$ satisfies the Bianchi identity as well and hence is (even a priori) an algebraic curvature tensor.} $\riem \in \Sym^2(\Lambda^2T_p^*M)$ by demanding $\riem(X,n,n,X) = 0$. Finally, we define $g \in J^1_p(M,\Sym^2(M))$, $H \in  J^1_p(M,\Lambda^3M)$, and $e = X+\xi \in J^1_p(TM \oplus T^*M)$ by requiring the following:
    \begin{equation*}
        \begin{split}
            &g_p = \extd t^2_p + h_p, \qquad H_p = H_\Sigma|_p + \extd t_p \wedge H^\perp_p, \qquad  e_p = [e_\Sigma + e^\perp_+ n_+ - e^\perp_- n_-]_p\\
            &\tilde{\nabla}g|_p = 0, \qquad\qquad\qquad \tilde{\nabla}_n H|_p = 0, \qquad\qquad\qquad\quad \tilde{\nabla}_n e|_p = 0, \\
           &\tilde{\nabla}_X H_p = \left[\tilde{\nabla}_X (H_\Sigma + \extd t \wedge  H^\perp)\right]_p, \qquad \tilde{\nabla}_X e|_p = \left[\tilde{\nabla}_X (e_\Sigma + e^\perp_+ n_+ - e^\perp_- n_-)\right]_p.
        \end{split}
    \end{equation*}
    These definitions imply the equations one obtains from Propositions \ref{puretyperiemprop} and \ref{mixedtyperiemprop} with the generalised Riemann tensor assumed to vanish - the non-trivial part of this statement being that it applies to the equations that come from the assumed vanishing of the components
    \begin{equation*}
        \genriem^{\Tilde{D}}(a, n_\mp, v,w), \qquad\quad \genriem^{\Tilde{D}}(a, n_\pm, n_\pm,w), \qquad\quad \genriem^{\Tilde{D}}(n_\mp, n_\pm, n_\pm,w)
    \end{equation*}
    which lie in neither of the bundles (\ref{genriem_puretypesubbundle_def}), (\ref{genriem_mixedtypesubbundle_def}), \eqref{lin_L:eq} and \eqref{quad_L:eq}. Thus, by Remarks \ref{genriemflatness_jet_rem} and \ref{genriemflatness_jet_rem2}, we can apply Corollary \ref{riemannian_genriemflatnesscor} to obtain that $H_p = 0, e_p = 0$, and $\riem_p = 0$ for all $p \in \Sigma$. In particular, the inital connection $\tilde{\nabla}$ in $\tilde T \to \Sigma$ is flat, and we can find a parallel local ON frame $\{e_i\}$ for $\tilde{T}$ on an open set $p \in U \subset \Sigma$. This defines the trivialising isomorphism 
    \begin{equation*}
        \phi \colon (\tilde{T}, \tilde{\nabla}) \longrightarrow (\underline{\reals}^{d}, \nabla)
    \end{equation*}   
    where $\tilde{\nabla} = \phi^* \nabla$. Then for all $X,Y \in \Gamma(T\Sigma)$, 
    \begin{equation*}
        \begin{split}
            \nabla_X \phi(Y) &= \phi(\Tilde{\nabla}_{X} Y) = \phi(\nabla^\Sigma_{X} Y - k(X, Y)n)
        \end{split}
    \end{equation*}
    so that the one-form $\alpha := \phi|_{T\Sigma} \in \Omega^1(\Sigma, \underline{\reals}^d)$ satisfies 
    \begin{equation*}
        \begin{split}
            \extd \alpha(X,Y) &= \nabla_X \alpha(Y) - \nabla_Y \alpha(X) - \alpha([X,Y]) \\
            &= \alpha\left(\nabla^\Sigma_{X} Y - k(X, Y) n\right) - \alpha\left(\nabla^\Sigma_{Y} X - k(Y, X) n\right) - \alpha([X,Y]) \\
            &= \alpha(\nabla^\Sigma_X Y - \nabla^\Sigma_Y X - [X,Y]) \\
            &= 0
        \end{split}
    \end{equation*}
    By assumption, $U$ is simply connected. Hence, we can find $\varphi \colon U \to \reals^{d}$ such that $\alpha = \extd \varphi$. Note that $\varphi$ is an immersion because $\alpha$ has, restricted to the fibre at any point, maximal rank. After possibly shrinking $U$, $\varphi$ becomes an embedding. Employing the splitting $F_\Sigma \colon E_\Sigma \cong \mathbb{T}\Sigma$ obtained from $\genmetHS$, we find with
    \begin{equation*}
        \Bar{\varphi} \coloneqq \begin{pmatrix}
            \varphi_* & 0 \\
            0 & (\varphi^*)^{-1} 
        \end{pmatrix}
    \end{equation*}
    that
    \begin{equation*}
        F\coloneqq \bar{\varphi} \circ \restr{F_\Sigma}{U} \colon E_U \longrightarrow \mathbb{T}\reals^{d}
    \end{equation*}
    is as desired. 
\end{proof}

\cleardoublepage
\section{Appendix A: The Generalised Riemann Tensor for the Canonical Connection}\label{genriemappendix}
Let $E\to M$ be an exact CA with semi-Riemannian generalised metric $\genmet$ and divergence operator $\divergence = \divergence^\genmet - \scalbrack{e, \cdot}$. Denote by $D$ the canonical generalised Levi Civita connection with divergence $\divergence$. We cite results proven in \cite{vicentethomas} that provide a decomposition of the generalised Riemann tensor in terms of the classical Riemann tensor $\riem$, the twist $H$, and the dilaton $e$.
\begin{proposition}\label{puretyperiemprop}
    The pure-type components of the generalised Riemann curvature are given by
    \begin{equation*}
        \begin{split}
            &\pm\genriem^{D}(a,b,v,w) \\
            &= \riem(a,b,v,w)  - \frac{1}{36} H^{(2)}(a,v,b,w) - \frac{1}{36}H^{(2)}(b,v,w,a) - \frac{1}{18} H^{(2)}(v,w,a,b) \\
            &\quad \pm \frac{1}{2(d-1)} \left\{ [D^0_v \chi_\pm^{e_\pm}](w,b,a) - [D^0_w \chi_\pm^{e_\pm}](v,b,a)\right. \\
            &\qquad\qquad\qquad \left. + [D^0_b \chi_\pm^{e_\pm}](a,v,w) - [D^0_a \chi_\pm^{e_\pm}](b,v,w) \right\} \\
            &\quad + \frac{1}{2(d-1)^2} \left\{2 \genmet(e_\pm,e_\pm) \big[\genmet(w,a)\genmet(v,b)-\genmet(v,a)\genmet(w,b)\big] \right. \\
            &\qquad\qquad\qquad + \genmet(a,e_\pm) \big[\genmet(w,b)\genmet(v,e_\pm)-\genmet(v,b)\genmet(w,e_\pm)\big] \\
            &\qquad\qquad\qquad + \left.\genmet(b,e_\pm) \big[\genmet(v,a)\genmet(w,e_\pm)-\genmet(w,a)\genmet(v,e_\pm)\big]\right\}
        \end{split}
    \end{equation*}
    Herein, $a,b,v,w \in \Gamma(E_\pm)$ arbitrary, we denoted $d \coloneqq \dim M$, and from now on we use the notation $\riem(a,b,v,w) \coloneqq \riem (\pi a,\pi b, \pi v, \pi w)$. 
\end{proposition}
\begin{corollary}
    The generalised scalar curvature is given by
    \begin{equation}\label{genscaleq}
        \genscal(\genmet, \divergence) = \rscal - \frac{\absolute{H}^2}{12} + \divergence^\genmet(e_+ - e_-) - \frac{1}{2}\absolute{e}_\genmet^2
    \end{equation}
    Herein, $\divergence = \divergence^\genmet - \scalbrack{e, \cdot}$.
\end{corollary}
\begin{proposition}\label{mixedtyperiemprop}
    The mixed-type components of the generalised Riemann curvature are given by
    \begin{equation*}
        \begin{split}
            &\pm2 \genriem^D(a,\Bar{b},v,w) \\
            &= \riem(a,\Bar{b},v,w) \mp \frac{1}{2}[\nabla_{a} H](\Bar{b},v,w) \pm \frac{1}{6}[\nabla_{\Bar{b}} H](a,v,w) \\
            &\quad - \frac{1}{12} H^{(2)}(\Bar{b},w,a,v) - \frac{1}{12} H^{(2)}(w,a,\Bar{b},v) - \frac{1}{6} H^{(2)}(a,\Bar{b},v,w) \\
            &\quad \pm \frac{1}{d-1}  [D^0_{\Bar{b}} \chi_\pm^{e_\pm}](a,v,w)            
        \end{split}
    \end{equation*}
    Herein, $a,v,w \in \Gamma(E_\pm)$ and $\Bar{b}\in \Gamma(E_\mp)$ arbitrary, and we denoted $d \coloneqq \dim M$.
\end{proposition}
\begin{corollary}\label{genriccor}
    The mixed-type components of the generalised Ricci curvature are given by
    \begin{equation}\label{genriceq}
        4\restr{\genric(\genmet,\divergence)}{E_\mp \times E_\pm}= 4\ric - H^{2} \mp 2\extd^* H +4 [\nabla\xi]^{\mathrm{sym}} \pm 4[\nabla g X]^{\mathrm{antisym}} \pm 2 H(\xi) 
    \end{equation}
    Herein, $\divergence = \divergence^\genmet -2 \scalbrack{X+\xi, \cdot}$.
\end{corollary}

\section{Appendix B: The Divergence of Mixed-Type Generalised Tensors}\label{mixedtypedivappendix}
In this section, we investigate the divergence of generalised tensors. The main result
is that the divergence of a mixed-type tensor with respect to a generalised Levi-Civita connection $D$ is 
invariant in the sense that it depends only on the divergence operator of $D$ and the generalised metric. 
This is relevant for the proper interpretation of the generalised Codazzi equations in the context of the generalised momentum constraint, see Corollary \ref{genmomentumconstrcor}. Indeed the constraints should only involve the 
fields of the considered gravitational theory, which involve the generalised metric and the divergence operator
but no further components of the generalised connection. 
\begin{definition}\label{divgentensordef}
    Let $E$ be an exact Courant algebroid, let $D \in \mathcal{D}(E)$ be a generalised connection on $E$. Then, we define the following divergence operator on the space of generalised tensors:
    \begin{equation}\label{conndivgentensordefeq}
        \divergence_D\colon \Gamma(E^{r+1}_s) \longrightarrow \Gamma(E^r_s); \qquad\quad \divergence_D(T) \coloneqq \tr(D T)
    \end{equation}
    or, in index notation,
    \begin{equation*}
        \divergence_D(T)^{B_1...B_r}_{C_1...C_s} = D_A T^{AB_1...B_r}_{C_1...C_s}
    \end{equation*}
    Let now $\genmet$ be a generalised metric on $E$. Then, we can define the canonical divergence operators
    \begin{equation}\label{metricdivgentensordefeq}
        \begin{split}
            \divergence^{\genmet,\pm} &\colon \Gamma(E_\mp \otimes (E_\pm)^r_s) \longrightarrow \Gamma((E_\pm)^r_s), \\
            \divergence^{\genmet,\pm}&(T) = \sigma_\pm \circ\tr(\nabla^\pm [\pi T \sigma_\pm]) \circ \pi,
        \end{split}
    \end{equation}
    where $\nabla^\pm$ are the Bismut connections from (\ref{nablapmdefeq}). In index notation, this reads
    \begin{equation*}
        \divergence^{\genmet,\pm}(T)^{B_1...B_r}_{C_1...C_s} = (\sigma_\pm)^{B_1}_{b_1} ... (\sigma_\pm)^{B_r}_{b_r} \nabla^\pm_a [\pi T \sigma_\pm]^{ab_1...b_r}_{c_1...c_s} (\pi \pi_\pm)^{c_1}_{C_1}...(\pi\pi_\pm)^{c_s}_{C_s}
    \end{equation*}
    Note that this formula defines tensors over $E$ as opposed to $E_\pm$, hence the appearance of the projection $\pi_\pm$ on the right hand side.
\end{definition}
The next Lemma investigates the meaning and compatibility of these definitions. Crucially, it asserts that the divergence operator $\divergence_D$ is not uniquely determined by the requirement that $D$ be generalised LC, i.e. $D\in \mathcal{D}^0(\genmet, \divergence)$. However, restricted to the right subspace of generalised tensors, this property is achieved.
\begin{lemma}\label{divgentensorlemma}
    Let $E$ be an exact Courant algebroid equipped with a generalised metric $\genmet \cong (g,F)$ and a compatible divergence operator $\divergence = \divergence^\genmet - \scalbrack{e, \cdot}$. Then, the divergence operator $\divergence_D$ from (\ref{conndivgentensordefeq}) depends non-trivially on the choice of divergence compatible generalised connection $D \in \mathcal{D}(\divergence)$. 
    
    However, if $D$ is assumed to be metric compatible and torsion-free, then the action of $\divergence_D$ on the \enquote{mixed-type} sub-spaces $\Gamma(E_\mp \otimes (E^*_\pm)^s)$ is independent of the choice of generalised connection. It then holds
    \begin{equation*}
        \divergence(T) = \divergence^{\genmet,\pm}(T) - \scalbrack{e, T}, \qquad\quad T \in \Gamma(E_\mp \otimes (E^*_\pm)^s)
    \end{equation*}
    where $\divergence^{\genmet,\pm}$ is as in (\ref{metricdivgentensordefeq}).
\end{lemma}
\begin{remark}
    Note that an analogous result applies to generalised tensors of arbitrary contravariant rank. This is due to compatibility of generalised connections with the inner product $\scalprodmap$, which provides the isomorphism $E \cong E^*$ used to identify these spaces.
\end{remark}
\begin{proof}
    We calculate, taking an orthonormal frame $\{e_A\}$ of $E$ and setting $\epsilon_A \coloneqq \scalbrack{e_A, e_A}$,
    \begin{equation*}
        \begin{split}
            &(\divergence_D\: T)(e_{A_1},...,e_{A_s}) \\
            & = \sum_B \epsilon_B \scalbrack{e_B, (D_{e_B} T)(e_{A_1},..., e_{A_s})} \\
            &= \sum_B \epsilon_B \left\{ \scalbrack{e_B, D_{e_B} [T(e_{A_1},...,e_{A_s})]} \right.\\
            &\qquad\left. - \scalbrack{e_B, T(D_{e_B} e_{A_1},...,e_{A_s}) + ... +  T(e_{A_1},...,D_{e_B}e_{A_s})}\right\} \\
            &= \divergence(T(e_{A_1},..., e_{A_s}))\\
            &\qquad -\sum_B \epsilon_B \scalbrack{e_B, T(D_{e_B} e_{A_1},...,e_{A_s}) + ... +  T(e_{A_1},...,D_{e_B}e_{A_s})} \\
            &= \divergence^\genmet(T(e_{A_1},..., e_{A_s})) - \scalbrack{e,T(e_{A_1},..., e_{A_s})}\\
            &\qquad -\sum_B \epsilon_B \scalbrack{e_B, T(D_{e_B} e_{A_1},...,e_{A_s}) + ... +  T(e_{A_1},...,D_{e_B}e_{A_s})} \\
        \end{split}
    \end{equation*}
    Note that, without further assumptions, this expression is in general not independent of the choice of divergence compatible generalised connection $D$. However, assuming $D \in \mathcal{D}^0(\genmet, \divergence)$ and $T \in \Gamma(E_+ \otimes (E^*_-)^s)$ (the case of flipped signs is analogous), the expression is uniquely defined, since the mixed-type operators $D^\pm_\mp$ are, cf.\ Lemma~\ref{genlcspacelemma}. Using $D_{\sigma_+ X} \sigma_- Y = \sigma_-(\nabla^-_X Y) $, cf.\ Lemma~\ref{cangenconnlemma}, we see that
    \begin{equation*}
        \begin{split}
            &(\divergence_D\: T)(e^-_{a_1},...,e^-_{a_s}) \\
            &= \divergence^\genmet(T(e^-_{a_1},..., e^-_{a_s})) - \scalbrack{e,T(e^-_{a_1},..., e^-_{a_s})}\\
            &\qquad -\sum_b \scalbrack{e^+_b, T(D_{e^+_b} e^-_{a_1},...,e^-_{a_s}) + ... +  T(e^-_{a_1},...,D_{e_b^+}e^-_{a_s})} \\
            &= \sum_b \Big\{g\left(e_b, \nabla_{e_b} [\pi T\sigma_-(e_{a_1},...,e_{a_s})]\right) \\
            &\qquad - g\left( e_b, \pi T \sigma_-(\nabla^-_{e_b} e_{a_1},...,e_{a_s}) + ... +  \pi T\sigma_-(e_{a_1},...,\nabla^-_{e_b}e_{a_s})\right)\Big\} \\
            &\qquad - \scalbrack{e,T(e_{A_1},..., e_{A_s})} \\
            &= \sum_b \Big\{g\left(e_b, \nabla_{e_b} \pi T \sigma_-)(e_{a_1},..., e_{a_s})\right) \\
            &\qquad + g\left(e_b,(\pi T \sigma_-)((\nabla- \nabla^-)_{e_b} e_{a_1},...,e_{a_s}) +... +  T(e_{a_1},...,(\nabla- \nabla^-)_{e_b}e_{a_s})\right)\Big\} \\
            &\qquad - \scalbrack{e,T(e^-_{a_1},..., e^-_{a_s})} \\
            &= \sum_b \Big\{g\left(e_b, \nabla_{e_b} \pi T \sigma_-)(e_{a_1},..., e_{a_s})\right) \\
            &\qquad + g\left(e_b,(\pi T \sigma_-)(H(e_b, e_{a_1})^\sharp,...,e_{a_s}) +... +  T(e_{a_1},...,H(e_b, e_{a_s})^\sharp)\right)\Big\} \\
            &\qquad - \scalbrack{e,T(e^-_{a_1},..., e^-_{a_s})} \\
            &= \sum_b g\left(\nabla^-_{e_b}  \pi T \sigma_-)(e_{a_1},..., e_{a_s}),e_b\right) - \scalbrack{e,T(e^-_{a_1},..., e^-_{a_s})} \\
            &= (\divergence^{\genmet,-} T)(e^-_{a_1},..., e^-_{a_s}) - \scalbrack{e,T}(e^-_{a_1},..., e^-_{a_s}) \\
        \end{split}
    \end{equation*}
    This proves the claim.
\end{proof}
\bibliographystyle{unsrturl}  
\bibliography{literatur}

V.\ Cort\'es: vicente.cortes@uni-hamburg.de

O.\ Schiller: 
oskar.schiller@uni-hamburg.de

Department of Mathematics, University of Hamburg,  Bundesstr.\ 55, 20146  Hamburg, Germany.
\cleardoublepage
\end{document}